\documentclass[11pt]{article}

\usepackage[utf8]{inputenc}
\usepackage[T1]{fontenc}
\usepackage{lmodern}
\usepackage[a4paper,margin=1in]{geometry}
\usepackage{microtype}
\usepackage{mathtools,amsmath,amssymb,amsthm,bm}
\usepackage{graphicx}
\usepackage{xcolor}
\usepackage{booktabs}
\usepackage{bbm}
\usepackage{enumitem}
\usepackage{hyperref}
\usepackage[nameinlink,noabbrev,capitalise]{cleveref}
\usepackage[linesnumbered,ruled,vlined]{algorithm2e}
\usepackage{caption}
\usepackage{subcaption}
\usepackage{url}
\hypersetup{colorlinks=true,linkcolor=darkblue,citecolor=darkblue,urlcolor=darkblue}
\usepackage[most]{tcolorbox}
\tcbset{colback=white, colframe=black!15, boxrule=0.4pt, rounded corners,
        left=6pt,right=6pt,top=6pt,bottom=6pt}

\usepackage{booktabs,tabularx,ragged2e,array,adjustbox,makecell}
\newcolumntype{Y}{>{\RaggedRight\arraybackslash}X}

\usepackage[most]{tcolorbox}
\tcbset{colback=white, colframe=black!15, boxrule=0.4pt,
        rounded corners, left=6pt,right=6pt,top=6pt,bottom=6pt,
        before skip=2pt, after skip=2pt}

\usepackage{booktabs,tabularx,ragged2e,array,adjustbox}

\newcommand{\St}{\mathrm{St}}
\DeclareMathOperator{\grad}{grad}

\definecolor{darkblue}{rgb}{0.05,0.05,0.6}

\theoremstyle{plain}
\newtheorem{theorem}{Theorem}[section]
\newtheorem{proposition}[theorem]{Proposition}
\newtheorem{lemma}[theorem]{Lemma}

\newtheorem*{proposition*}{Proposition}
\newtheorem*{theorem*}{Theorem}
\newtheorem*{corollary*}{Corollary}
\newtheorem*{remark*}{Remark}
\theoremstyle{definition}

\newtheorem{remark}[theorem]{Remark}
\newtheorem{assumption}[theorem]{Assumption}

\DeclareMathOperator{\diag}{diag}
\DeclareMathOperator{\argmin}{argmin}

\newcommand{\KL}{\mathrm{KL}}
\newcommand{\softmax}{\mathrm{softmax}}

\newcommand{\gravidy}{\textsc{Gravidy}}
\newcommand{\gravidypos}{\textsc{Gravidy-Pos}}
\newcommand{\gravidydelta}{\textsc{Gravidy-$\Delta$}}
\newcommand{\gravidybox}{\textsc{Gravidy-Box}}
\newcommand{\gravidyst}{\textsc{Gravidy-St}}

\title{Constrained Optimization via Constraint-Induced Geometry\\
\textit{Implicit Feasible Dynamics and Optimality from Stationarity}}
\author{Valentin Leplat}
\date{\today}

\begin{document}
\maketitle

\begin{abstract}
We introduce \gravidy, a geometry-aware framework for constrained optimization in which the constraints are encoded directly into the optimization dynamics. The resulting motion remains feasible by construction. The precise geometric mechanism depends on the constraint set: componentwise reparameterizations and their induced Hessian geometries for the nonnegative orthant and box constraints, Fisher--Shahshahani and KL geometry for the simplex, and canonical Riemannian geometry for the Stiefel manifold. Within this framework, we derive feasible continuous-time flows and implicit discretizations adapted to each geometry. On the vector domains, the implicit updates admit exact Bregman-proximal interpretations. This yields monotone descent and convergence guarantees for convex objectives, a linear contraction under relative strong convexity, and a Kurdyka--Lojasiewicz convergence analysis for nonconvex problems under the stated compact-interiority, decrease, and relative-error assumptions. We also show that convergent trajectories generated from the interior recover the Karush--Kuhn--Tucker conditions on the orthant, simplex, and box. The same conclusion holds for convergent implicit sequences generated from the interior when their stepsizes are bounded away from zero. On the Stiefel manifold, stationarity of the canonical Riemannian gradient is equivalent to the usual first-order optimality condition. The resulting algorithms combine large implicit outer steps with problem-adapted Newton, modified Gauss--Newton, Newton--KKT, and Newton--Krylov inner solvers. Numerical experiments on nonnegative, simplex-constrained, box-constrained, and orthogonality-constrained problems show rapid convergence to high accuracy in a small number of outer iterations while preserving feasibility of the accepted outer iterates. A sparse elastic-obstacle experiment further illustrates how the orthant construction can exploit large structured systems directly.
\end{abstract}


\section*{Dedication}
In memory of Jean-Marie Souriau (1922--2012), whose insistence that geometry should guide dynamics inspired this work.

\section{Introduction}

We consider smooth constrained optimization problems of the form
\[
\min_{x\in Q} \Phi(x),
\]
where \(Q\subset\mathbb{R}^n\) is a structured feasible set, such as the nonnegative orthant, a box, the probability simplex, or a manifold such as the Stiefel manifold. In many classical methods, feasibility is enforced externally, for instance by projection, by barrier terms, or by ad hoc feasible updates. These approaches are often effective, but they do not always expose the geometry induced by the constraints themselves. Our goal in this paper is to start from the opposite direction: rather than correcting the dynamics after the fact, we build the constraints directly into the geometry of the flow.

The central idea is to associate with the feasible set a \emph{constraint-induced geometry} that shapes the optimization dynamics. In the orthant and box settings, this geometry is generated by a componentwise reparameterization \(x=g(u)\), with Jacobian factors that encode how the boundary is approached. On the simplex, the corresponding dynamics are governed by the Fisher--Shahshahani tangent geometry~\cite{Shahshahani1979}, which preserves mass and acts only along feasible directions. On the Stiefel manifold, feasibility is expressed intrinsically through the canonical Riemannian structure. These mechanisms are not identical, and we do not force them into a single rigid metric model. What remains common, however, is the optimization principle: the constraints are translated into the feasible motion itself, and the resulting gradient-type dynamics evolve inside the admissible set by construction.

This viewpoint leads naturally to feasible continuous-time dynamics and, more importantly for computation, to implicit discrete schemes. The use of implicit discretization is not a secondary numerical detail. It is a structural consequence of the flow interpretation. Backward-Euler-type steps provide robust large-step updates on vector domains, while feasible Cayley-type updates play the analogous role on the Stiefel manifold. In this way, the continuous and discrete levels remain aligned: the geometry shapes the flow, and the discretization respects that geometry.

A second key feature of the framework is the relation between limiting
stationarity and constrained first-order optimality. On the vector domains,
the dynamics are generated from interior points. The accompanying parameter
or dual variables therefore retain the boundary information that may be lost
in a purely primal zero-velocity condition. We show that convergent
trajectories generated from the interior recover the KKT conditions on the
orthant, simplex, and box. The same conclusion holds for convergent implicit
sequences generated from the interior when their stepsizes are bounded away
from zero. Complementarity and the boundary inequalities
emerge from the behavior of the geometry-generating variables as the boundary
is approached. On the Stiefel manifold, the vanishing of the canonical
Riemannian gradient is equivalent to the usual equality-constrained
first-order condition. Stationarity is therefore not only a dynamical
endpoint: together with the constraint-induced geometry, it recovers the
optimality system of the original problem.

We refer to the resulting family of methods as \gravidy. The name originates from the vanishing-Jacobian mechanism that appears in the orthant and box settings, while the broader framework also includes the simplex and Stiefel cases through their own native geometries. The unification is therefore conceptual rather than formulaic. The paper is unified at the level of optimization principle, not at the level of a single universal metric tensor.

\paragraph{Contributions.}
The main contributions of the paper are the following.
\begin{itemize}[itemsep=2pt,topsep=3pt]

\item 
We introduce a geometry-driven framework in which the feasible dynamics are
built directly from the constraint structure. The framework covers orthant
and box constraints through componentwise reparameterizations and their
induced positive operators, the simplex through Fisher--Shahshahani tangent
geometry, and the Stiefel manifold through its canonical Riemannian geometry.

\item 
We establish the connection between the limiting behavior of the feasible
dynamics and the corresponding first-order optimality conditions. For the
orthant, simplex, and box, convergent trajectories generated from the interior
satisfy the KKT conditions at their limits. The same conclusion holds for
convergent implicit sequences generated from the interior when their
stepsizes are bounded away from zero. The boundary signs and complementarity
relations are recovered from the parameter or dual dynamics. On the Stiefel
manifold, the vanishing of the canonical Riemannian gradient is equivalent to
the standard constrained first-order condition.

\item 
We derive implicit discretizations that preserve feasibility and remain
faithful to the underlying geometry. The vector-domain steps admit
Bregman-proximal formulations, including a KL-proximal step on the simplex.
For the Stiefel manifold, we use an implicit Cayley update that preserves
orthogonality. The resulting nonlinear systems are solved using damped Newton
or modified Gauss--Newton methods, Newton--KKT solvers on the simplex, and
Newton--Krylov or dense Newton methods on the Stiefel manifold.

\item 
For the exact vector-domain steps, we establish monotone descent and
last-iterate convergence guarantees in the convex case, together with a
linear contraction under relative strong convexity. We also give a
Kurdyka--\L{}ojasiewicz analysis for the nonconvex setting under the stated
compact-interiority, decrease, and relative-error assumptions. Numerical
experiments on four constraint geometries show that the implicit methods
produce large reductions of the objective and stationarity measures in a
small number of outer iterations. The additional elastic-obstacle experiment
shows that the orthant construction can also exploit a large sparse structured
problem directly, without forming a dense Hessian or a matrix square root.
\end{itemize}

\subsection{Related work and positioning}
\label{sec:related}

\paragraph{Continuous-time and dynamical-systems viewpoints.}
The interpretation of optimization algorithms through differential equations has a long history, from early ODE-based approaches to minimization \cite{boggs1977algorithm,brown1989some,zghier1981use} to modern analyses of first-order methods, accelerated dynamics, and variational formulations \cite{su2014differential,wibisono2016variational,scieur2017integration,shi2021understanding,luo2022differential,suh2022continuous}. Dynamical viewpoints have also been developed for proximal-gradient and operator-splitting methods \cite{attouch1996dynamical,attouch2014dynamical,franca2018admm,hassan2021proximal}. Our work belongs to this line of research, but focuses on feasible flows whose geometry is induced directly by the constraints.

\paragraph{Implicit discretization, proximal steps, and regularized Newton viewpoints.}
The backward Euler step is closely related to proximal regularization and to
the proximal point method \cite{rockafellar1976monotone}. This connection is
particularly useful here because it allows the vector-domain outer steps to
be analyzed through descent and resolvent arguments, while the inner
nonlinear systems can be solved by damped Newton or
Levenberg--Marquardt-type procedures
\cite{higham1999trust,attouch2011continuous,
doi:10.1080/10556788.2020.1712602}. The stability motivation comes from the
classical A-stability of backward Euler
\cite{dahlquist1963special}. In the convex
vector-domain setting, this is reflected by convergence results that do not
require an upper restriction on the outer stepsize. The Stiefel construction
has a different, structure-preserving interpretation: the Cayley transform
maps a skew-symmetric generator to an orthogonal update, and the resulting
scheme coincides with the trapezoidal rule when the skew field is frozen.
Implicitness is therefore both numerical and geometric, but its precise role
depends on the constraint class.

\paragraph{Manifold optimization and orthogonality constraints.}
Optimization on smooth manifolds is now a mature area, with well-developed geometric algorithms and theory \cite{AbsilMahonySepulchre2008,Boumal2023}. For orthogonality constraints, the Stiefel manifold provides the canonical setting, and the geometry of tangent spaces, retractions, and transports has been studied in depth, including in the foundational treatment of orthogonality-constrained algorithms by Edelman, Arias, and Smith \cite{EDELMAN1998}. Feasible Cayley-type updates are also classical in practice, notably in the method of Wen and Yin \cite{WenYin2013}. Our contribution does not compete with this literature at the level of general manifold formalism. Rather, it places the Stiefel case inside the same optimization philosophy as the vector-domain problems considered in this paper: feasibility is built into the dynamics, and the discrete update is derived from an implicit flow viewpoint rather than introduced only as a retraction mechanism.


\paragraph{Bregman and information-geometric viewpoints.}
Hessian--Riemannian gradient flows induced by Legendre functions have been
studied in convex programming, including their continuous-time convergence,
boundary behavior, dual trajectories, and relation to Bregman distances
\cite{AlvarezBolteBrahic2004,attouch2004singular}. Once the potential \(h\)
is defined by
\[
\nabla h(x)=g^{-1}(x),
\]
the primal orthant and box flows considered here fit within this general
class.

Our starting point and main purpose are nevertheless different. We begin with
a concrete componentwise feasibility map \(x=g(u)\), derive the associated
metric and mobility from this map, and discretize the resulting dynamics
implicitly in the free variable \(u\). Each implicit step is then treated as
a nonlinear computational problem, for which we develop Newton and
regularized Gauss--Newton inner solvers. We also study how convergent implicit
sequences whose stepsizes are bounded away from zero recover the KKT
conditions, and develop corresponding constructions for the simplex and
Stiefel manifold.

Bregman and information-geometric ideas have also played an important role in
constrained and non-Euclidean optimization, especially when the natural
geometry is not Euclidean. In applications such as nonnegative factorization
and inverse problems, KL and \(\beta\)-divergence geometries are particularly
relevant \cite{Fevotte_betadiv,VL2021}. Information-geometric and
probability-space flows have also been used to study accelerated dynamics and
optimization over distributions
\cite{WangWuchen2020,Taghvaei2019,pmlr_v97_liu19i}.

\paragraph{Algorithms for nonnegative and bound-constrained least squares.}
For nonnegative least squares and related problems, classical active-set, block-pivoting, and projection-based methods remain important baselines \cite{lawson1995solving,10.2307/2153286,bro1997fast}. Bound-constrained smooth optimization has also benefited from highly effective quasi-Newton schemes, such as limited-memory methods \cite{doi:10.1137/0916069}. Our goal is not to replace this mature algorithmic literature by a single universal solver. Instead, we propose a complementary viewpoint: solve the unconstrained dynamics in geometry-generating variables, recover feasibility by construction, and use implicit steps to obtain robust updates for the resulting constrained motion.

\paragraph{Positioning of the present work.}
The present paper differs from prior work in three main ways. First, it develops a single optimization principle that covers positivity, box, simplex, and Stiefel constraints, while remaining honest about the fact that these cases arise from different geometric realizations. Second, it links feasible dynamics to first-order optimality by showing how the geometry itself encodes complementarity, supportwise gradient equalization, or manifold stationarity. Third, it derives practical implicit algorithms from this continuous-time structure and analyzes them with optimization tools, rather than treating the geometric update and the convergence analysis as separate layers.

The remainder of the paper is organized as follows. Section~\ref{sec:preliminaries} introduces the
notation and preliminary material used throughout the paper. Section~\ref{sec:riemannian-lens}
presents the geometric viewpoint underlying the different constraint classes.
Section~\ref{sec:gravidy} introduces the \gravidy{} framework and derives the corresponding
feasible flows for the orthant, simplex, box, and Stiefel settings. Section~\ref{sec:kkt-from-flow}
studies how limiting stationarity recovers the appropriate first-order
optimality conditions. Section~\ref{sec:algorithms} derives the implicit algorithms and discusses
the associated inner solvers. Section~\ref{sec:global-convergence} establishes the convergence results,
while Section~\ref{sec:experiments} presents the numerical experiments. Finally, Section~\ref{sec:conclusion}
summarizes the main conclusions and discusses several directions for future
work, including adaptive implicit schemes and a possible extension of the
Lie-algebra viewpoint suggested by the Stiefel construction.


\section{Preliminaries: Notation, Assumptions, Flows, and Implicit Stability}
\label{sec:preliminaries}

This section introduces the notation and the continuous-time models used throughout the paper. We also recall the linear notion of A-stability and explain what it gives, and what it does not give, for the implicit schemes considered later.

\subsection{Notation and Standing Assumptions}
\label{sec:notation}

\paragraph{Notation.}
We work in finite-dimensional Euclidean spaces. For vectors $x,y\in\mathbb{R}^n$, we write
\[
\langle x,y\rangle=x^\top y,
\qquad
\|x\|=\|x\|_2.
\]
For a matrix $M$, $\|M\|_2$ denotes its spectral norm and $\operatorname{spec}(M)$ its spectrum. We use
\[
\mathbb{R}_+^n:=\{x\in\mathbb{R}^n:x\geq 0\},
\qquad
\Delta_n:=\{x\in\mathbb{R}_+^n:\mathbf{1}^\top x=1\},
\]
and, for vectors $\ell,b\in\mathbb{R}^n$ with $\ell_i<b_i$,
\[
[\ell,b]:=\{x\in\mathbb{R}^n:\ell\leq x\leq b\}.
\]
The Stiefel manifold is
\[
\St(n,p):=\{X\in\mathbb{R}^{n\times p}:X^\top X=I_p\}.
\]
We also denote by
\[
\Pi:=I-\frac{1}{n}\mathbf{1}\mathbf{1}^\top
\]
the orthogonal projector onto $\mathbf{1}^\perp$.

\paragraph{Optimization problems.}
On the vector domains, we consider
\[
\min_{x\in\mathcal C}\ \Phi(x),
\qquad
\mathcal C\in\{\mathbb{R}_+^n,\Delta_n,[\ell,b]\}.
\]
On the Stiefel manifold, we write
\[
\min_{X\in\St(n,p)}\ \Phi(X).
\]
A main example on the vector domains is the least-squares objective
\[
\Phi(x)=\frac12\|Ax-c\|^2,
\qquad
\nabla\Phi(x)=A^\top(Ax-c),
\qquad
Q:=A^\top A\succeq 0.
\]

\paragraph{Constraint representations.}
For the orthant and box constraints, we use componentwise maps from $\mathbb{R}^n$ onto the interior of the feasible set.
\begin{itemize}[itemsep=3pt,topsep=4pt]
\item \textbf{Nonnegative orthant.} Two examples are the exponential and softplus maps,
\[
 g(u)=\exp(u),
 \qquad
 g(u)=\log(1+\exp(u)),
\]
where all operations are componentwise. In both cases, $g_i'(u_i)>0$ for every finite $u_i$, and
\[
 g_i(u_i)\to 0,
 \qquad
 g_i'(u_i)\to 0,
 \qquad
 \text{as }u_i\to-\infty.
\]

\item \textbf{Probability simplex.} We use the softmax map
\[
 x_i=\frac{e^{u_i}}{\sum_{j=1}^n e^{u_j}}.
\]
Its Jacobian is
\[
J_{\mathrm{sm}}(u)=\diag(x)-xx^\top.
\]
For $x\in\operatorname{ri}(\Delta_n)$, this matrix is symmetric positive semidefinite, has rank $n-1$, and satisfies
\[
\ker J_{\mathrm{sm}}(u)=\operatorname{span}\{\mathbf{1}\}.
\]
The softmax map is unchanged when a multiple of $\mathbf{1}$ is added to $u$.

\item \textbf{Box constraints.} We use
\[
 g(w)=\ell+(b-\ell)\odot\sigma(w),
 \qquad
 \sigma(t)=\frac{1}{1+e^{-t}},
\]
with componentwise operations. Its Jacobian is
\[
J_g(w)=\diag\!\big((b-\ell)\odot\sigma(w)\odot(1-\sigma(w))\big).
\]
It is diagonal and positive definite for every finite $w$. Its $i$th diagonal entry tends to zero when $w_i\to-\infty$ or $w_i\to+\infty$, which corresponds to $x_i\to\ell_i$ or $x_i\to b_i$.

\item \textbf{Stiefel constraint.} We do not use a global Euclidean reparameterization. The dynamics is defined directly on $\St(n,p)$.
\end{itemize}
Thus finite orthant and box parameters always represent interior points. Boundary points are obtained only as limits of diverging parameter sequences. The same observation holds for boundary points of the simplex in softmax coordinates.

\paragraph{Standing assumptions.}
We use the following assumptions only when they are required by the corresponding result.
\begin{enumerate}[label=\textbf{(A\arabic*)},itemsep=3pt,topsep=4pt]
\item \label{ass:smooth}
The objective is continuously differentiable on an open neighborhood of the feasible set, and its gradient is locally Lipschitz. When Hessians or Newton systems are used, the objective is twice continuously differentiable on the region under consideration.

\item \label{ass:Lipschitz}
When a global Lipschitz constant is needed, $\nabla\Phi$ is $L$-Lipschitz on the relevant level set or on another explicitly specified region containing the iterates.

\item \label{ass:coercive}
The objective is bounded below on the feasible set. Boundedness of the iterates, boundedness of a relevant level set, or coercivity is assumed separately whenever it is needed. These properties are not used interchangeably.

\item \label{ass:g}
The orthant and box maps are $C^2$ and componentwise strictly increasing, with $g_i'(z_i)>0$ for finite parameters. Hence their Jacobians are diagonal and positive definite in parameter space. The softmax case is treated either on the quotient by $\operatorname{span}\{\mathbf{1}\}$, after fixing a gauge, or directly in the primal variable.

\item \label{ass:stiefel}
For the Stiefel problem, $\Phi$ is the restriction of a $C^2$ function defined on an open neighborhood of $\St(n,p)$. We use the canonical Riemannian metric specified below.

\item \label{ass:KL}
In the nonconvex convergence results, the function under consideration satisfies the Kurdyka--\L{}ojasiewicz property at every point of the relevant compact level set \cite{bolte2007lojasiewicz,AttouchBolteSvaiter2013,BolteSabachTeboulle2014}.
\end{enumerate}
Additional assumptions, such as convexity, relative strong convexity, or positive definiteness of a local Hessian, will be stated with the results that use them.

\subsection{Continuous-Time Viewpoints}
\label{subsec:flows}

We now introduce the flows used in the paper. They share the same goal, namely to generate descent while keeping the motion feasible, but their constructions are different.

\paragraph{Euclidean gradient flow.}
For a smooth function $\Phi:\mathbb{R}^n\to\mathbb{R}$, the gradient flow is
\[
\dot x(t)=-\nabla\Phi(x(t)).
\]
Along every solution,
\[
\frac{d}{dt}\Phi(x(t))
=
\langle\nabla\Phi(x(t)),\dot x(t)\rangle
=
-\|\nabla\Phi(x(t))\|^2
\leq 0.
\]

\paragraph{Riemannian gradient flow.}
Let $(\mathcal M,\mathsf g)$ be a Riemannian manifold. The Riemannian gradient is defined by
\[
\mathsf g_y(\grad\Phi(y),\xi)=D\Phi(y)[\xi]
\qquad
\text{for every }\xi\in T_y\mathcal M.
\]
The corresponding flow is
\[
\dot y(t)=-\grad\Phi(y(t)),
\]
and it satisfies
\[
\frac{d}{dt}\Phi(y(t))
=-\|\grad\Phi(y(t))\|_y^2
\leq 0.
\]

On $\St(n,p)$, we use the canonical metric
\[
\mathsf g_X(\xi,\zeta)
=
\operatorname{tr}\!\left(\xi^\top\left(I-\frac12XX^\top\right)\zeta\right).
\]
If $G=\nabla\Phi(X)$ denotes the ambient Euclidean gradient, then
\begin{equation}
\grad\Phi(X)=G-XG^\top X=A(X)X,
\qquad
A(X):=GX^\top-XG^\top.
\label{eq:canonical-stiefel-gradient}
\end{equation}
The matrix $A(X)$ is skew-symmetric, and $-A(X)X$ is the canonical Riemannian gradient direction used later for the Cayley update; see \cite{EDELMAN1998}.

\paragraph{Reparameterized feasible flows.}
Let $g:\mathbb{R}^n\to\operatorname{int}(\mathcal C)$ be one of the componentwise orthant or box maps introduced above, and let $x=g(z)$. The parameter flow used in this paper is
\begin{equation}
\dot z(t)=-\nabla\Phi(g(z(t))).
\label{eq:prelim-parameter-flow}
\end{equation}
This is not the Euclidean gradient flow of $\Phi\circ g$. Its geometric interpretation will be given in Section~\ref{sec:riemannian-lens}. By the chain rule, the induced primal flow is
\begin{equation}
\dot x(t)
=
-J_g(z(t))\nabla\Phi(x(t)).
\label{eq:prelim-primal-mobility-flow}
\end{equation}
Since $J_g(z)$ is positive definite for finite $z$, we have
\[
\frac{d}{dt}\Phi(x(t))
=
-\nabla\Phi(x(t))^\top J_g(z(t))\nabla\Phi(x(t))
\leq 0.
\]
The diagonal entries of $J_g$ decrease when the corresponding primal variables approach an active face. This reduces the primal speed in the normal direction. As shown in Section~\ref{sec:kkt-from-flow}, the boundary sign conditions require the asymptotic behavior of the parameter dynamics in addition to the vanishing primal speed.

\paragraph{Simplex flow.}
The softmax Jacobian is singular on the full parameter space because softmax is invariant under translations along $\mathbf{1}$. We therefore write the simplex dynamics directly in the primal variable:
\begin{equation}
\dot x(t)
=
-\big(\diag(x(t))-x(t)x(t)^\top\big)\nabla\Phi(x(t)).
\label{eq:prelim-simplex-flow}
\end{equation}
The right-hand side belongs to $\mathbf{1}^\perp$, since
\[
\mathbf{1}^\top\big(\diag(x)-xx^\top\big)=0.
\]
Hence $\mathbf{1}^\top x(t)$ is constant. Moreover, the componentwise form
\[
\dot x_i=-x_i\big(\nabla_i\Phi(x)-x^\top\nabla\Phi(x)\big)
\]
shows that the simplex is invariant and that a trajectory starting in $\operatorname{ri}(\Delta_n)$ remains in the relative interior for every finite time for which the solution exists. Finally,
\[
\frac{d}{dt}\Phi(x(t))
=
-\nabla\Phi(x(t))^\top
\big(\diag(x(t))-x(t)x(t)^\top\big)
\nabla\Phi(x(t))
\leq 0.
\]

\paragraph{What is common.}
The vector reparameterizations, the simplex flow, and the Stiefel flow are not the same construction. The common point is simpler: the constraint determines the admissible motion, and the objective decreases along that motion.

\subsection{Implicit Discretization and Linear Stability}
\label{sec:stability}

We next recall the linear notion of A-stability. We then explain how it motivates the implicit steps used in the paper. A-stability is a property of a time-discretization method applied to a linear test equation. By itself, it does not prove nonlinear convergence, objective decrease, or accuracy for a large stepsize.

\subsubsection{Dahlquist's test equation}

Consider
\[
\dot y=\lambda y,
\qquad
\Re(\lambda)<0.
\]
A one-step method gives
\[
y_{k+1}=R(z)y_k,
\qquad
z=\eta\lambda,
\]
where $R$ is its stability function. Following the classical definition of Dahlquist \cite{dahlquist1963special}, the method is A-stable if
\[
|R(z)|\leq 1
\qquad
\text{for every }z\text{ such that }\Re(z)\leq 0.
\]
For the three classical methods used as references here,
\begin{center}
\begin{tabular}{lcc}
\toprule
Method & Stability function $R(z)$ & A-stable \\
\midrule
Explicit Euler & $1+z$ & No \\
Backward Euler & $\dfrac{1}{1-z}$ & Yes \\
Trapezoidal rule & $\dfrac{1+z/2}{1-z/2}$ & Yes \\
\bottomrule
\end{tabular}
\end{center}

For the quadratic function
\[
\Phi(x)=\frac12x^\top Qx,
\qquad
Q\succ0,
\]
the Euclidean gradient flow is $\dot x=-Qx$. Explicit Euler gives
\[
x_{k+1}=(I-\eta Q)x_k,
\]
which requires $0<\eta<2/\lambda_{\max}(Q)$ for linear stability. Backward Euler gives
\[
x_{k+1}=(I+\eta Q)^{-1}x_k,
\]
which is linearly stable for every $\eta>0$. It is also the proximal point step
\[
x_{k+1}
=
\arg\min_x
\left\{
\Phi(x)+\frac{1}{2\eta}\|x-x_k\|^2
\right\}.
\]

\subsubsection{Local linear models for the feasible flows}

The same spectral argument can be used locally, but only under the assumptions stated below.

\paragraph{Orthant and box.}
Consider the least-squares objective and the parameter flow \eqref{eq:prelim-parameter-flow}. Let $z^\star$ be a finite equilibrium and set
\[
x^\star=g(z^\star),
\qquad
D^\star:=J_g(z^\star)\succ0.
\]
Since $z^\star$ is finite, equilibrium means $\nabla\Phi(x^\star)=0$. The linearized parameter dynamics is
\[
\dot{\delta z}=-QD^\star\delta z.
\]
The matrix $QD^\star$ is similar to
\[
(D^\star)^{1/2}Q(D^\star)^{1/2}\succeq0.
\]
Thus the eigenvalues of the generator $-QD^\star$ are real and nonpositive. Equivalently, in primal coordinates,
\[
\dot{\delta x}=-D^\star Q\delta x.
\]
This calculation concerns finite interior equilibria. It does not describe a boundary KKT point, because such a point corresponds to a diverging parameter sequence and the derivative of the mobility also enters a boundary linearization.

\paragraph{Simplex.}
Let $x^\star\in\operatorname{ri}(\Delta_n)$ be a stationary point of \eqref{eq:prelim-simplex-flow}. Then
\[
\nabla\Phi(x^\star)=\tau^\star\mathbf{1}
\]
for some scalar $\tau^\star$. For a least-squares objective, the linearization on $\mathbf{1}^\perp$ is
\[
\dot{\delta x}=-M^\star Q\delta x,
\qquad
M^\star:=\diag(x^\star)-x^\star x^{\star\top}.
\]
Indeed, $M(x)\mathbf{1}=0$ for every $x$, so the derivative of $M(x)$ applied to the stationary gradient $\tau^\star\mathbf{1}$ vanishes. The operator $M^\star$ is positive definite on $\mathbf{1}^\perp$, and $M^\star Q$ is similar on this tangent space to
\[
(M^\star)^{1/2}Q(M^\star)^{1/2}\succeq0.
\]
Hence the tangent generator has real nonpositive eigenvalues. As in the orthant and box cases, this argument is an interior local model and is not a boundary stability theorem.

\paragraph{Stiefel manifold.}
At a Riemannian stationary point $X^\star$, the linearization of the gradient flow in local tangent coordinates is
\[
\dot\xi=-\operatorname{Hess}\Phi(X^\star)[\xi].
\]
At a local minimizer, the Riemannian Hessian is positive semidefinite, so the linearized generator is dissipative. The Cayley update used later is not backward Euler applied to this linearized equation. When the skew field $A(X)$ is frozen, it reduces to the trapezoidal rule for $\dot X=-A X$. This gives a useful linear stability interpretation, but it does not by itself imply descent of the nonlinear objective.

The conclusions above explain why implicit steps are attractive near stable equilibria. The nonlinear convergence results and the treatment of boundary limits require separate arguments.

\subsubsection{Implicit steps and proximal interpretations}
\label{sec:glance}

The implicit constructions used later have two main forms: Bregman steps on the vector domains and a Cayley equation on the Stiefel manifold.

\paragraph{Euclidean and Bregman steps.}
Backward Euler applied to the Euclidean gradient flow gives the Euclidean proximal point step displayed above. More generally, let $h$ be a Legendre function and consider the mirror flow
\[
\frac{d}{dt}\nabla h(x(t))=-\nabla\Phi(x(t)).
\]
Backward Euler in the dual variable gives
\[
\nabla h(x_{k+1})-\nabla h(x_k)+\eta\nabla\Phi(x_{k+1})=0,
\]
with a normal-space term when the feasible set contains an affine constraint. This is the first-order condition of
\[
x_{k+1}
\in
\arg\min_{x\in\mathcal C}
\left\{
D_h(x,x_k)+\eta\Phi(x)
\right\}.
\]
When $\Phi$ is convex and the subproblem has a minimizer, this minimizer is unique if $h$ is strictly convex.

\paragraph{Reparameterized orthant and box steps.}
Backward Euler applied to \eqref{eq:prelim-parameter-flow} gives
\begin{equation}
z_{k+1}-z_k+\eta\nabla\Phi(g(z_{k+1}))=0.
\label{eq:prelim-reparameterized-BE}
\end{equation}
Section~\ref{sec:riemannian-lens} shows that every componentwise increasing diffeomorphism $g$ induces a separable potential $h$ such that
\[
\nabla h(x)=g^{-1}(x).
\]
Consequently, \eqref{eq:prelim-reparameterized-BE} is exactly the dual backward Euler equation
\[
\nabla h(x_{k+1})-\nabla h(x_k)+\eta\nabla\Phi(x_{k+1})=0.
\]
Thus it is a Bregman proximal step in the primal variable $x$. In the parameter variable $z$, it is implemented as a nonlinear root equation. Its Jacobian need not be symmetric, even though the equivalent primal Bregman Hessian is symmetric.

For the exponential orthant map, $h$ is the negative entropy on $\mathbb R_+^n$. For the logistic box map, $h$ is a two-sided entropy. On the simplex, the same construction with the negative entropy and the affine mass constraint gives the KL-proximal step.

\paragraph{Stiefel Cayley step.}
The Stiefel algorithm uses the implicit equation
\[
X_{k+1}
=
\left(I+\frac{\eta}{2}A(X_{k+1})\right)^{-1}
\left(I-\frac{\eta}{2}A(X_{k+1})\right)X_k.
\]
At every solution, the Cayley factor is orthogonal because $A(X_{k+1})$ is skew-symmetric. Hence $X_k\in\St(n,p)$ implies $X_{k+1}\in\St(n,p)$. We use this equation as a structure-preserving implicit step. We do not identify it with a proximal map.

\paragraph{Conclusion.}
On the vector domains, the reparameterized backward Euler equations are exact Bregman steps when written in the primal variable. On the Stiefel manifold, the implicit step is a feasible Cayley equation. A-stability explains the absence of a linear stability cap in the corresponding frozen linear models. Descent and nonlinear convergence require separate arguments.

\section{A Geometric Lens: Constraint-Induced Metrics and Feasible Flows}
\label{sec:riemannian-lens}

This section gives the precise geometric interpretation of the flows introduced in Section~\ref{sec:preliminaries}. The main point is that a constraint representation does more than guarantee feasibility. It also defines a metric in which the corresponding dynamics is a gradient flow. For the orthant and box, the metric is induced by the componentwise change of variables. On the simplex, it is the Fisher--Shahshahani metric generated by the negative entropy. On the Stiefel manifold, it is the canonical Riemannian metric.

The geometries are therefore not identical, but they follow the same
construction. A metric determines the admissible gradient direction, and its
inverse acts as the mobility of the primal flow. For the orthant and box maps
used below, the metric becomes singular, or equivalently the mobility
vanishes, as an active face is approached. More generally, this conclusion
requires the derivative of the feasibility map to tend to zero at the
corresponding endpoint. This boundary behavior will be used in
Section~\ref{sec:kkt-from-flow}.

\subsection{Riemannian Gradient Flows}
\label{subsec:riem-template}

Let $\mathcal C\subset\mathbb R^n$ be open, and let $G(x)\succ0$ be a smooth metric tensor. The inner product at $x$ is
\[
\langle v,w\rangle_x=v^\top G(x)w.
\]
The Riemannian gradient of a smooth function $f$ is
\[
\grad_G f(x)=G(x)^{-1}\nabla f(x),
\]
and the corresponding gradient flow is
\begin{equation}
\dot x(t)=-G(x(t))^{-1}\nabla f(x(t)).
\label{eq:general-riem-flow}
\end{equation}
Along every solution,
\[
\frac{d}{dt}f(x(t))
=-\nabla f(x(t))^\top G(x(t))^{-1}\nabla f(x(t))
=-\|\grad_G f(x(t))\|_{x(t)}^2
\leq0.
\]
Thus $G^{-1}$ is the mobility operator of the flow.

A direct backward Euler discretization of \eqref{eq:general-riem-flow} is
\begin{equation}
\frac{x_{k+1}-x_k}{\eta}
+G(x_{k+1})^{-1}\nabla f(x_{k+1})=0.
\label{eq:implicit-euler-riem}
\end{equation}
If $G$ is constant, this is the first-order condition of the metric proximal problem
\[
x_{k+1}\in\argmin_x
\left\{
f(x)+\frac{1}{2\eta}\|x-x_k\|_G^2
\right\}.
\]
For a state-dependent metric, \eqref{eq:implicit-euler-riem} is not, in general, the optimality condition of one proximal subproblem. A sounder proximal interpretation is obtained when the metric is Hessian, as shown next.

\subsection{Componentwise Reparameterizations and Hessian Geometry}
\label{sec:jacobian-nonpullback}

Let
\[
\mathcal C=I_1\times\cdots\times I_n,
\]
where each $I_i$ is an open interval. Assume that
\[
g_i:\mathbb R\to I_i
\]
is a $C^2$ increasing diffeomorphism, and set
\[
x=g(u),
\qquad
J_g(u)=\diag(g_1'(u_1),\ldots,g_n'(u_n)).
\]
The following result identifies the geometry induced by this change of variables.

\begin{proposition}[Hessian geometry induced by a componentwise map]
\label{prop:reparam-hessian-geometry}
Choose $\bar x_i\in I_i$ and define
\begin{equation}
h(x)=\sum_{i=1}^n h_i(x_i),
\qquad
h_i(x_i)=\int_{\bar x_i}^{x_i}g_i^{-1}(s)\,ds.
\label{eq:induced-potential}
\end{equation}
Then $h$ is strictly convex on $\mathcal C$ and
\begin{equation}
\nabla h(x)=g^{-1}(x),
\qquad
\nabla^2h(x)=J_g(g^{-1}(x))^{-1}.
\label{eq:induced-hessian}
\end{equation}
Moreover, the pullback of the Hessian metric $\nabla^2h(x)$ through $x=g(u)$ is
\begin{equation}
J_g(u)^\top\nabla^2h(g(u))J_g(u)=J_g(u).
\label{eq:pullback-induced-metric}
\end{equation}
\end{proposition}

\begin{proof}
Since $h_i'(x_i)=g_i^{-1}(x_i)$, the inverse function theorem gives
\[
h_i''(x_i)
=(g_i^{-1})'(x_i)
=\frac{1}{g_i'(g_i^{-1}(x_i))}>0.
\]
This proves \eqref{eq:induced-hessian} and the strict convexity of $h$. Since $J_g(u)$ is diagonal and symmetric,
\[
J_g(u)^\top\nabla^2h(g(u))J_g(u)
=J_g(u)J_g(u)^{-1}J_g(u)
=J_g(u),
\]
which proves \eqref{eq:pullback-induced-metric}.
\end{proof}

Proposition~\ref{prop:reparam-hessian-geometry} gives two equivalent descriptions of the same geometry. In the primal variable $x$, the metric and its inverse mobility are
\[
G_x(x)=\nabla^2h(x),
\qquad
M_g(x)=G_x(x)^{-1}=J_g(g^{-1}(x)).
\]
In the parameter variable $u$, the metric is
\[
G_u(u)=J_g(u).
\]
This is not the pullback of the Euclidean metric, which would be $J_g(u)^\top J_g(u)$. It is the pullback of the Hessian metric induced by $h$.

Let $\Psi(u)=\Phi(g(u))$. Since
\[
\nabla\Psi(u)=J_g(u)^\top\nabla\Phi(g(u)),
\]
the Riemannian gradient of $\Psi$ for the metric $G_u(u)=J_g(u)$ is
\begin{equation}
\grad_{G_u}\Psi(u)
=G_u(u)^{-1}\nabla\Psi(u)
=\nabla\Phi(g(u)).
\label{eq:grad-def}
\end{equation}
The corresponding parameter flow is
\begin{equation}
\dot u(t)=-\nabla\Phi(g(u(t))),
\label{eq:riem-flow-u}
\end{equation}
and its primal image is
\begin{equation}
\dot x(t)
=-J_g(u(t))\nabla\Phi(x(t))
=-M_g(x(t))\nabla\Phi(x(t)).
\label{eq:non-pullback-flow}
\end{equation}
Equivalently, using \eqref{eq:induced-hessian},
\[
\frac{d}{dt}\nabla h(x(t))=-\nabla\Phi(x(t)).
\]
Hence the reparameterized dynamics is exactly a mirror, or Hessian--Riemannian, gradient flow in the primal variable.

The objective decreases along this flow:
\[
\frac{d}{dt}\Phi(x(t))
=-\nabla\Phi(x(t))^\top M_g(x(t))\nabla\Phi(x(t))
\leq0.
\]
The metric is defined on the interior of $\mathcal C$. Whenever
$g_i'(u_i)\to0$ at the relevant endpoint---as for the exponential, softplus,
and logistic maps considered here---the corresponding entry of $M_g(x)$ tends
to zero as that coordinate approaches an active face, and the primal motion
slows down in that coordinate. A zero primal velocity at the boundary is not
sufficient by itself to obtain the KKT sign condition. The additional argument
based on the parameter dynamics is given in
Section~\ref{sec:kkt-from-flow}.

\paragraph{Orthant example.}
For $g_i(u_i)=e^{u_i}$,
\[
h_i(x_i)=x_i\log x_i-x_i,
\qquad
h_i''(x_i)=\frac{1}{x_i},
\qquad
M_g(x)=\diag(x).
\]
The primal flow is
\[
\dot x=-x\odot\nabla\Phi(x).
\]
Other increasing diffeomorphisms, such as the softplus map, generate a different separable potential through \eqref{eq:induced-potential}.

\paragraph{Box example.}
Let $d_i=b_i-\ell_i$ and
\[
x_i=\ell_i+d_i\sigma(w_i).
\]
Then
\[
w_i=\log\frac{x_i-\ell_i}{b_i-x_i}.
\]
The induced potential can be chosen as
\begin{equation}
h_{\Box}(x)
=\sum_{i=1}^n
\left[
(x_i-\ell_i)\log\frac{x_i-\ell_i}{d_i}
+(b_i-x_i)\log\frac{b_i-x_i}{d_i}
\right].
\label{eq:box-induced-potential}
\end{equation}
It satisfies
\[
\nabla_i h_{\Box}(x)=w_i,
\qquad
\nabla_{ii}^2h_{\Box}(x)
=\frac{d_i}{(x_i-\ell_i)(b_i-x_i)}.
\]
Therefore,
\[
M_g(x)
=\diag\!\left(
\frac{(x_i-\ell_i)(b_i-x_i)}{b_i-\ell_i}
\right)_{i=1}^n.
\]
Thus the logistic reparameterization is exactly the mirror geometry generated by the two-sided entropy \eqref{eq:box-induced-potential}. It is not the logarithmic-barrier geometry, which would give a different mobility.

\subsection{Implicit Step on Product Domains}

Applying backward Euler to \eqref{eq:riem-flow-u} gives
\begin{equation}
u_{k+1}-u_k+\eta\nabla\Phi(x_{k+1})=0,
\qquad
x_{k+1}=g(u_{k+1}).
\label{eq:bregman-dual-implicit}
\end{equation}
Since $u=\nabla h(x)$, this is equivalent to
\[
\nabla h(x_{k+1})-\nabla h(x_k)
+\eta\nabla\Phi(x_{k+1})=0.
\]
Whenever the corresponding minimization problem has a solution, this is the first-order condition of
\begin{equation}
x_{k+1}\in\argmin_{x\in\mathcal C}
\left\{
D_h(x,x_k)+\eta\Phi(x)
\right\}.
\label{eq:bregman-prox}
\end{equation}
If $\Phi$ is convex, the objective in \eqref{eq:bregman-prox} is strictly convex, and the solution is unique whenever it exists. For a nonconvex objective, \eqref{eq:bregman-dual-implicit} is an implicit critical-point equation, and a selected root need not be a global minimizer of \eqref{eq:bregman-prox}.

In the parameter variable, the same step is the root equation
\begin{equation}
F_k(u)
:=u-u_k+\eta\nabla\Phi(g(u))=0.
\label{eq:BE-nonpullback}
\end{equation}
If $\Phi\in C^2$, then
\begin{equation}
J_{F_k}(u)
=I+\eta\nabla^2\Phi(g(u))J_g(u)
=\bigl(\nabla^2h(x)+\eta\nabla^2\Phi(x)\bigr)J_g(u).
\label{eq:parameter-residual-jacobian}
\end{equation}
For a convex objective, the first factor in the last expression is positive definite. Hence $J_{F_k}(u)$ is nonsingular. The matrix is generally nonsymmetric, and the identity term alone does not give a uniform bound on its Euclidean condition number. This point will be important when we discuss the inner solvers in Section~\ref{sec:algorithms}.

\subsection{Fisher--Shahshahani Geometry on the Simplex}
\label{subsec:canonical-metrics}

On $\operatorname{ri}(\Delta_n)$, consider the negative entropy
\[
h_{\Delta}(x)=\sum_{i=1}^n x_i\log x_i.
\]
For tangent vectors $v,w\in T_x\Delta_n=\{z:\mathbf 1^\top z=0\}$, the Hessian metric is
\[
\langle v,w\rangle_x^{\mathrm{FS}}
=\sum_{i=1}^n\frac{v_iw_i}{x_i}.
\]
Let
\[
M_{\Delta}(x)=\diag(x)-xx^\top.
\]
For every $q\in\mathbb R^n$ and every $\xi\in T_x\Delta_n$,
\[
\left\langle M_{\Delta}(x)q,\xi\right\rangle_x^{\mathrm{FS}}
=q^\top\xi.
\]
Thus $M_{\Delta}(x)$ is the inverse metric operator acting on ambient gradients. The Riemannian gradient flow is
\[
\dot x=-M_{\Delta}(x)\nabla\Phi(x),
\]
which is the replicator flow introduced in Section~\ref{sec:preliminaries}.

Backward Euler in the entropy dual variable, with the affine constraint $\mathbf 1^\top x=1$, gives
\[
\log x_{k+1}-\log x_k
+\eta\nabla\Phi(x_{k+1})
+\nu_{k+1}\mathbf 1=0.
\]
This is the optimality condition of
\[
x_{k+1}\in\argmin_{x\in\Delta_n}
\left\{
\KL(x\|x_k)+\eta\Phi(x)
\right\}.
\]
If $x_k\in\operatorname{ri}(\Delta_n)$ and $\Phi$ is smooth and finite on $\Delta_n$, every minimizer is in the relative interior. Boundary points may still appear as limits of the iterates.

\subsection{Canonical Geometry on the Stiefel Manifold}

On $\St(n,p)$, we use the canonical metric introduced in Section~\ref{sec:preliminaries}. If $G=\nabla\Phi(X)$ and
\[
A(X)=GX^\top-XG^\top,
\]
then
\[
\grad\Phi(X)=A(X)X.
\]
The Riemannian gradient flow is
\[
\dot X=-A(X)X,
\]
and
\[
\frac{d}{dt}\Phi(X(t))
=-\|\grad\Phi(X(t))\|_{X(t)}^2
\leq0.
\]

The implicit Stiefel step used later is the Cayley equation
\[
\left(I+\frac{\eta}{2}A(X_{k+1})\right)X_{k+1}
=
\left(I-\frac{\eta}{2}A(X_{k+1})\right)X_k.
\]
Every solution is feasible because the Cayley factor is orthogonal. This step is not a Bregman proximal map. When $A$ is frozen, it is the trapezoidal rule for the linear equation $\dot X=-AX$.

\paragraph{Summary.}
The orthant, box, and simplex flows are Hessian--Riemannian or mirror flows on the interiors of their feasible sets. Their inverse metrics are the mobility operators that appear in the primal dynamics. The Stiefel flow is Riemannian directly on the manifold. This gives a common geometric interpretation while keeping the different boundary and discretization mechanisms explicit.

\section{The \gravidy\ Framework: Constraint-Induced Feasible Dynamics}
\label{sec:gravidy}

We now collect the constructions of the previous section under the name \gravidy. The framework consists of a feasible gradient flow and an implicit discretization built from a metric adapted to the constraint. The metric is not determined by the feasible set alone. For the orthant and box, it also depends on the chosen reparameterization. For the simplex and Stiefel manifold, we use the standard entropy and canonical metrics.

The common principle is the following:
\begin{enumerate}[label=\textbf{(G\arabic*)},itemsep=3pt,topsep=4pt]
\item choose coordinates or an intrinsic metric that represent the constraint;
\item form the corresponding feasible gradient flow;
\item discretize this flow implicitly while preserving feasibility.
\end{enumerate}
The four versions used in this paper are described below.

\subsection{The Four \gravidy\ Dynamics}
\label{subsec:family}

\paragraph{\gravidypos\ on the nonnegative orthant.}
Let $x=g(u)$, where $g_i:\mathbb R\to(0,+\infty)$ is an increasing diffeomorphism. The induced potential $h$ is defined by
\[
\nabla h(x)=g^{-1}(x).
\]
The parameter and primal flows are
\begin{equation}
\dot u=-\nabla\Phi(x),
\qquad
\dot x=-J_g(u)\nabla\Phi(x).
\label{eq:jac-induced-flow}
\end{equation}
For the exponential map,
\[
\dot x=-x\odot\nabla\Phi(x),
\qquad
h(x)=\sum_i(x_i\log x_i-x_i).
\]
The implicit step can be written either in the parameter variable,
\[
u_{k+1}-u_k+\eta\nabla\Phi(x_{k+1})=0,
\qquad x_{k+1}=g(u_{k+1}),
\]
or as the Bregman step generated by $h$. Finite iterates are positive, while zero components may appear in the limit.

\paragraph{\gravidydelta\ on the simplex.}
For $x\in\operatorname{ri}(\Delta_n)$, the inverse Fisher--Shahshahani metric is
\[
M_{\Delta}(x)=\diag(x)-xx^\top.
\]
The flow is
\begin{equation}
\dot x=-M_{\Delta}(x)\nabla\Phi(x).
\label{eq:replicator-maintext}
\end{equation}
It preserves nonnegativity and $\mathbf 1^\top x=1$. The implicit step is
\[
x_{k+1}\in\argmin_{x\in\Delta_n}
\left\{
\KL(x\|x_k)+\eta\Phi(x)
\right\}.
\]
Starting from the relative interior, every finite iterate remains in the relative interior under the smoothness assumptions stated earlier.

\paragraph{\gravidybox\ under box constraints.}
For
\[
x=\ell+(b-\ell)\odot\sigma(w),
\]
the primal flow is
\[
\dot x=-M_{\Box}(x)\nabla\Phi(x),
\]
where
\[
M_{\Box}(x)
=\diag\!\left(
\frac{(x_i-\ell_i)(b_i-x_i)}{b_i-\ell_i}
\right)_{i=1}^n.
\]
This is the mirror flow generated by the two-sided entropy \eqref{eq:box-induced-potential}. Its implicit Bregman step is exactly the backward Euler equation in the logistic variable $w$. Finite iterates stay in the open box, and active bounds may be reached only in the limit.

\paragraph{\gravidyst\ on the Stiefel manifold.}
For $X\in\St(n,p)$, let
\[
A(X)=\nabla\Phi(X)X^\top-X\nabla\Phi(X)^\top.
\]
Under the canonical metric,
\[
\grad\Phi(X)=A(X)X,
\qquad
\dot X=-A(X)X.
\]
The implicit update is
\[
\left(I+\frac{\eta}{2}A(X_{k+1})\right)X_{k+1}
=
\left(I-\frac{\eta}{2}A(X_{k+1})\right)X_k.
\]
Every exact solution satisfies $X_{k+1}^\top X_{k+1}=I_p$.

\subsection{Properties Shared by the Four Constructions}
\label{subsec:desiderata}

\paragraph{Feasibility.}
For the orthant and box, feasibility follows from the maps $g$. For the simplex, the flow is tangent and the KL step satisfies the mass constraint. For the Stiefel manifold, the continuous flow is tangent and the Cayley factor is orthogonal. Thus feasibility is part of the dynamics and not a correction applied afterwards.

\paragraph{Descent of the continuous flows.}
On the vector domains, let $M(x)$ denote the corresponding mobility. Then
\[
\frac{d}{dt}\Phi(x(t))
=-\nabla\Phi(x(t))^\top M(x(t))\nabla\Phi(x(t))
\leq0.
\]
On the Stiefel manifold,
\[
\frac{d}{dt}\Phi(X(t))
=-\|\grad\Phi(X(t))\|_{X(t)}^2
\leq0.
\]
The discrete descent properties require separate arguments. For the vector Bregman steps, they follow from the proximal subproblem when the selected point is a minimizer. For the nonlinear Cayley step, descent is not automatic and will be enforced by the acceptance rule used later.

\paragraph{First-order optimality.}
The relation with optimality must be stated carefully. On the Stiefel manifold, a stationary point of the Riemannian flow is exactly a first-order critical point. On the orthant, simplex, and box, a zero of the degenerate primal vector field at the boundary is not sufficient. Section~\ref{sec:kkt-from-flow} proves instead that convergent trajectories generated from the interior satisfy the corresponding KKT conditions. The same conclusion holds for convergent implicit iterates generated from the interior when their stepsizes are bounded away from zero.

\paragraph{Implicit equations.}
The vector-domain methods lead to Bregman optimality systems. For the orthant and box, these systems are implemented through the parameter residual
\[
F_k(u)=u-u_k+\eta\nabla\Phi(g(u)).
\]
For a twice differentiable objective,
\[
J_{F_k}(u)=I+\eta\nabla^2\Phi(g(u))J_g(u).
\]
For convex $\Phi$, this Jacobian is nonsingular, but it is generally nonsymmetric and need not be well conditioned in the Euclidean norm. On the simplex, the primal Hessian of the KL subproblem is
\[
\diag(1/x)+\eta\nabla^2\Phi(x),
\]
which is positive definite in the convex case. On the Stiefel manifold, the Cayley equation gives a nonlinear matrix residual whose Fr\'echet derivative is used by the Newton solvers. These systems are studied in detail in Section~\ref{sec:algorithms}.

\subsection{From Geometry to Optimality}
\label{subsec:energy-kkt-brief}

The role of the framework can be summarized as
\begin{center}
\fbox{%
\parbox{0.82\linewidth}{%
\[
\begin{aligned}
\text{constraint and representation}
&\;\Longrightarrow\;\text{metric and mobility}\\
&\;\Longrightarrow\;\text{feasible descent flow}\\
&\;\Longrightarrow\;\text{implicit feasible step}\\
&\;\Longrightarrow\;\text{convergent limit}\\
&\;\Longrightarrow\;\text{first-order optimality.}
\end{aligned}
\]
}}
\end{center}
For the vector domains, the last implication uses the asymptotic behavior of the dual or parameter variables. For the Stiefel manifold, it follows directly from Riemannian stationarity. This is the precise form of the unifying principle used in the remainder of the paper.

\section{Stationarity of the Dynamics and First-Order Optimality}
\label{sec:kkt-from-flow}

We now study the relation between the feasible dynamics and the first-order optimality conditions of the constrained problem. The argument is simple in the interior, where the mobility is nonsingular. More care is needed at the boundary. Indeed, the primal velocity may vanish only because the mobility vanishes. For example, for the orthant flow
\[
\dot x_i=-x_i\nabla_i\Phi(x),
\]
we have \(\dot x_i=0\) whenever \(x_i=0\), independently of the sign of \(\nabla_i\Phi(x)\). Hence, a zero of the primal vector field is not sufficient to obtain the KKT conditions.

The missing sign condition is recovered from the variables that generate the geometry. If a trajectory converges to an active face, a gradient with the wrong sign would move the corresponding parameter away from that face. The same argument applies to the implicit iterations. On the simplex, it is convenient to use the logarithms of the coordinates, while on the Stiefel manifold there is no boundary and the usual Riemannian stationarity condition applies directly.

Throughout this section, \(\Phi\) satisfies Assumption~\ref{ass:smooth}. For the orthant, simplex, and box, the continuous trajectories and the discrete iterates are initialized in the interior of the feasible set.

\subsection{Nonnegative Orthant}
\label{subsec:kkt-orthant}

We first consider
\[
\min_{x\geq 0}\Phi(x).
\]
Let \(x_i=g_i(u_i)\), where each \(g_i:\mathbb R\to(0,+\infty)\) is a \(C^1\) increasing diffeomorphism such that
\[
\lim_{u_i\to-\infty}g_i(u_i)=0,
\qquad
\lim_{u_i\to-\infty}g_i'(u_i)=0,
\qquad
 g_i'(u_i)>0 \quad\text{for every finite }u_i.
\]
The parameter and primal dynamics are
\begin{equation}
\dot u=-\nabla\Phi(x),
\qquad
x=g(u),
\label{eq:orthant-u-flow-kkt}
\end{equation}
and
\begin{equation}
\dot x_i=-g_i'(u_i)\nabla_i\Phi(x).
\label{eq:orthant-x-flow-kkt}
\end{equation}

\begin{proposition}[Convergent orthant trajectories satisfy the KKT conditions]
\label{prop:kkt-orthant-flow}
Let \(x(t)=g(u(t))\) solve \eqref{eq:orthant-u-flow-kkt} for \(t\geq 0\). Assume that
\[
x(t)\longrightarrow x^\star\in\mathbb R_+^n
\qquad\text{as }t\to+\infty.
\]
Then
\[
x^\star\geq 0,
\qquad
\nabla\Phi(x^\star)\geq 0,
\qquad
x_i^\star\nabla_i\Phi(x^\star)=0
\quad\text{for all }i.
\]
Hence \(x^\star\) satisfies the KKT conditions of the orthant-constrained problem.
\end{proposition}

\begin{proof}
Fix an index \(i\).

Assume first that \(x_i^\star>0\). Since \(g_i\) is a diffeomorphism, \(u_i(t)\) converges to the finite value \(g_i^{-1}(x_i^\star)\). Suppose that \(\nabla_i\Phi(x^\star)\neq 0\). By continuity, \(\dot u_i(t)=-\nabla_i\Phi(x(t))\) then has a fixed nonzero sign and is bounded away from zero for all sufficiently large \(t\). This is impossible because \(u_i(t)\) converges. Therefore,
\[
\nabla_i\Phi(x^\star)=0.
\]

Assume now that \(x_i^\star=0\). Then \(u_i(t)\to-\infty\). Suppose that \(\nabla_i\Phi(x^\star)<0\). There exist \(\varepsilon>0\) and \(T>0\) such that
\[
\nabla_i\Phi(x(t))\leq-\varepsilon
\qquad\text{for every }t\geq T.
\]
It follows that
\[
\dot u_i(t)=-\nabla_i\Phi(x(t))\geq\varepsilon
\qquad\text{for every }t\geq T.
\]
Hence \(u_i(t)\geq u_i(T)+\varepsilon(t-T)\), which contradicts \(u_i(t)\to-\infty\). We conclude that
\[
\nabla_i\Phi(x^\star)\geq 0.
\]

Thus, \(\nabla_i\Phi(x^\star)=0\) when \(x_i^\star>0\), while \(\nabla_i\Phi(x^\star)\geq0\) when \(x_i^\star=0\). These are exactly the KKT conditions.
\end{proof}

For the exponential map, \(x_i=e^{u_i}\), the primal flow is
\[
\dot x_i=-x_i\nabla_i\Phi(x).
\]
The identity \(x_i\nabla_i\Phi(x)=0\) gives complementarity. The proof above is still needed to obtain the inequality \(\nabla_i\Phi(x^\star)\geq0\) at a zero component.

\begin{proposition}[Convergent implicit orthant iterates satisfy the KKT conditions]
\label{prop:kkt-orthant-discrete}
Let \((u_k)\) be generated by
\begin{equation}
 u_{k+1}-u_k+\eta_k\nabla\Phi(x_{k+1})=0,
 \qquad x_{k+1}=g(u_{k+1}),
\label{eq:orthant-implicit-kkt}
\end{equation}
where \(\eta_k\geq\underline\eta>0\). If
\[
x_k\longrightarrow x^\star\in\mathbb R_+^n,
\]
then \(x^\star\) satisfies the orthant KKT conditions.
\end{proposition}

\begin{proof}
Fix an index \(i\). If \(x_i^\star>0\), then \(u_{k,i}\to g_i^{-1}(x_i^\star)\), and therefore \(u_{k+1,i}-u_{k,i}\to0\). From \eqref{eq:orthant-implicit-kkt},
\[
\eta_k\nabla_i\Phi(x_{k+1})=-(u_{k+1,i}-u_{k,i}).
\]
Since \(\eta_k\geq\underline\eta\), we obtain \(\nabla_i\Phi(x^\star)=0\).

If \(x_i^\star=0\), then \(u_{k,i}\to-\infty\). Suppose that \(\nabla_i\Phi(x^\star)<0\). For all sufficiently large \(k\),
\[
\nabla_i\Phi(x_{k+1})\leq-\varepsilon
\]
for some \(\varepsilon>0\). Hence
\[
u_{k+1,i}-u_{k,i}
=-\eta_k\nabla_i\Phi(x_{k+1})
\geq\underline\eta\varepsilon.
\]
The sequence \((u_{k,i})\) is then eventually increasing by at least a fixed positive amount, which contradicts \(u_{k,i}\to-\infty\). Therefore, \(\nabla_i\Phi(x^\star)\geq0\). The KKT conditions follow.
\end{proof}

For \(g_i(u_i)=e^{u_i}\), \eqref{eq:orthant-implicit-kkt} can also be written as
\[
x_{k+1,i}
=x_{k,i}\exp\!\bigl(-\eta_k\nabla_i\Phi(x_{k+1})\bigr).
\]
This form makes the boundary sign condition easy to interpret. A negative limiting gradient would eventually multiply \(x_{k,i}\) by a factor strictly larger than one, and the iterates could not converge to zero.

\subsection{Simplex}
\label{subsec:kkt-simplex}

We now consider
\[
\min_{x\in\Delta_n}\Phi(x).
\]
The \gravidydelta\ flow is
\begin{equation}
\dot x=-\bigl(\diag(x)-xx^\top\bigr)\nabla\Phi(x).
\label{eq:simplex-flow-kkt}
\end{equation}
Equivalently,
\begin{equation}
\dot x_i=-x_i\bigl(\nabla_i\Phi(x)-\tau(x)\bigr),
\qquad
\tau(x):=x^\top\nabla\Phi(x).
\label{eq:simplex-replicator-kkt}
\end{equation}
If \(x(0)\in\operatorname{ri}(\Delta_n)\), then \(x_i(t)>0\) for every \(i\) and every finite \(t\). We may therefore write
\begin{equation}
\frac{d}{dt}\log x_i(t)
=-\bigl(\nabla_i\Phi(x(t))-\tau(x(t))\bigr).
\label{eq:simplex-log-flow-kkt}
\end{equation}

\begin{proposition}[Convergent simplex trajectories satisfy the KKT conditions]
\label{prop:kkt-simplex-flow}
Let \(x(t)\) solve \eqref{eq:simplex-flow-kkt} with \(x(0)\in\operatorname{ri}(\Delta_n)\). Assume that
\[
x(t)\longrightarrow x^\star\in\Delta_n.
\]
Set
\[
\tau^\star:=x^{\star\top}\nabla\Phi(x^\star).
\]
Then
\[
x_i^\star>0
\quad\Longrightarrow\quad
\nabla_i\Phi(x^\star)=\tau^\star,
\]
and
\[
x_i^\star=0
\quad\Longrightarrow\quad
\nabla_i\Phi(x^\star)\geq\tau^\star.
\]
Hence \(x^\star\) satisfies the KKT conditions of the simplex-constrained problem.
\end{proposition}

\begin{proof}
If \(x_i^\star>0\), then \(\log x_i(t)\) converges to \(\log x_i^\star\). By continuity, the right-hand side of \eqref{eq:simplex-log-flow-kkt} converges to
\[
-\bigl(\nabla_i\Phi(x^\star)-\tau^\star\bigr).
\]
This limit must be zero. Otherwise, \(\log x_i(t)\) would eventually increase or decrease at a rate bounded away from zero and could not converge. Thus,
\[
\nabla_i\Phi(x^\star)=\tau^\star.
\]

Assume now that \(x_i^\star=0\). Then \(\log x_i(t)\to-\infty\). Suppose that \(\nabla_i\Phi(x^\star)<\tau^\star\). There exist \(\varepsilon>0\) and \(T>0\) such that
\[
\nabla_i\Phi(x(t))-\tau(x(t))\leq-\varepsilon
\qquad\text{for every }t\geq T.
\]
By \eqref{eq:simplex-log-flow-kkt},
\[
\frac{d}{dt}\log x_i(t)\geq\varepsilon
\qquad\text{for every }t\geq T,
\]
which contradicts \(\log x_i(t)\to-\infty\). Therefore,
\[
\nabla_i\Phi(x^\star)\geq\tau^\star.
\]

Define
\[
\nu^\star:=-\tau^\star,
\qquad
\lambda_i^\star:=\nabla_i\Phi(x^\star)-\tau^\star.
\]
Then
\[
\nabla\Phi(x^\star)+\nu^\star\mathbf 1-\lambda^\star=0,
\qquad
\lambda^\star\geq0,
\qquad
\lambda_i^\star x_i^\star=0.
\]
Together with \(x^\star\in\Delta_n\), this is the simplex KKT system.
\end{proof}

\begin{proposition}[Convergent KL-proximal iterates satisfy the KKT conditions]
\label{prop:kkt-simplex-discrete}
Let \(x_0\in\operatorname{ri}(\Delta_n)\), and let
\begin{equation}
 x_{k+1}\in\argmin_{x\in\Delta_n}
 \left\{\KL(x\|x_k)+\eta_k\Phi(x)\right\},
\label{eq:simplex-kl-kkt}
\end{equation}
where \(\eta_k\geq\underline\eta>0\). Assume that
\[
x_k\longrightarrow x^\star\in\Delta_n.
\]
Then \(x^\star\) satisfies the simplex KKT conditions.
\end{proposition}

\begin{proof}
We first note that \(x_{k+1}\) is interior. Suppose, on the contrary, that \(x_{k+1,i}=0\) for some \(i\). Moving a small mass from a positive coordinate to the \(i\)-th coordinate changes the KL term by an amount containing \(s\log s\), while the change in \(\Phi\) is of order \(s\). The objective therefore decreases for all sufficiently small \(s>0\), which is a contradiction.

The first-order condition of \eqref{eq:simplex-kl-kkt} is
\begin{equation}
\log x_{k+1}-\log x_k
+\eta_k\nabla\Phi(x_{k+1})
+c_{k+1}\mathbf 1=0,
\label{eq:simplex-kl-opt-kkt}
\end{equation}
where the scalar \(c_{k+1}\) contains the multiplier of \(\mathbf 1^\top x=1\), as well as the constant arising from the derivative of the entropy.

Let
\[
S:=\{i:x_i^\star>0\}.
\]
The set \(S\) is nonempty. Choose \(j\in S\). Since \(x_{k,j}\to x_j^\star>0\),
\[
\log\frac{x_{k+1,j}}{x_{k,j}}\longrightarrow0.
\]
Subtracting the \(j\)-th equation in \eqref{eq:simplex-kl-opt-kkt} from the \(i\)-th one gives
\begin{equation}
\log\frac{x_{k+1,i}}{x_{k,i}}
-
\log\frac{x_{k+1,j}}{x_{k,j}}
=-\eta_k\bigl(\nabla_i\Phi(x_{k+1})-\nabla_j\Phi(x_{k+1})\bigr).
\label{eq:simplex-ratio-kkt}
\end{equation}

If \(i\in S\), both logarithmic ratios in \eqref{eq:simplex-ratio-kkt} tend to zero. Since \(\eta_k\geq\underline\eta\),
\[
\nabla_i\Phi(x^\star)=\nabla_j\Phi(x^\star).
\]
Thus the gradient is constant on \(S\). Denote this common value by \(\gamma^\star\).

Assume that \(i\notin S\) and \(\nabla_i\Phi(x^\star)<\gamma^\star\). For all sufficiently large \(k\),
\[
\nabla_i\Phi(x_{k+1})-\nabla_j\Phi(x_{k+1})\leq-\varepsilon
\]
for some \(\varepsilon>0\). Since the logarithmic ratio associated with \(j\) tends to zero, \eqref{eq:simplex-ratio-kkt} yields, for all sufficiently large \(k\),
\[
\log\frac{x_{k+1,i}}{x_{k,i}}
\geq \frac{\underline\eta\varepsilon}{2}>0.
\]
Hence \(x_{k+1,i}\geq qx_{k,i}\) for some \(q>1\) and all sufficiently large \(k\). This is impossible because \(x_{k,i}>0\) and \(x_{k,i}\to0\). Therefore,
\[
\nabla_i\Phi(x^\star)\geq\gamma^\star
\qquad\text{for every }i\notin S.
\]

Since the gradient is equal to \(\gamma^\star\) on \(S\), we also have
\[
\gamma^\star=x^{\star\top}\nabla\Phi(x^\star)=\tau^\star.
\]
The KKT conditions now follow exactly as in Proposition~\ref{prop:kkt-simplex-flow}.
\end{proof}

\subsection{Box Constraints}
\label{subsec:kkt-box}

We next consider
\[
\min_{\ell\leq x\leq b}\Phi(x),
\qquad \ell_i<b_i \quad\text{for all }i.
\]
We use the logistic representation
\[
x_i=\ell_i+(b_i-\ell_i)\sigma(w_i),
\qquad
\sigma(t)=\frac{1}{1+e^{-t}}.
\]
The parameter dynamics is
\begin{equation}
\dot w=-\nabla\Phi(x),
\qquad
x=g(w),
\label{eq:box-w-flow-kkt}
\end{equation}
and the primal dynamics is
\begin{equation}
\dot x_i=-m_i(x)\nabla_i\Phi(x),
\qquad
m_i(x)=\frac{(x_i-\ell_i)(b_i-x_i)}{b_i-\ell_i}.
\label{eq:box-flow-kkt}
\end{equation}
The mobility is positive in the interior and vanishes at both endpoints.

\begin{proposition}[Convergent box trajectories satisfy the KKT conditions]
\label{prop:kkt-box-flow}
Let \(x(t)=g(w(t))\) solve \eqref{eq:box-w-flow-kkt}. Assume that
\[
x(t)\longrightarrow x^\star\in[\ell,b].
\]
Then, for every \(i\),
\[
\ell_i<x_i^\star<b_i
\quad\Longrightarrow\quad
\nabla_i\Phi(x^\star)=0,
\]
\[
x_i^\star=\ell_i
\quad\Longrightarrow\quad
\nabla_i\Phi(x^\star)\geq0,
\]
and
\[
x_i^\star=b_i
\quad\Longrightarrow\quad
\nabla_i\Phi(x^\star)\leq0.
\]
Hence \(x^\star\) satisfies the KKT conditions of the box-constrained problem.
\end{proposition}

\begin{proof}
Fix an index \(i\). If \(\ell_i<x_i^\star<b_i\), then \(w_i(t)\) converges to the finite value \(g_i^{-1}(x_i^\star)\). As in the orthant case, \(\dot w_i(t)=-\nabla_i\Phi(x(t))\) can converge to zero only, and therefore \(\nabla_i\Phi(x^\star)=0\).

If \(x_i^\star=\ell_i\), then \(w_i(t)\to-\infty\). If \(\nabla_i\Phi(x^\star)<0\), then \(\dot w_i(t)\geq\varepsilon\) for all sufficiently large \(t\), for some \(\varepsilon>0\). This contradicts \(w_i(t)\to-\infty\). Hence \(\nabla_i\Phi(x^\star)\geq0\).

If \(x_i^\star=b_i\), then \(w_i(t)\to+\infty\). If \(\nabla_i\Phi(x^\star)>0\), then \(\dot w_i(t)\leq-\varepsilon\) for all sufficiently large \(t\). This contradicts \(w_i(t)\to+\infty\). Hence \(\nabla_i\Phi(x^\star)\leq0\).

Define
\[
\alpha_i^\star=
\begin{cases}
\nabla_i\Phi(x^\star),&x_i^\star=\ell_i,\\
0,&x_i^\star>\ell_i,
\end{cases}
\qquad
\beta_i^\star=
\begin{cases}
-\nabla_i\Phi(x^\star),&x_i^\star=b_i,\\
0,&x_i^\star<b_i.
\end{cases}
\]
Then
\[
\nabla\Phi(x^\star)-\alpha^\star+\beta^\star=0,
\qquad
\alpha^\star\geq0,
\qquad
\beta^\star\geq0,
\]
and
\[
\alpha_i^\star(x_i^\star-\ell_i)=0,
\qquad
\beta_i^\star(b_i-x_i^\star)=0.
\]
This is the box KKT system.
\end{proof}

\begin{proposition}[Convergent implicit box iterates satisfy the KKT conditions]
\label{prop:kkt-box-discrete}
Let \((w_k)\) be generated by
\begin{equation}
 w_{k+1}-w_k+\eta_k\nabla\Phi(x_{k+1})=0,
 \qquad x_{k+1}=g(w_{k+1}),
\label{eq:box-implicit-kkt}
\end{equation}
where \(\eta_k\geq\underline\eta>0\). If
\[
x_k\longrightarrow x^\star\in[\ell,b],
\]
then \(x^\star\) satisfies the box KKT conditions.
\end{proposition}

\begin{proof}
Fix an index \(i\). If \(\ell_i<x_i^\star<b_i\), then \(w_{k,i}\) converges to a finite value, and \(w_{k+1,i}-w_{k,i}\to0\). Equation~\eqref{eq:box-implicit-kkt} and the lower bound on \(\eta_k\) give \(\nabla_i\Phi(x^\star)=0\).

If \(x_i^\star=\ell_i\), then \(w_{k,i}\to-\infty\). A negative value of \(\nabla_i\Phi(x^\star)\) would imply
\[
w_{k+1,i}-w_{k,i}\geq\underline\eta\varepsilon>0
\]
for all sufficiently large \(k\), which is impossible. Thus \(\nabla_i\Phi(x^\star)\geq0\).

If \(x_i^\star=b_i\), then \(w_{k,i}\to+\infty\). A positive value of \(\nabla_i\Phi(x^\star)\) would imply
\[
w_{k+1,i}-w_{k,i}\leq-\underline\eta\varepsilon<0
\]
for all sufficiently large \(k\), which is also impossible. Thus \(\nabla_i\Phi(x^\star)\leq0\). The KKT conditions follow from Proposition~\ref{prop:kkt-box-flow}.
\end{proof}

\subsection{Stiefel Manifold}
\label{subsec:kkt-stiefel}

We finally consider
\[
\min_{X\in\St(n,p)}\Phi(X).
\]
Let \(G(X)=\nabla\Phi(X)\) be the Euclidean gradient. With the canonical metric on \(\St(n,p)\), the Riemannian gradient is
\begin{equation}
\grad\Phi(X)
=G(X)-XG(X)^\top X
=A(X)X,
\label{eq:canonical-stiefel-gradient-kkt}
\end{equation}
where
\[
A(X):=G(X)X^\top-XG(X)^\top
\]
is skew-symmetric. Thus the Stiefel flow used in this paper is
\[
\dot X=-\grad\Phi(X)=-A(X)X.
\]

\begin{proposition}[Stationarity and first-order optimality on \(\St(n,p)\)]
\label{prop:kkt-stiefel-flow}
For \(X^\star\in\St(n,p)\), the following statements are equivalent:
\begin{enumerate}[label=\textup{(\roman*)},itemsep=2pt,topsep=3pt]
\item \(\grad\Phi(X^\star)=0\);
\item \(A(X^\star)X^\star=0\);
\item there exists a symmetric matrix \(\Lambda^\star\in\mathbb R^{p\times p}\) such that
\[
\nabla\Phi(X^\star)=X^\star\Lambda^\star.
\]
\end{enumerate}
Hence Riemannian stationarity is equivalent to the KKT condition associated with \(X^\top X=I_p\).
\end{proposition}

\begin{proof}
The equivalence of (i) and (ii) follows from \eqref{eq:canonical-stiefel-gradient-kkt}. Assume (ii), and write \(G^\star=G(X^\star)\). Then
\[
G^\star-X^\star G^{\star\top}X^\star=0.
\]
Multiplying on the left by \(X^{\star\top}\) gives
\[
X^{\star\top}G^\star=G^{\star\top}X^\star,
\]
so \(X^{\star\top}G^\star\) is symmetric. Setting
\[
\Lambda^\star:=X^{\star\top}G^\star
\]
gives \(G^\star=X^\star\Lambda^\star\), and therefore (iii) holds.

Conversely, if \(G^\star=X^\star\Lambda^\star\) with \(\Lambda^\star\) symmetric, then
\[
A(X^\star)X^\star
=X^\star\Lambda^\star-X^\star\Lambda^{\star\top}=0.
\]
Thus (iii) implies (ii).
\end{proof}

As a direct consequence, if a solution of the Riemannian gradient flow converges to \(X^\star\), then \(X^\star\) is a first-order critical point. Indeed, by continuity, a nonzero limiting Riemannian gradient would prevent the trajectory from converging.

\begin{proposition}[Convergent implicit Cayley iterates satisfy the first-order condition]
\label{prop:kkt-stiefel-discrete}
Let \((X_k)\subset\St(n,p)\) satisfy
\begin{equation}
\left(I+\frac{\eta_k}{2}A(X_{k+1})\right)X_{k+1}
=
\left(I-\frac{\eta_k}{2}A(X_{k+1})\right)X_k,
\label{eq:stiefel-cayley-kkt}
\end{equation}
where \(\eta_k\geq\underline\eta>0\). If
\[
X_k\longrightarrow X^\star\in\St(n,p),
\]
then
\[
\grad\Phi(X^\star)=0.
\]
\end{proposition}

\begin{proof}
Equation~\eqref{eq:stiefel-cayley-kkt} gives
\[
\frac{X_{k+1}-X_k}{\eta_k}
=-\frac12 A(X_{k+1})(X_{k+1}+X_k).
\]
Since \(X_k\to X^\star\), we have \(X_{k+1}-X_k\to0\). The lower bound on \(\eta_k\) implies that the left-hand side tends to zero. Passing to the limit gives
\[
A(X^\star)X^\star=0.
\]
Proposition~\ref{prop:kkt-stiefel-flow} now yields \(\grad\Phi(X^\star)=0\).
\end{proof}

\subsection{Summary}
\label{subsec:kkt-summary}

The same conclusion is obtained in the four settings, but the boundary argument is different from the manifold argument. For the orthant and the box, a wrong sign of the gradient would move the parameter away from the active face. On the simplex, the logarithmic dynamics gives the same conclusion and also shows that the gradient is constant on the support. On the Stiefel manifold, there is no boundary, and stationarity is exactly the usual Riemannian first-order condition.

The main point is that the vanishing of the primal velocity alone is not sufficient when the mobility degenerates at the boundary. The KKT inequalities follow from the asymptotic behavior of the geometry-generating variables. This is the precise form of the connection between the \gravidy\ dynamics and first-order optimality.

\section{Algorithms: Implicit Schemes and Inner Solvers}
\label{sec:algorithms}

We now describe how the implicit equations introduced in the previous sections are solved in practice. The outer step and the inner solve must be distinguished. The outer step is fixed by the geometry: it is a Bregman proximal equation on the vector domains and a Cayley equation on the Stiefel manifold. The inner method only computes an approximate solution of this equation.

Throughout this section, the inner iterations are stopped when the residual of the implicit equation is sufficiently small. The convergence analysis in Section~\ref{sec:global-convergence} will therefore have to account for the remaining inner error. We first treat the orthant and box together, since they lead to the same parameter residual. We then discuss the simplex and the Stiefel manifold.

\subsection{Orthant and Box: Backward Euler in the Parameter Variable}
\label{subsec:orthant-box-inner}

Let $x=g(\zeta)$, where $\zeta=u$ for the orthant and $\zeta=w$ for the box. At outer iteration $k$, backward Euler gives
\begin{equation}
F_k(\zeta)
:=\zeta-\zeta_k+\eta_k\nabla\Phi(g(\zeta))=0.
\label{eq:parameter-inner-residual}
\end{equation}
As shown in Section~\ref{sec:riemannian-lens}, this is also the optimality equation of the corresponding Bregman proximal step in the primal variable. If $\Phi\in C^2$, then
\begin{equation}
J_k(\zeta)
:=DF_k(\zeta)
=I+\eta_k\nabla^2\Phi(g(\zeta))D(\zeta),
\qquad
D(\zeta):=J_g(\zeta).
\label{eq:parameter-inner-jacobian}
\end{equation}
For the componentwise maps used here, $D(\zeta)$ is diagonal and positive for every finite $\zeta$.

For a convex objective, $H(\zeta):=\nabla^2\Phi(g(\zeta))\succeq0$. In this case,
\begin{equation}
J_k(\zeta)
=D(\zeta)^{-1/2}
\left(I+\eta_kD(\zeta)^{1/2}H(\zeta)D(\zeta)^{1/2}\right)
D(\zeta)^{1/2}.
\label{eq:parameter-jacobian-similarity}
\end{equation}
Hence $J_k(\zeta)$ is nonsingular. This similarity relation does not imply that $J_k(\zeta)$ is well conditioned in the Euclidean norm. In particular, the scaling may become unbalanced when some entries of $D(\zeta)$ are very small.

\paragraph{Newton direction through a symmetric positive definite system.}
The Newton equation
\[
J_k(\zeta)h=-F_k(\zeta)
\]
can be solved without forming the nonsymmetric matrix $J_k(\zeta)$. Let
\[
h=D(\zeta)^{-1/2}v.
\]
Multiplying the Newton equation by $D(\zeta)^{1/2}$ gives
\begin{equation}
\left(I+\eta_kD(\zeta)^{1/2}H(\zeta)D(\zeta)^{1/2}\right)v
=-D(\zeta)^{1/2}F_k(\zeta).
\label{eq:parameter-newton-spd}
\end{equation}
For convex $\Phi$, the matrix in \eqref{eq:parameter-newton-spd} is symmetric positive definite. It can be solved by Cholesky in a dense implementation, or by conjugate gradients when Hessian--vector products are available. Recovering $h$ from $v$ requires division by $D(\zeta)^{1/2}$, so care is needed when the iterate is very close to an active face.

\paragraph{Levenberg--Marquardt regularization.}
A second possibility is to minimize the residual merit function
\[
m_k(\zeta)=\frac12\|F_k(\zeta)\|^2.
\]
The Levenberg--Marquardt, or regularized Gauss--Newton, direction solves
\begin{equation}
\left(J_k(\zeta)^\top J_k(\zeta)+\lambda I\right)h
=-J_k(\zeta)^\top F_k(\zeta),
\qquad \lambda>0.
\label{eq:parameter-lm}
\end{equation}
The coefficient matrix is symmetric positive definite for every $\lambda>0$. The regularization improves robustness when the linearized equation is poorly scaled, although the use of normal equations may square the condition number of $J_k(\zeta)$. In our implementation, \eqref{eq:parameter-lm} is applied matrix-free and solved by conjugate gradients. A backtracking step ensures that $m_k$ decreases. The damping parameter may be kept fixed or adapted according to the success of the trial step.

\begin{algorithm}[H]
\DontPrintSemicolon
\caption{Parameter backward Euler with a regularized Gauss--Newton inner solve}
\label{alg:MGN-orthant-box}
\KwIn{$\zeta_k$, $\eta_k>0$, an initial damping $\lambda>0$, and an inner tolerance $\varepsilon_k$}
$\zeta\leftarrow\zeta_k$\;
\For{$j=0,1,\ldots$}{
  $x\leftarrow g(\zeta)$; $F\leftarrow \zeta-\zeta_k+\eta_k\nabla\Phi(x)$\;
  \lIf{$\|F\|\leq\varepsilon_k$}{\textbf{break}}
  Form the actions of $J=I+\eta_k\nabla^2\Phi(x)J_g(\zeta)$ and $J^\top$\;
  Solve $(J^\top J+\lambda I)h=-J^\top F$ approximately\;
  Choose $\alpha\in(0,1]$ by backtracking until
  $\|F_k(\zeta+\alpha h)\|<\|F_k(\zeta)\|$\;
  $\zeta\leftarrow\zeta+\alpha h$\;
  Optionally decrease $\lambda$ after a successful full step and increase it after a rejected trial\;
}
\KwOut{$\zeta_{k+1}\leftarrow\zeta$, $x_{k+1}\leftarrow g(\zeta_{k+1})$}
\end{algorithm}

\paragraph{The two maps used in the experiments.}
For the orthant,
\[
g_i(u_i)=e^{u_i},
\qquad
D(u)=\diag(x).
\]
For the box,
\[
g_i(w_i)=\ell_i+(b_i-\ell_i)\sigma(w_i),
\qquad
D(w)=\diag\!\left(
\frac{(x_i-\ell_i)(b_i-x_i)}{b_i-\ell_i}
\right)_{i=1}^n.
\]
Both maps keep every finite iterate in the interior. The mobility tends to zero when a component approaches an active face.

\begin{remark}[Exact and inexact inner solves]
If $F_k(\zeta_{k+1})=0$, the resulting primal point is the exact Bregman proximal step. In practice, the residual is only reduced below a tolerance $\varepsilon_k$. The sequence $(\varepsilon_k)$ must be chosen consistently with the inexact convergence conditions stated later. Neither a small Newton correction nor a small gradient norm alone is a substitute for the residual test in \eqref{eq:parameter-inner-residual}.
\end{remark}

\begin{remark}
    For the dense least-squares experiments, we cache
\(H=A^\top A\) and \(c=A^\top b\). The transformed Newton system is solved
by preconditioned conjugate gradients with the Jacobi preconditioner
\[
P_{ii}=1+\eta d_iH_{ii},
\]
where \(d_i=x_i\) for the orthant and \(d_i=g_i'(z_i)\) for the box.
For the regularized Gauss--Newton system, we use the exact diagonal
preconditioner
\[
P_{ii}
=
1+2\eta d_iH_{ii}
+\eta^2d_i^2\sum_j H_{ij}^2+\lambda.
\]
The original matrix-free actions based on \(A\) and \(A^\top\) remain
available for problems where forming \(H\) is undesirable.
\end{remark}

\subsection{Simplex: KL-Proximal Step and Newton--KKT Solver}
\label{sec:simplex-klprox}

Let $x_k\in\operatorname{ri}(\Delta_n)$. The implicit entropy step is
\begin{equation}
 x_{k+1}
 \in\argmin_{x\in\Delta_n}
 \left\{
 \KL(x\|x_k)+\eta_k\Phi(x)
 \right\},
 \qquad
 \KL(x\|x_k)=\sum_{i=1}^n x_i\log\frac{x_i}{x_{k,i}}.
\label{eq:kl-prox-prob}
\end{equation}
The omitted affine terms in the entropy Bregman divergence vanish because both vectors have unit mass. If $\Phi$ is continuous and convex on $\Delta_n$, the problem has a unique minimizer. The minimizer belongs to $\operatorname{ri}(\Delta_n)$ when $x_k$ is strictly positive and $\Phi$ has a finite gradient on the simplex. The entropy term is not a barrier in value, but its derivative tends to $-\infty$ when a component approaches zero.

Let
\[
R_k(x)=\KL(x\|x_k)+\eta_k\Phi(x).
\]
The equality-constrained first-order system is
\begin{equation}
\log x-\log x_k+\eta_k\nabla\Phi(x)+\nu\mathbf 1=0,
\qquad
\mathbf 1^\top x=1.
\label{eq:kl-implicit-kkt}
\end{equation}
The constant vector arising from differentiating $\sum_i x_i\log x_i$ is absorbed into $\nu$. No nonnegativity multiplier is needed at a finite inner iterate because the solution is strictly positive.

Assume first that $\Phi$ is convex and twice differentiable. Set
\begin{equation}
q_k(x)=\log x-\log x_k+\eta_k\nabla\Phi(x),
\qquad
K_k(x)=\diag(1/x)+\eta_k\nabla^2\Phi(x).
\label{eq:simplex-newton-objects}
\end{equation}
Then $K_k(x)\succ0$ on the simplex interior. Its smallest eigenvalue is bounded away from zero, but its largest eigenvalue may become large when a component of $x$ is small. Thus positive definiteness should not be confused with uniform good conditioning.

At a feasible inner iterate, the Newton correction solves
\begin{equation}
\begin{bmatrix}
K_k(x)&\mathbf 1\\
\mathbf 1^\top&0
\end{bmatrix}
\begin{bmatrix}
\Delta x\\ \Delta\nu
\end{bmatrix}
=-
\begin{bmatrix}
q_k(x)\\0
\end{bmatrix}.
\label{eq:simplex-newton-kkt-system}
\end{equation}
The current value of $\nu$ need not be stored, since adding a constant multiple of $\mathbf 1$ to $q_k(x)$ does not change $\Delta x$. Let
\[
y=K_k(x)^{-1}q_k(x),
\qquad
z=K_k(x)^{-1}\mathbf 1.
\]
Then
\begin{equation}
\Delta\nu=-\frac{\mathbf 1^\top y}{\mathbf 1^\top z},
\qquad
\Delta x=-y-z\Delta\nu,
\label{eq:simplex-schur-step}
\end{equation}
and $\mathbf 1^\top\Delta x=0$. A fraction-to-the-boundary rule preserves strict positivity, and a line search on $R_k$ globalizes the inner Newton iteration.

\begin{algorithm}[H]
\DontPrintSemicolon
\caption{Newton--KKT inner solver for the simplex KL-proximal step}
\label{alg:KLprox}
\KwIn{$x_k\in\operatorname{ri}(\Delta_n)$, $\eta_k>0$, and an inner tolerance $\varepsilon_k$}
$x\leftarrow x_k$\;
\For{$j=0,1,\ldots$}{
  Compute $q\leftarrow\log x-\log x_k+\eta_k\nabla\Phi(x)$\;
  Set $\widehat\nu\leftarrow-\frac1n\mathbf 1^\top q$\;
  \lIf{$\|q+\widehat\nu\mathbf 1\|\leq\varepsilon_k$}{\textbf{break}}
  $K\leftarrow\diag(1/x)+\eta_k\nabla^2\Phi(x)$\;
  Solve $Ky=q$ and $Kz=\mathbf 1$\;
  $\Delta\nu\leftarrow-(\mathbf 1^\top y)/(\mathbf 1^\top z)$;
  $\Delta x\leftarrow-y-z\Delta\nu$\;
  Choose $\alpha\in(0,1]$ by a fraction-to-the-boundary rule and backtracking on $R_k$\;
  $x\leftarrow x+\alpha\Delta x$\;
}
\KwOut{$x_{k+1}\leftarrow x$}
\end{algorithm}

For least squares, $\nabla^2\Phi=A^\top A$, so the two systems with $K_k(x)$ can be solved by Cholesky in the dense case. Diagonal scaling or preconditioned conjugate gradients may be preferable when $x$ is close to the boundary or when $A$ is sparse. A failed Cholesky factorization in the convex least-squares setting indicates numerical loss of definiteness rather than a failure of the mathematical model.

\paragraph{Hessian-free fixed-point variant.}
Appendix~\ref{app:kl-fixed-point} gives a relaxed Picard iteration for \eqref{eq:kl-implicit-kkt}. Its search direction is a descent direction for $R_k$ unless the fixed-point equation is already satisfied. This provides a useful Hessian-free alternative, although we do not claim the same local rate as Newton--KKT.

\paragraph{Reduced-logit formulation.}
A second alternative works in the quotient coordinates
\[
v_i=u_i-u_n,
\qquad i=1,\ldots,n-1,
\]
with $x(v)=\softmax([v;0])$. Define
\[
P=\begin{bmatrix}I_{n-1}&-\mathbf 1_{n-1}\end{bmatrix}
\in\mathbb R^{(n-1)\times n}.
\]
Subtracting the $n$th equation of \eqref{eq:kl-implicit-kkt} from the first $n-1$ equations gives the gauge-invariant residual
\begin{equation}
F_k^{\Delta}(v)
=v-v_k+\eta_kP\nabla\Phi(x(v)).
\label{eq:reduced-logit-residual}
\end{equation}
Its Jacobian is
\begin{equation}
J_k^{\Delta}(v)
=I_{n-1}
+\eta_kP\nabla^2\Phi(x(v))J_v(v),
\qquad
J_v(v)=J_{\softmax}([v;0])_{:,1:n-1}.
\label{eq:reduced-logit-jacobian}
\end{equation}
A regularized Gauss--Newton solve may then be applied to $F_k^{\Delta}(v)=0$. The subtraction by the reference component is essential. Without it, the residual depends on the chosen softmax gauge and is not equivalent to the KL-proximal equation. The complete variant is given in Appendix~\ref{app:kl-mgn-variant}.

\subsection{Stiefel Manifold: Implicit Cayley Equation}
\label{sec:stiefel-ics}

Let
\[
G(X)=\nabla\Phi(X),
\qquad
A(X)=G(X)X^\top-XG(X)^\top.
\]
Under the canonical metric used in this paper,
\[
\grad\Phi(X)=A(X)X.
\]
Given $X_k\in\St(n,p)$ and $c_k=\eta_k/2$, the implicit Cayley equation is
\begin{equation}
F_k(Y)
:=\left(I+c_kA(Y)\right)Y
-\left(I-c_kA(Y)\right)X_k
=0.
\label{eq:ics-fixedpoint}
\end{equation}
If $Y$ is an exact root, then
\[
Y=Q(Y)X_k,
\qquad
Q(Y)=\left(I+c_kA(Y)\right)^{-1}
       \left(I-c_kA(Y)\right).
\]
Since $A(Y)$ is skew-symmetric, $Q(Y)$ is orthogonal. Hence every exact root of \eqref{eq:ics-fixedpoint} belongs to $\St(n,p)$.

When the skew field is frozen, \eqref{eq:ics-fixedpoint} is the trapezoidal rule for the linear system $\dot X=-AX$. The corresponding Cayley factor is norm preserving for every $\eta_k>0$. For a state-dependent field, however, this observation is only a stability motivation. It does not imply that the nonlinear equation has a solution for every stepsize, or that every computed root decreases $\Phi$. We therefore combine the inner solve with an outer acceptance test.

\paragraph{Fr\'echet derivative.}
Assume that $G$ is differentiable and denote its directional derivative by
\[
\mathcal H_Y[H]=DG(Y)[H].
\]
Then
\begin{align}
DA(Y)[H]
={}&\mathcal H_Y[H]Y^\top+G(Y)H^\top
-HG(Y)^\top-Y\mathcal H_Y[H]^\top,
\label{eq:stiefel-dA-general}\\
DF_k(Y)[H]
={}&\left(I+c_kA(Y)\right)H
+c_kDA(Y)[H](Y+X_k).
\label{eq:stiefel-dF-general}
\end{align}
For the quadratic model used in the experiments, $G$ is linear and $\mathcal H_Y[H]=G(H)$. Appendix~\ref{app:stiefel-frechet} gives the resulting matrix-free action and its vectorized form.

\paragraph{Newton--Krylov inner solve.}
At an inner iterate $Y$, Newton's equation is
\begin{equation}
DF_k(Y)[H]=-F_k(Y).
\label{eq:stiefel-newton-equation}
\end{equation}
GMRES can solve this equation using only the action \eqref{eq:stiefel-dF-general}. The dominant additional operation is a Hessian--vector product $DG(Y)[H]$. For the quadratic model, this reduces to applying the matrices $Q^{(j)}$ columnwise.

The factor $I+c_kA(Y)$ provides a natural left preconditioner. Since $A(Y)=UV^\top-VU^\top$ with $U=G(Y)$ and $V=Y$, its inverse can be applied through a Woodbury formula using a $2p\times2p$ system. The cost is $O(np^2+p^3)$, in addition to the derivative action.

\paragraph{Dense Newton solve.}
When $np$ is moderate, the Jacobian of $F_k$ can be assembled by applying \eqref{eq:stiefel-dF-general} to the canonical basis of $\mathbb R^{n\times p}$. The resulting $np\times np$ system is then solved by LU. This is useful for validation and for small dense problems, but its memory and factorization costs limit its range of application.

\paragraph{Feasibility of approximate inner iterates.}
Exact roots are feasible, but an ambient Newton iterate need not be. Two
implementations are possible. The first performs the Newton iterations in the
ambient space and returns a point only after the residual is sufficiently
small. The second starts from a feasible Cayley predictor and applies a polar
or QR projection after each trial correction. This keeps the inner iterates
feasible, but it changes the exact Newton map. We therefore treat it as a
safeguarded projected Newton variant and do not claim an automatic quadratic
rate for it. In either case, one final polar projection is applied before the
outer acceptance test, after which the implicit residual and the objective are
recomputed. Thus every accepted outer iterate is feasible, while the ambient
Newton corrections themselves remain unchanged. Appendix~\ref{app:f-ics}
describes the feasible predictor and the projected safeguard.

\paragraph{Outer acceptance.}
Let $Y_k$ be the point returned by the inner solver. We accept it only if its residual is below the prescribed tolerance and
\begin{equation}
\Phi(Y_k)
\leq
\Phi(X_k)-c_1\eta_k\|\grad\Phi(X_k)\|_F^2,
\qquad c_1\in(0,\tfrac12).
\label{eq:stiefel-acceptance}
\end{equation}
If either test fails, the step is rejected and recomputed with a smaller $\eta_k$. After a successful step, $\eta_k$ may be increased. The inequality \eqref{eq:stiefel-acceptance} is an algorithmic safeguard; it is not a consequence of A-stability.

\begin{algorithm}[H]
\DontPrintSemicolon
\caption{Implicit Cayley step with a Newton--Krylov inner solve}
\label{alg:ICS}
\KwIn{$X_k\in\St(n,p)$, $\eta_k>0$, residual tolerance $\varepsilon_k$, and $c_1\in(0,\tfrac12)$}
\Repeat{the residual and decrease tests are satisfied}{
  Choose an initial guess $Y$, for example the feasible Cayley predictor of Appendix~\ref{app:f-ics}\;
  \For{$j=0,1,\ldots$}{
    Compute $F_k(Y)$\;
    \lIf{$\|F_k(Y)\|_F\leq\varepsilon_k$}{\textbf{break}}
    Solve $DF_k(Y)[H]=-F_k(Y)$ approximately by GMRES\;
    Use a residual line search for $Y+H$; optionally apply the feasibility safeguard described in Appendix~\ref{app:f-ics}\;
    After any projection, recompute the implicit residual and the objective\;
  }
  Set $Y_k\leftarrow\operatorname{Proj}_{\St}(Y)$, and recompute $F_k(Y_k)$ and $\Phi(Y_k)$\;
  \eIf{$\|F_k(Y_k)\|_F\leq\varepsilon_k$ and \eqref{eq:stiefel-acceptance} holds}{
    accept $X_{k+1}\leftarrow Y_k$ and optionally increase $\eta_k$\;
  }{
    decrease $\eta_k$ and repeat the inner solve\;
  }
}
\KwOut{$X_{k+1}$}
\end{algorithm}

Appendix~\ref{app:stiefel-frechet} gives the exact derivative used in the Newton equation. Appendix~\ref{app:f-ics} describes the feasible predictor and the projected safeguard. Appendix~\ref{app:dense-nr} gives the dense Newton alternative.

\subsection{Implementation Consistency}
\label{subsec:implementation-consistency}

The algorithms above describe the equations that must be solved. Three details are particularly important when implementing them.

First, the stopping criterion of an outer constrained method must measure constrained stationarity. The Euclidean norm $\|\nabla\Phi(x)\|$ is not suitable when the solution lies on the boundary or, on the simplex, when the optimal gradient is a nonzero constant vector. The experiments therefore use the projected KKT residual on the vector
domains and $\|\grad\Phi(X)\|_F$ on the Stiefel manifold.

Second, the reduced-logit simplex residual must contain the gradient differences in \eqref{eq:reduced-logit-residual}. This makes the equation independent of the softmax gauge.

Third, the Stiefel Newton--Krylov and dense Newton variants must use the
residual derivative at the current inner point $Y$, namely
\eqref{eq:stiefel-dF-general}, with $DA(Y)[H]$ evaluated as in
\eqref{eq:stiefel-dA-general}.
Freezing the derivative terms at $X_k$ produces a quasi-Newton approximation, not the exact Fr\'echet derivative. Such an approximation may be useful in practice, but it must be identified as such and analyzed separately.

\section{Global Convergence of Implicit Gradient--Flow Discretizations}
\label{sec:global-convergence}

We now study the convergence of the implicit schemes on the orthant, the simplex, and the box. By Sections~\ref{sec:jacobian-nonpullback} and \ref{sec:gravidy}, these three methods are exact Bregman proximal steps in the primal variable. This common formulation gives a short analysis that applies to the three geometries.

The Stiefel Cayley step is not a Bregman proximal method, and it is not covered by the results below. Its feasibility and inner Newton equations were studied in Section~\ref{sec:algorithms}. A separate manifold convergence analysis would require additional conditions on the accepted Cayley steps.

We first treat exact outer steps. We then state precisely what is needed when the nonlinear subproblem is solved only approximately. The convex results allow arbitrary positive stepsizes. In the nonconvex case, additional assumptions are needed: the selected point must satisfy a genuine descent condition, the stepsizes must remain bounded away from zero and infinity, and the iterates must stay in a compact subset of the interior. These restrictions will be made explicit.

All proofs are given in Appendix~\ref{app:proofs-global}.

\subsection{Common Bregman Formulation}
\label{subsec:setup}

Let \(\mathcal X\) be a closed convex set contained in a
finite-dimensional affine space, and let
\[
\mathcal C:=\operatorname{ri}(\mathcal X).
\]
All derivatives in this section are taken along directions that remain in
\(\operatorname{aff}(\mathcal X)\), and we use the usual Euclidean inner
product on these direction vectors. The symbols \(\nabla f\) and
\(\nabla h\) denote the corresponding gradients.

For the orthant and box constraints,
\(\operatorname{aff}(\mathcal X)=\mathbb R^n\), so these are the usual
Euclidean gradients. On the simplex, the admissible directions satisfy
\(\mathbf1^\top d=0\), and therefore
\[
\nabla f(x)=\Pi\nabla\Phi(x)\in\mathbf1^\perp.
\]
In particular, the condition \(\nabla f(x)=0\) on the simplex means that
\(\nabla\Phi(x)\) is a constant multiple of \(\mathbf1\), rather than that
the ambient gradient vanishes.

Let \(h:\mathcal C\to\mathbb R\) be a \(C^2\) Legendre function. Whenever
needed, we assume that \(h\) has a continuous extension to the comparison
points in \(\mathcal X\). For \(x\in\mathcal X\) and \(y\in\mathcal C\),
define
\begin{equation}
D_h(x,y)
:=h(x)-h(y)-\langle \nabla h(y),x-y\rangle.
\label{eq:global-bregman-distance}
\end{equation}

Starting from $x_k\in\mathcal C$, the exact implicit step is
\begin{equation}
 x_{k+1}\in\argmin_{x\in\mathcal X}
 \left\{f(x)+\frac{1}{\eta_k}D_h(x,x_k)\right\},
 \qquad \eta_k>0.
\label{eq:global-bregman-step}
\end{equation}
We assume that the minimizer exists and belongs to $\mathcal C$. For the
compact simplex and box problems, existence and interiority follow under the
stated continuity, smoothness, and entropy assumptions. On the unbounded
orthant, and in the nonconvex setting, existence must be assumed or ensured
separately, for example by coercivity of the Bregman subproblem. If $f$ is
convex, the minimizer is unique because $D_h(\cdot,x_k)$ is strictly convex.

The first-order condition of \eqref{eq:global-bregman-step} is
\begin{equation}
\nabla h(x_{k+1})-\nabla h(x_k)+\eta_k\nabla f(x_{k+1})=0.
\label{eq:global-dual-step}
\end{equation}
This is exactly the backward-Euler equation in the dual variable $\nabla h(x)$.

We will repeatedly use the three-point identity
\begin{equation}
\langle \nabla h(y)-\nabla h(x),u-y\rangle
=
D_h(u,x)-D_h(u,y)-D_h(y,x).
\label{eq:three-point}
\end{equation}

\subsection{Convex Objectives}
\label{subsec:convex-global}

We first assume that $f$ is convex on $\mathcal X$.

\begin{proposition}[One-step Bregman inequality]
\label{prop:one-step}
Let $x_{k+1}$ satisfy \eqref{eq:global-bregman-step}. Then, for every $x\in\mathcal X$ for which $D_h(x,x_k)$ is finite,
\begin{equation}
\eta_k\bigl(f(x_{k+1})-f(x)\bigr)
\le
D_h(x,x_k)-D_h(x,x_{k+1})-D_h(x_{k+1},x_k).
\label{eq:one-step}
\end{equation}
In particular,
\begin{equation}
 f(x_{k+1})+\frac{1}{\eta_k}D_h(x_{k+1},x_k)
 \le f(x_k),
\label{eq:convex-descent}
\end{equation}
so the objective values are nonincreasing for every $\eta_k>0$.
\end{proposition}

The next theorem gives both the convergence estimate and the role of the stepsizes.

\begin{theorem}[Convex convergence and last-iterate rate]
\label{thm:global-convex}
Assume that $f$ is convex and that $x^\star\in\argmin_{x\in\mathcal X}f(x)$ satisfies $D_h(x^\star,x_0)<\infty$. Let
\[
S_K:=\sum_{k=0}^{K}\eta_k.
\]
Then
\begin{equation}
\sum_{k=0}^{K}\eta_k\bigl(f(x_{k+1})-f^\star\bigr)
\le D_h(x^\star,x_0),
\label{eq:convex-weighted-sum}
\end{equation}
and, since the objective values are nonincreasing,
\begin{equation}
 f(x_{K+1})-f^\star
 \le \frac{D_h(x^\star,x_0)}{S_K}.
\label{eq:convex-last-iterate}
\end{equation}
Moreover, $D_h(x^\star,x_k)$ is nonincreasing. Consequently, if $S_K\to\infty$, then
\[
f(x_k)\to f^\star.
\]
Every cluster point at which $f$ is continuous is therefore a minimizer.
\end{theorem}

For a constant stepsize $\eta_k\equiv\eta>0$, \eqref{eq:convex-last-iterate} gives
\[
f(x_{K+1})-f^\star
\le \frac{D_h(x^\star,x_0)}{\eta(K+1)}.
\]
There is no upper stability bound on $\eta$ in this exact convex result. A large value of $\eta$ may make the inner problem harder to solve, but it does not invalidate the descent inequality.

\subsection{Relative Strong Convexity}
\label{subsec:relative-strong}

\begin{assumption}[Relative strong convexity]
\label{ass:rel-strong}
There exists $\mu>0$ such that, for every $x\in\mathcal C$ and $y\in\mathcal X$,
\begin{equation}
 f(y)
 \ge f(x)+\langle \nabla f(x),y-x\rangle+\mu D_h(y,x).
\label{eq:relative-strong-convexity}
\end{equation}
\end{assumption}

The orientation of the Bregman distance in \eqref{eq:relative-strong-convexity} is important. It allows the comparison point $y=x^\star$ to lie on the boundary, provided $D_h(x^\star,x)$ is finite.

\begin{theorem}[Linear contraction in Bregman distance]
\label{thm:linear}
Assume that $f$ satisfies Assumption~\ref{ass:rel-strong}. Then every exact step satisfies
\begin{equation}
(1+\eta_k\mu)D_h(x^\star,x_{k+1})
+D_h(x_{k+1},x_k)
\le D_h(x^\star,x_k).
\label{eq:relative-strong-contraction}
\end{equation}
Hence
\begin{equation}
D_h(x^\star,x_{K+1})
\le
\left(\prod_{k=0}^{K}\frac{1}{1+\eta_k\mu}\right)
D_h(x^\star,x_0).
\label{eq:relative-strong-product}
\end{equation}
For a constant stepsize $\eta_k\equiv\eta>0$, the contraction factor is
\[
q=\frac{1}{1+\eta\mu}\in(0,1).
\]
If $h$ is uniformly strongly convex on a convex region containing $x^\star$ and the iterates, then $x_k\to x^\star$ linearly in norm as well.
\end{theorem}

This result is stronger and simpler than a Euclidean comparison based on upper and lower Hessian bounds: no condition-number factor is needed in the Bregman-distance contraction. If $x^\star\in\mathcal C$, $f$ is also relatively smooth with constant $L_{\rm rel}$, and the two orientations of $D_h$ are uniformly comparable on the relevant set, then the objective gap is $R$-linearly bounded by the same product in \eqref{eq:relative-strong-product}.

\subsection{Inexact Inner Solves}
\label{subsec:inexact}

The algorithms in Section~\ref{sec:algorithms} solve the nonlinear equation only approximately. We therefore record the exact effect of a nonzero inner residual. Let the computed point $x_{k+1}\in\mathcal C$ satisfy
\begin{equation}
 r_{k+1}
 :=\nabla h(x_{k+1})-\nabla h(x_k)
 +\eta_k\nabla f(x_{k+1}).
\label{eq:dual-residual-inexact}
\end{equation}

\begin{proposition}[Fundamental inequality with an inner residual]
\label{prop:inexact}
If $f$ is convex, then, for every admissible comparison point $x$,
\begin{align}
\eta_k\bigl(f(x_{k+1})-f(x)\bigr)
&\le
D_h(x,x_k)-D_h(x,x_{k+1})-D_h(x_{k+1},x_k)
\nonumber\\
&\quad +\langle r_{k+1},x_{k+1}-x\rangle.
\label{eq:inexact-fundamental}
\end{align}
Assume, in addition, that $\|x_{k+1}-x^\star\|\le R$ for all $k$. Then
\begin{equation}
\sum_{k=0}^{K}\eta_k\bigl(f(x_{k+1})-f^\star\bigr)
\le
D_h(x^\star,x_0)
+R\sum_{k=0}^{K}\|r_{k+1}\|_*.
\label{eq:inexact-convex-sum}
\end{equation}
Consequently, if $\sum_k\|r_{k+1}\|_*<\infty$, the best weighted iterate retains the rate $O(1/S_K)$. If the implementation also enforces $f(x_{k+1})\le f(x_k)$, then the same bound holds for the last iterate.
\end{proposition}

A summable residual does not automatically preserve the exact linear factor in Theorem~\ref{thm:linear}. Under relative strong convexity, the correct recursion is
\begin{equation}
(1+\eta_k\mu)D_h(x^\star,x_{k+1})
\le
D_h(x^\star,x_k)
+\langle r_{k+1},x_{k+1}-x^\star\rangle.
\label{eq:inexact-linear-recursion}
\end{equation}
Thus a linear rate requires a correspondingly strong control of the residual, for example a geometrically decreasing tolerance. A fixed inner tolerance generally leads to a final accuracy floor rather than exact convergence.

\subsection{Nonconvex Objectives}
\label{subsec:nonconvex}

We now drop convexity. There is an important distinction with the convex case. The implicit equation \eqref{eq:global-dual-step} may have several roots, and a root found by Newton's method need not minimize the Bregman subproblem. The first result below therefore assumes that $x_{k+1}$ is a global minimizer of \eqref{eq:global-bregman-step}. We give a directly verifiable inexact alternative afterward.

Assume that the iterates remain in a compact convex set $\mathcal K\Subset\mathcal C$ on which
\begin{equation}
\sigma_h I\preceq\nabla^2h(x)\preceq L_h I,
\qquad x\in\mathcal K,
\label{eq:nonconvex-h-bounds}
\end{equation}
and that
\begin{equation}
0<\underline\eta\le\eta_k\le\overline\eta<\infty.
\label{eq:nonconvex-step-bounds}
\end{equation}
These two-sided stepsize bounds have different roles: the upper bound gives a uniform decrease in the step norm, while the lower bound converts a small step into a small gradient.

\begin{proposition}[Sufficient decrease and asymptotic stationarity]
\label{prop:noncvx-decrease}
Assume that $f$ is $C^1$, bounded below on $\mathcal K$, and that each $x_{k+1}$ is a global minimizer of \eqref{eq:global-bregman-step}. Then
\begin{equation}
 f(x_k)-f(x_{k+1})
 \ge \frac{1}{\eta_k}D_h(x_{k+1},x_k)
 \ge \frac{\sigma_h}{2\overline\eta}\|x_{k+1}-x_k\|^2.
\label{eq:nonconvex-sufficient-decrease}
\end{equation}
Consequently,
\[
\sum_{k=0}^{\infty}\|x_{k+1}-x_k\|^2<\infty,
\qquad
\|x_{k+1}-x_k\|\to0.
\]
Moreover,
\begin{equation}
\|\nabla f(x_{k+1})\|
\le \frac{L_h}{\underline\eta}\|x_{k+1}-x_k\|,
\label{eq:nonconvex-relative-error}
\end{equation}
so
\begin{equation}
\nabla f(x_k)\to0.
\label{eq:nonconvex-gradient-zero}
\end{equation}
In addition,
\begin{equation}
\sum_{k=0}^{\infty}\eta_k\|\nabla f(x_{k+1})\|^2
\le \frac{2L_h^2}{\sigma_h}\bigl(f(x_0)-\inf_{\mathcal K}f\bigr).
\label{eq:gradient-square-sum}
\end{equation}
\end{proposition}

\begin{remark}[Meaning of asymptotic stationarity]
\label{rem:relative-gradient-nonconvex}
The gradient in \eqref{eq:nonconvex-gradient-zero} is the gradient relative
to the affine hull of \(\mathcal X\), as specified in
Section~\ref{subsec:setup}. Therefore, on the simplex,
\[
\nabla f(x_k)\to0
\qquad\Longleftrightarrow\qquad
\Pi\nabla\Phi(x_k)\to0.
\]
This does not require the ambient gradient \(\nabla\Phi(x_k)\) to vanish.

Moreover, the assumption
\[
\mathcal K\Subset\mathcal C
\]
means that the iterates remain in a compact subset of the relative interior.
Consequently, this nonconvex result does not cover convergence to an active
inequality boundary. For the orthant and box, the limiting condition is the
usual equation \(\nabla\Phi(x^\star)=0\) because the limit is interior.
Boundary limits and their KKT interpretation are treated separately in
Section~5.
\end{remark}

\begin{theorem}[KL convergence and rates]
\label{thm:KL-conv}
Assume the conditions of Proposition~\ref{prop:noncvx-decrease}, and assume that $f$ has the Kurdyka--\L{}ojasiewicz property on $\mathcal K$. Then the sequence has finite length,
\[
\sum_{k=0}^{\infty}\|x_{k+1}-x_k\|<\infty,
\]
and converges to a critical point $x^\star\in\mathcal K$.

If the KL inequality at $x^\star$ holds with exponent $\theta\in[0,1)$, then the objective residual $s_k:=f(x_k)-f(x^\star)$ has the usual rates:
\begin{itemize}
\item $\theta=0$: finite termination;
\item $0<\theta\le\tfrac12$: linear convergence;
\item $\tfrac12<\theta<1$: $s_k=O\bigl(k^{-1/(2\theta-1)}\bigr)$.
\end{itemize}
\end{theorem}

The same KL proof applies to an inexact outer step when the implementation enforces two uniform conditions:
\begin{equation}
 f(x_k)-f(x_{k+1})\ge a\|x_{k+1}-x_k\|^2,
\qquad
 \|\nabla f(x_{k+1})\|\le b\|x_{k+1}-x_k\|,
\label{eq:accepted-inexact-conditions}
\end{equation}
for some constants $a,b>0$. The first condition is a sufficient-decrease test. The second follows, for example, if the normalized inner residual
\[
 e_{k+1}:=\nabla f(x_{k+1})
 +\frac{1}{\eta_k}\bigl(\nabla h(x_{k+1})-\nabla h(x_k)\bigr)
\]
satisfies $\|e_{k+1}\|\le\gamma\|x_{k+1}-x_k\|$ with a fixed $\gamma$. These conditions are stronger than merely requiring a small absolute residual, but they are the conditions needed for a complete nonconvex convergence proof.

\begin{theorem}[Linear objective rate under a PL inequality]
\label{thm:PL-linear}
Assume the conditions of Proposition~\ref{prop:noncvx-decrease}, let $\eta_k\equiv\eta>0$, define
\[
f^\star:=\min_{x\in\mathcal K} f(x),
\]
and suppose that
\begin{equation}
 f(x)-f^\star
 \le \frac{1}{2\mu_{\rm PL}}\|\nabla f(x)\|^2,
 \qquad x\in\mathcal K,
\label{eq:PL-inequality}
\end{equation}
for some $\mu_{\rm PL}>0$. Then
\begin{equation}
 f(x_{k+1})-f^\star
 \le
 \frac{1}{1+(\mu_{\rm PL}\sigma_h/L_h^2)\eta}
 \bigl(f(x_k)-f^\star\bigr).
\label{eq:PL-rate}
\end{equation}
Thus the objective values converge $Q$-linearly for every fixed $\eta>0$ for which the exact subproblem is well defined.
\end{theorem}

\subsection{Application to the Three Vector Geometries}
\label{subsec:specializations}

\paragraph{Nonnegative orthant.}
For the exponential map,
\[
h(x)=\sum_i(x_i\log x_i-x_i),
\qquad
D_h(x,y)=\sum_i\left[x_i\log\frac{x_i}{y_i}-x_i+y_i\right].
\]
The backward-Euler equation in $u=\log x$ is exactly \eqref{eq:global-dual-step}. Other increasing maps, such as softplus, generate the potential defined in \eqref{eq:induced-potential} and are covered by the same argument.

\paragraph{Box constraints.}
The logistic map generates the two-sided entropy $h_{\Box}$ in \eqref{eq:box-induced-potential}. It does not generate the logarithmic barrier. The Gravidy-Box step is therefore the Bregman proximal step associated with this two-sided entropy.

\paragraph{Simplex.}
On the affine hull of $\Delta_n$, take
\[
h_{\Delta}(x)=\sum_i x_i\log x_i.
\]
Then $D_h$ is the KL divergence and \eqref{eq:global-bregman-step} is exactly the KL-proximal step in Section~\ref{sec:simplex-klprox}.

For these three entropy geometries, $h$ extends continuously to the boundary and $D_h(x^\star,x_0)$ is finite for every feasible $x^\star$ when $x_0$ is interior. Hence the convex objective-value result in Theorem~\ref{thm:global-convex} also covers boundary minimizers. By contrast, the nonconvex KL theorem assumes that the iterates stay in a compact subset of the interior. It does not prove convergence to an active boundary point. If a vector-domain sequence is known to converge to the boundary, the KKT conclusion follows from Section~\ref{sec:kkt-from-flow}.

\paragraph{Stiefel manifold.}
The results of this section do not apply to the implicit Cayley equation. An exact Cayley root is feasible, and the accepted algorithm decreases the objective, but these two facts alone are not sufficient for a global stationarity theorem. The local convergence of the inner Newton solve follows from the usual assumptions on the Fr\'echet derivative, as discussed in Appendices~\ref{app:stiefel-frechet}--\ref{app:dense-nr}.

\subsection{Relation with Linear Stability}
\label{subsec:linear-stability-convergence}

For the quadratic model
\[
f(x)=\tfrac12x^\top Qx-c^\top x,
\qquad Q\succeq0,
\]
and a fixed SPD metric $G$, backward Euler gives
\[
x_{k+1}=(G+\eta Q)^{-1}(Gx_k+\eta c).
\]
If $x^\star$ satisfies $Qx^\star=c$, then
\[
x_{k+1}-x^\star=(G+\eta Q)^{-1}G(x_k-x^\star).
\]
The eigenvalues of the iteration matrix lie in $(0,1]$ for every $\eta>0$, and they lie in $(0,1)$ when $Q$ is positive definite on the affine space under consideration. This is the optimization counterpart of backward Euler's unconditional stability on a dissipative linear system. The nonlinear Bregman results above should not be called A-stability statements; what they show is that the exact proximal step remains descending for every positive stepsize.

\subsection{Summary of the Convergence Guarantees}
\label{subsec:convergence-summary}

The main conclusions are the following.
\begin{itemize}
\item For convex objectives, the exact orthant, box, and simplex steps are descending for every $\eta_k>0$, and the last objective value satisfies the bound \eqref{eq:convex-last-iterate}.
\item Relative strong convexity gives the direct contraction \eqref{eq:relative-strong-contraction}, with no upper stepsize restriction.
\item An inexact inner solve preserves the convex $O(1/S_K)$ rate under a summable dual residual and bounded iterates. It does not automatically preserve the exact linear factor.
\item In the nonconvex case, the global-minimizer interpretation of the outer step, two-sided stepsize bounds, and a compact interior region are needed. Under the KL property, these assumptions give convergence to one critical point.
\item The Stiefel Cayley scheme lies outside this Bregman analysis and should be treated separately.
\end{itemize}

\section{Numerical Experiments}
\label{sec:experiments}

This section evaluates the four variants of \gravidy\ on representative constrained least-squares problems. Our main goal is to assess three points: the reduction of the objective and stationarity residual with respect to the number of outer iterations, the accuracy reached within a fixed iteration budget, and the preservation of feasibility. Since every implicit outer step requires an inner nonlinear solve, we report both iteration counts and wall-clock times. The main advantage we observe is the rapid convergence to high accuracy in terms of outer iterations.

\subsection{Experimental protocol}
\label{subsec:experimental-protocol}

All experiments are run for ten trials. Within each geometry, every method receives an independent copy of the same feasible initial point for a given trial. The problem data are fixed across the ten trials, while the initial points vary with the trial seed, except on the simplex where all methods start from the barycenter
\[
x_0=\frac{1}{n}\mathbf 1.
\]
This interior initialization is used throughout the simplex experiments.

For the three vector least-squares problems, the right-hand side is generated from a feasible reference point $x^\star$, so the reference optimal value is $\Phi^\star=0$.

For the vector problems, we measure first-order stationarity through the projected-gradient residual
\begin{equation}
\mathcal R_{\mathcal C}(x)
:=
\left\|x-\Pi_{\mathcal C}\bigl(x-\nabla\Phi(x)\bigr)\right\|_2,
\label{eq:numerical-kkt-residual}
\end{equation}
where \(\mathcal C\) is the orthant, simplex, or box. On the Stiefel manifold, we report the canonical Riemannian-gradient norm
\[
\mathcal R_{\mathrm{St}}(X):=\|\grad\Phi(X)\|_F
\]
and the feasibility error \(\|X^\top X-I\|_F\). The stopping tolerance is \(10^{-8}\). The direct Newton-type inner solves use tolerance \(10^{-10}\), the Armijo parameter is \(10^{-4}\), and the backtracking factor is \(1/2\). For projected Newton--Krylov on the Stiefel manifold, the nonlinear Cayley residual is tested against \(10^{-12}(1+\|X_k\|_F)\). GMRES uses the adaptive relative tolerance
\[
\min\!\left\{10^{-1},\max\!\left(10^{-6},0.1\|F_k(Y)\|_F\right)\right\},
\]
with zero absolute tolerance, at most 200 GMRES iterations, and at most five inner Newton steps. The curves show averages over the ten trials. Table~\ref{tab:definitive-summary} reports medians; the stationarity column contains the smallest residual reached within the allotted iteration budget.

The experiments were run with Python~3.13.5, NumPy~2.3.2, SciPy~1.16.1, Matplotlib~3.10.5, and OpenBLAS~0.3.33 on macOS ARM64. The package versions and all experimental parameters are pinned in the accompanying configuration files.

\subsection{Test problems and competing methods}
\label{subsec:test-problems}

\paragraph{Nonnegative orthant.}
We consider the nonnegative least-squares problem
\[
\min_{x\geq0}\ \frac12\|Ax-b\|_2^2
\]
with \(m=n=120\). The ground-truth vector is sparse, with density \(0.15\), and the data are generated so that the optimal value is known. We use \(\eta=300\) and at most 400 outer iterations. We compare \gravidypos\ with Newton and MGN inner solves against accelerated projected gradient, projected Barzilai--Borwein, and multiplicative updates~\cite{nesterov2018lectures,BirginMartinezRaydan2000,lee1999learning}.

\paragraph{Simplex.}
We solve
\[
\min_{x\in\Delta_n}\ \frac12\|Ax-b\|_2^2
\]
with \(m=n=40\) and a matrix $A$ with prescribed condition number \(10^3\). We use \(\eta=100\) and at most 400 outer iterations. The two \gravidydelta\ variants are the KL-proximal Newton--KKT method and the reduced-logit MGN method. They are compared with projected gradient, accelerated projected gradient, and entropic mirror descent~\cite{nesterov2018lectures,BeckTeboulle2003,KivinenWarmuth1997,NemirovskyYudin1983}.

\paragraph{Box constraints.}
We consider a box-constrained least-squares problem
\[
\min_{\ell\leq x\leq b}\ \frac12\|Ax-c\|_2^2
\]
with \(m=n=120\). We use \(\eta=300\) and at most 400 outer iterations. The Newton and MGN variants of \gravidybox\ are compared with accelerated projected gradient and projected Barzilai--Borwein~\cite{nesterov2018lectures,BirginMartinezRaydan2000}.

\paragraph{Stiefel manifold.}
We minimize
\[
\Phi(X)
=
\frac12\sum_{j=1}^{p}
\left\langle Q^{(j)}x^{(j)},x^{(j)}\right\rangle,
\qquad
X\in\St(n,p),
\]
with \(n=200\), \(p=2\), and symmetric positive definite matrices \(Q^{(j)}\) of condition number \(10^3\). The two implicit Cayley variants use either projected Newton--Krylov or dense Newton for the inner equation. We compare them with the feasible Cayley method of Wen and Yin~\cite{WenYin2013} and a Riemannian gradient method with QR retraction. The \gravidyst\ methods are limited to 500 outer iterations, while the two first-order competitors are allowed up to 50,000 iterations.

\subsection{Results}
\label{subsec:numerical-results}

\begin{figure}[ht!]
\centering
\includegraphics[width=.485\linewidth]{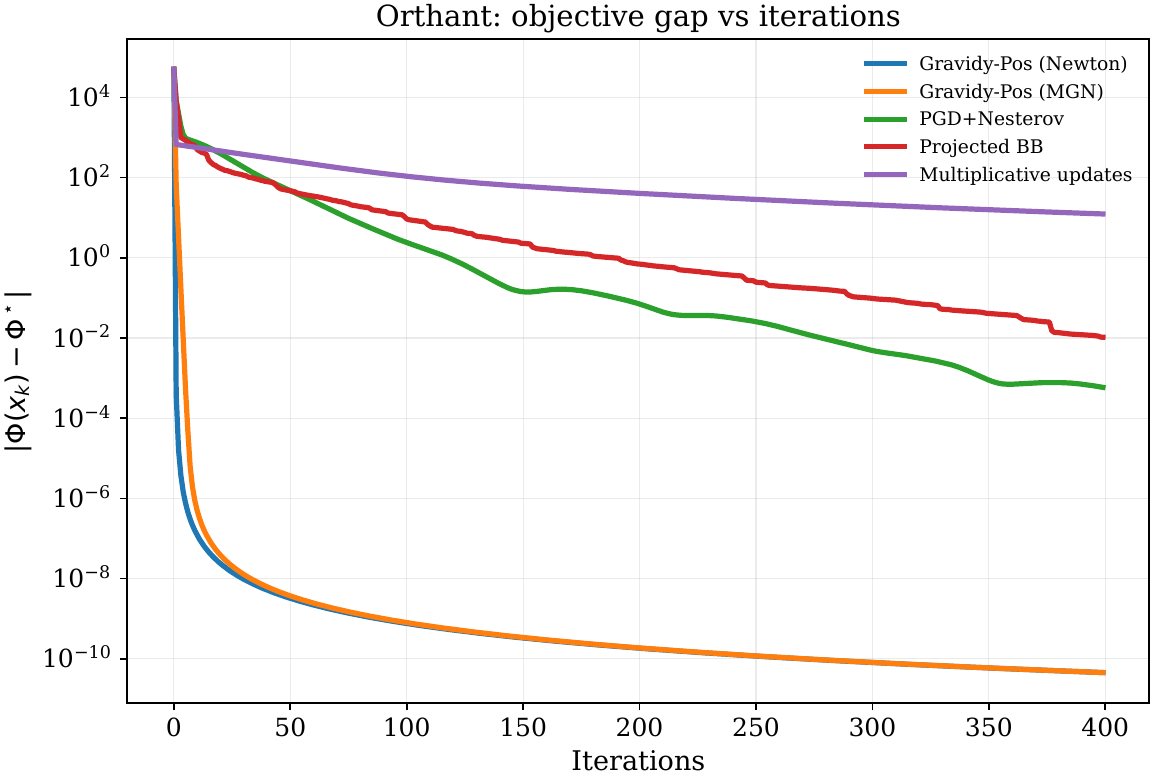}\hfill
\includegraphics[width=.485\linewidth]{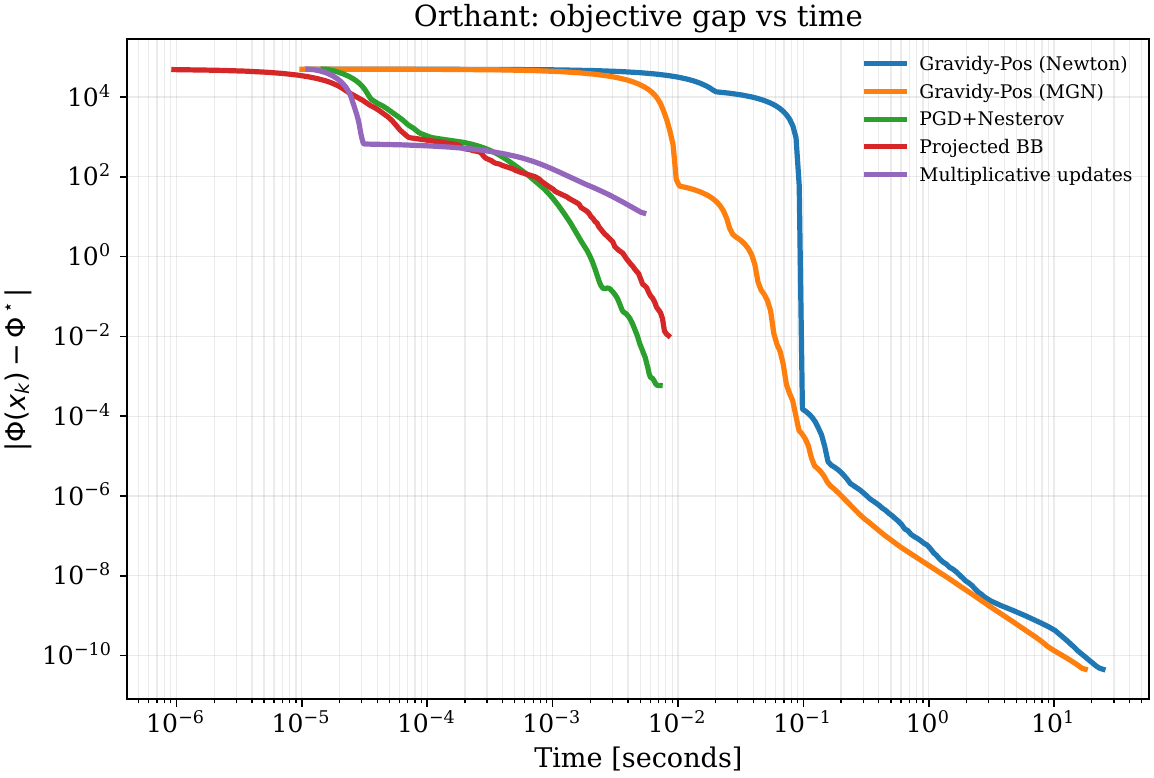}
\caption{Nonnegative least squares, \(m=n=120\). Left: objective gap versus outer iterations. Right: objective gap versus wall-clock time. The curves are averages over ten trials.}
\label{fig:pos-definitive}
\end{figure}

\paragraph{Orthant.}
Figure~\ref{fig:pos-definitive} shows that both \gravidypos\ variants reduce the objective gap by more than fourteen orders of magnitude. Most of this decrease occurs during the first few outer iterations, after which the two curves continue steadily toward high accuracy. The Newton and MGN variants have almost identical outer behavior and reach median stationarity residuals of approximately \(2\times10^{-6}\). Within the same 400-iteration budget, the best projected baseline remains above \(4.9\times10^{-2}\). The implicit methods therefore require more work per outer iteration, but each iteration produces a much stronger reduction of the objective and stationarity residual.

\begin{figure}[ht!]
\centering
\includegraphics[width=.485\linewidth]{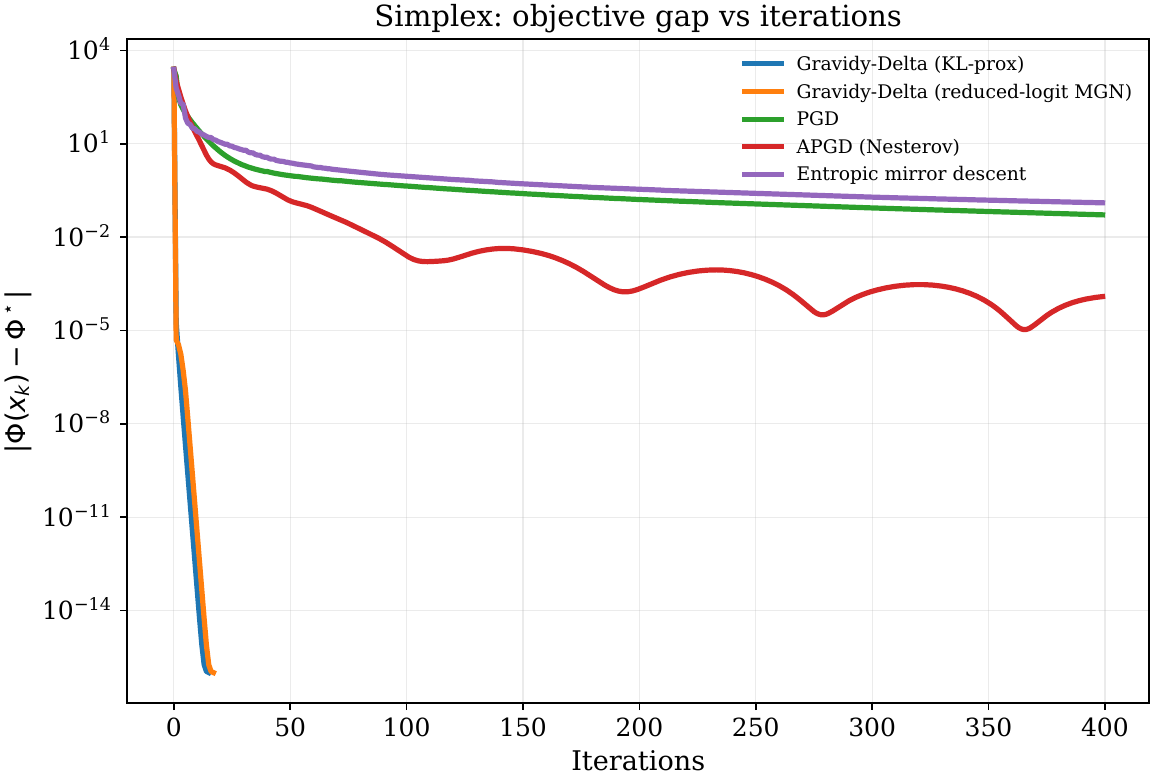}\hfill
\includegraphics[width=.485\linewidth]{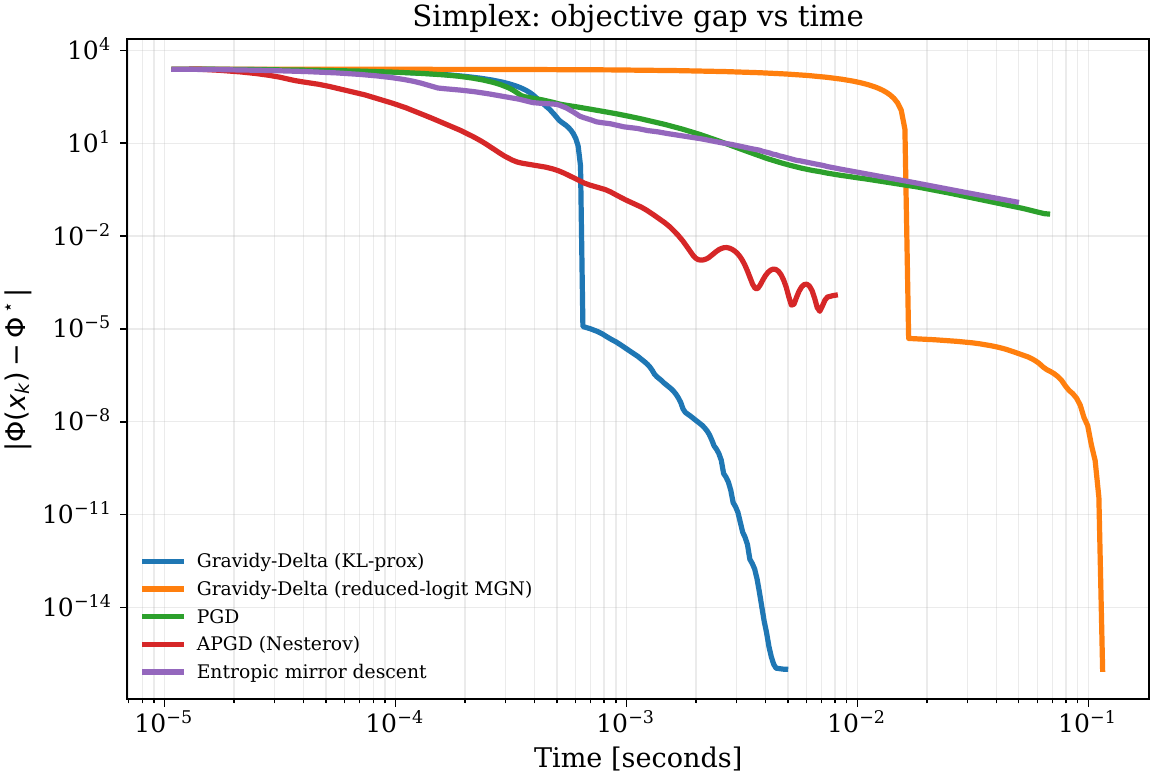}
\caption{Simplex-constrained least squares, \(m=n=40\), with condition number \(10^3\). Left: objective gap versus outer iterations. Right: objective gap versus wall-clock time. All methods start from \(x_0=\frac1n\mathbf 1\).}
\label{fig:simplex-definitive}
\end{figure}

\paragraph{Simplex.}
The simplex experiment gives the clearest illustration of the benefit of the implicit geometric step. As shown in Figure~\ref{fig:simplex-definitive}, both \gravidydelta\ variants reach an objective gap close to machine precision in fewer than twenty outer iterations. The median recorded indices at termination are 15 for the KL-proximal method and 17 for reduced-logit MGN. Both methods reach the stationarity tolerance in all ten trials, with median final residuals below \(9\times10^{-9}\). The KL-proximal Newton--KKT implementation is also very fast in wall-clock time. By comparison, the explicit projected and mirror methods remain several orders of magnitude less accurate after 400 iterations.

\begin{figure}[ht!]
\centering
\includegraphics[width=.485\linewidth]{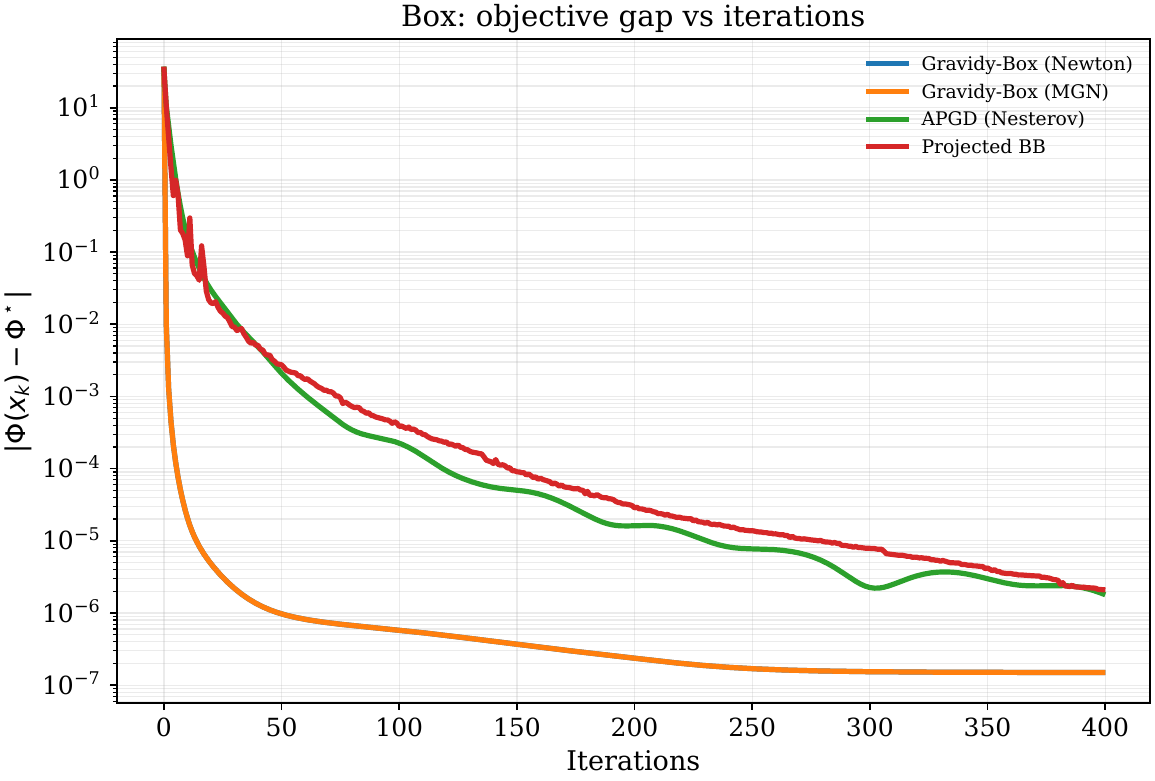}\hfill
\includegraphics[width=.485\linewidth]{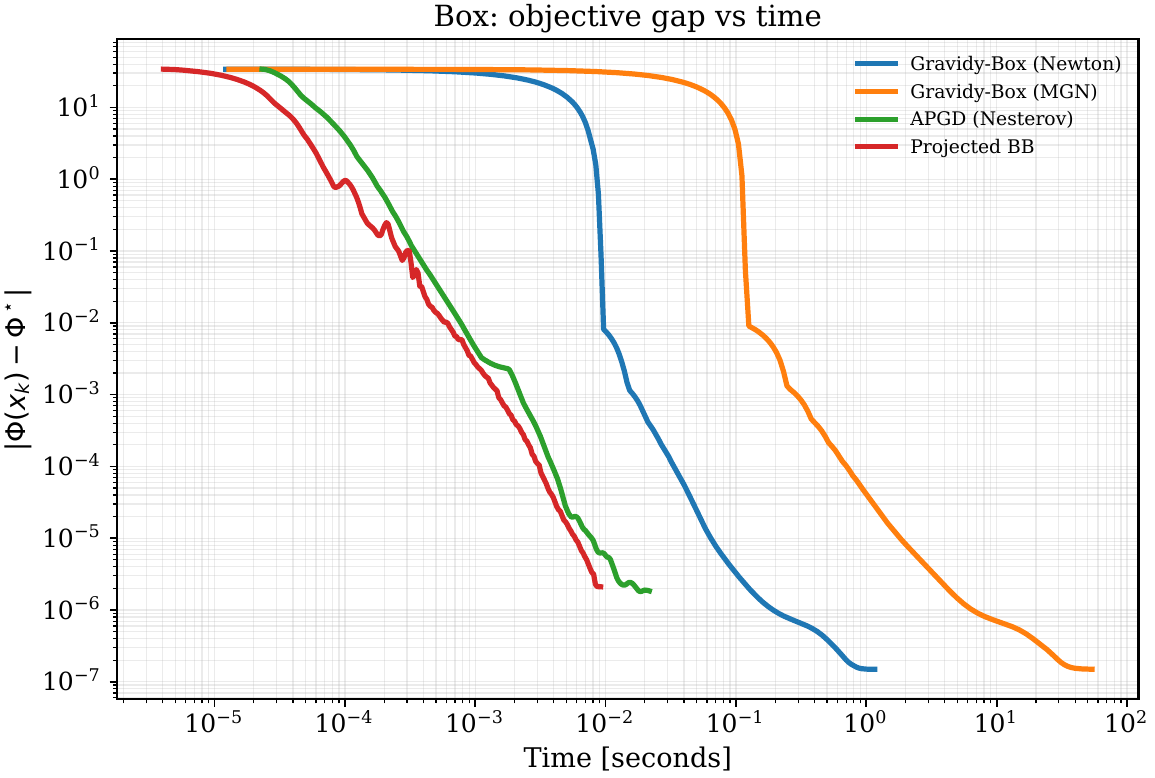}
\caption{Box-constrained least squares, \(m=n=120\). Left: objective gap versus outer iterations. Right: objective gap versus wall-clock time. The two \gravidybox\ variants follow nearly identical outer trajectories.}
\label{fig:box-definitive}
\end{figure}

\paragraph{Box constraints.}
Figure~\ref{fig:box-definitive} shows a similarly strong reduction for \gravidybox. The two implicit variants produce essentially indistinguishable outer trajectories and reach objective gaps of order \(10^{-7}\). Their best median stationarity residual is \(4.57\times10^{-6}\), compared with \(4.39\times10^{-5}\) for projected Barzilai--Borwein and \(1.25\times10^{-4}\) for accelerated projected gradient. Hence the implicit geometry again provides about one additional order of magnitude in first-order accuracy within the same outer-iteration budget. At this dimension, the direct Newton inner solve is considerably more economical than MGN while producing the same outer sequence.

\begin{figure}[ht!]
\centering
\includegraphics[width=.485\linewidth]{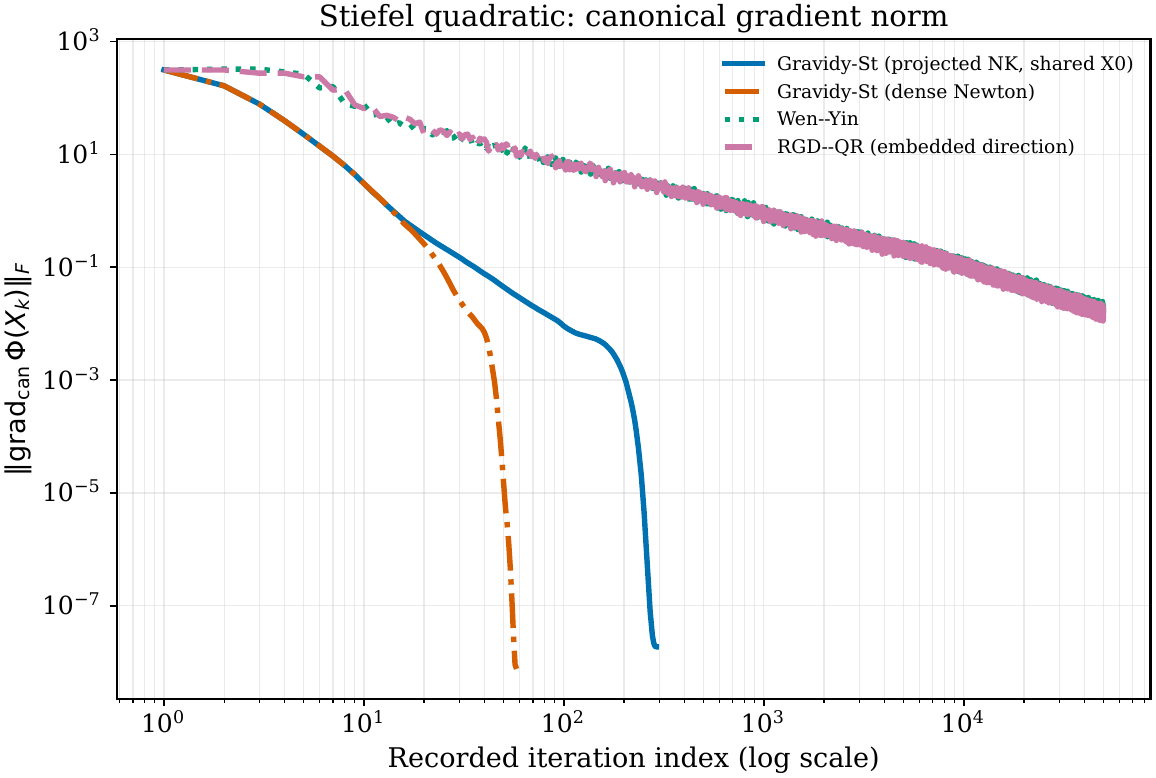}\hfill
\includegraphics[width=.485\linewidth]{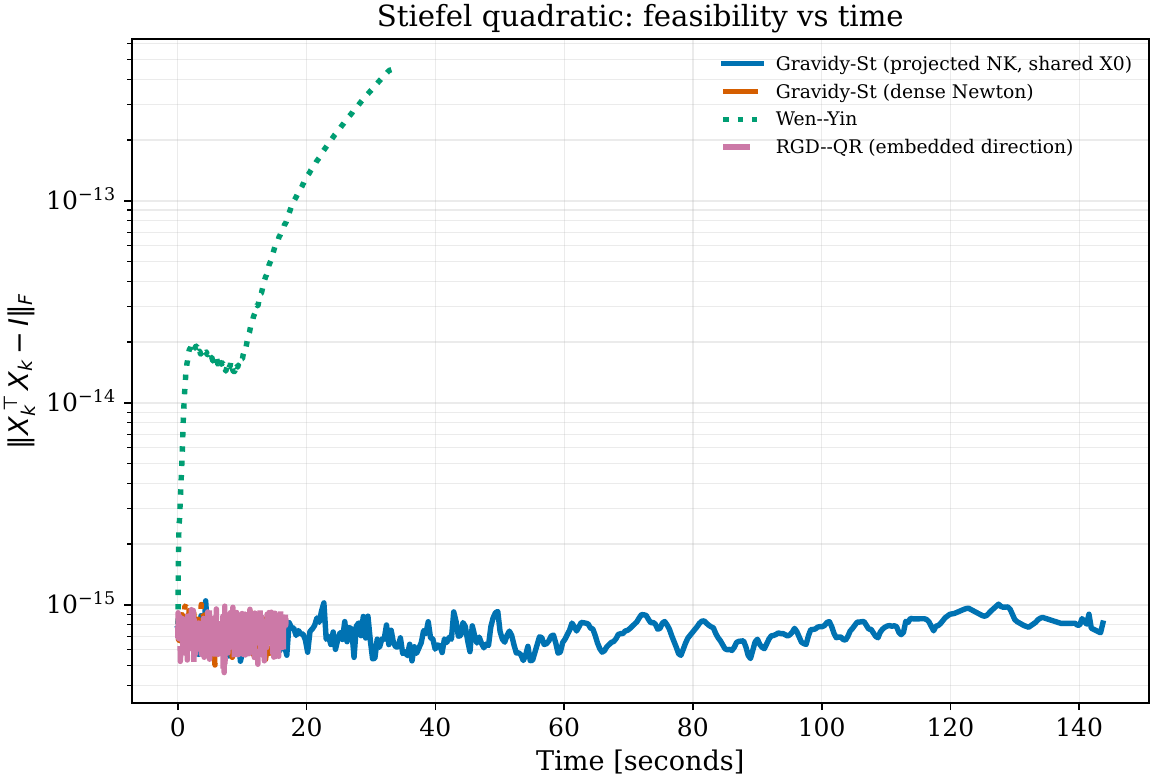}
\caption{Stiefel quadratic problem on \(\St(200,2)\). Left: canonical Riemannian-gradient norm versus recorded iteration index. Right: orthogonality error versus wall-clock time. The implicit Cayley methods reach high stationarity accuracy while preserving feasibility at machine precision.}
\label{fig:stiefel-definitive}
\end{figure}

\paragraph{Stiefel manifold.}
The Stiefel results are reported in Figure~\ref{fig:stiefel-definitive}. Dense Newton reduces the canonical gradient norm from order \(10^2\) to below \(10^{-8}\) in a median of about 50 outer iterations. Projected Newton--Krylov reaches \(1.71\times10^{-8}\) in a median of about 245 iterations. The two first-order competitors remain near \(10^{-2}\) after 50,000 recorded iterations. The feasibility plot confirms that all geometric methods preserve orthogonality to numerical precision. In particular, the median feasibility errors of the two \gravidyst\ variants are below \(8\times10^{-16}\).

\begin{table}[ht!]
\centering
\caption{Definitive numerical results. The values are medians over ten trials. ``Best stat.'' denotes the smallest KKT residual reached within the iteration budget; on the Stiefel manifold it denotes the smallest canonical Riemannian-gradient norm.}
\label{tab:definitive-summary}
\scriptsize
\renewcommand{\arraystretch}{1.08}
\begin{adjustbox}{max width=\linewidth}
\begin{tabular}{llrrrr}
\toprule
Geometry & Method & Best stat. & Feasibility & Recorded index & Wall time [s] \\
\midrule
Orthant
& \gravidypos\ (Newton) & $2.013\times10^{-6}$ & $0$ & 399.0 & 23.279 \\
& \gravidypos\ (MGN) & $\mathbf{1.975\times10^{-6}}$ & $0$ & 399.0 & 17.003 \\
& PGD+Nesterov & $1.526\times10^{-1}$ & $0$ & 399.0 & 0.007 \\
& Projected BB & $4.923\times10^{-2}$ & $0$ & 399.0 & 0.008 \\
& Multiplicative updates & $7.820$ & $0$ & 399.0 & 0.005 \\
\midrule
Simplex
& \gravidydelta\ (KL-prox) & $8.722\times10^{-9}$ & $0$ & 15.0 & 0.004 \\
& \gravidydelta\ (reduced-logit MGN) & $\mathbf{8.378\times10^{-9}}$ & $0$ & 17.0 & 0.114 \\
& PGD & $5.951\times10^{-1}$ & $0$ & 399.0 & 0.065 \\
& APGD (Nesterov) & $1.229\times10^{-1}$ & $0$ & 399.0 & 0.007 \\
& Entropic mirror descent & $6.336\times10^{-1}$ & $0$ & 399.0 & 0.048 \\
\midrule
Box
& \gravidybox\ (Newton) & $4.571\times10^{-6}$ & $0$ & 399.0 & 1.112 \\
& \gravidybox\ (MGN) & $\mathbf{4.571\times10^{-6}}$ & $0$ & 399.0 & 53.744 \\
& APGD (Nesterov) & $1.250\times10^{-4}$ & $0$ & 399.0 & 0.007 \\
& Projected BB & $4.393\times10^{-5}$ & $0$ & 399.0 & 0.009 \\
\midrule
Stiefel
& \gravidyst\ (projected NK) & $1.709\times10^{-8}$ & $7.729\times10^{-16}$ & 244.5 & 115.753 \\
& \gravidyst\ (dense Newton) & $\mathbf{8.247\times10^{-9}}$ & $4.973\times10^{-16}$ & 49.5 & 12.268 \\
& Wen--Yin & $7.725\times10^{-3}$ & $4.632\times10^{-13}$ & 49999.0 & 32.913 \\
& RGD--QR & $7.848\times10^{-3}$ & $8.013\times10^{-16}$ & 49999.0 & 16.404 \\
\bottomrule
\end{tabular}
\end{adjustbox}
\end{table}

\subsection{Main observations}
\label{subsec:main-observations}

The experiments lead to three main observations. First, the implicit methods are particularly effective in terms of outer iterations. On the simplex and Stiefel manifold, high stationarity accuracy is reached in only a few tens or a few hundreds of iterations, while the first-order competitors remain far from the same accuracy after much larger iteration budgets. The orthant and box experiments show the same behavior in the objective gap, with a very sharp initial decrease followed by steady convergence toward high accuracy.

Second, the geometry is respected numerically. The vector-domain methods remain feasible by construction, and the Cayley updates preserve the Stiefel constraint to machine precision. This is consistent with the continuous and discrete constructions developed in Sections~\ref{sec:gravidy} and~\ref{sec:algorithms}.

Third, the choice of inner solver matters mainly for computational cost, not for the outer trajectory. Newton and MGN produce almost identical curves on the orthant and box, while Newton--KKT and reduced-logit MGN have very similar iteration counts on the simplex. For the dimensions considered here, the direct Newton variants are generally the most economical. The matrix-free Newton--Krylov formulation remains useful for larger Stiefel problems where assembling the dense Jacobian is no longer practical.

\subsection{A Sparse Elastic Obstacle Problem}
\label{subsec:elastic-obstacle}

We finally consider a structured sparse problem arising from the deformation
of an elastic membrane above a prescribed obstacle. This experiment has two
purposes. First, it illustrates the geometry of the computed solution through
a simple physical example. Second, it shows that the orthant construction can
be implemented without forming a dense least-squares matrix or a dense
Hessian.

\paragraph{Continuous problem.}
Let
\[
\Omega=[0,3\pi]\times[0,3\pi],
\]
and consider the obstacle
\[
\phi(x,y)
=
\max\{0,\sin(x)\}\max\{0,\sin(y)\}.
\]
The elastic obstacle problem consists in finding a membrane
\(w:\Omega\to\mathbb R\) of minimum elastic energy that remains above
\(\phi\). With homogeneous Dirichlet boundary conditions, the problem can be
written as
\begin{equation}
\min_{w\in H_0^1(\Omega)}
\left\{
\frac12\int_{\Omega}\|\nabla w(\xi)\|_2^2\,d\xi
\;:\;
w\geq\phi\ \text{a.e. in }\Omega
\right\}.
\label{eq:eop-continuous}
\end{equation}
The obstacle vanishes on the boundary of the domain, so these boundary
conditions are compatible with the inequality \(w\geq\phi\).

We discretize the interior of \(\Omega\) on an \(N\times N\) uniform grid,
with
\[
h=\frac{3\pi}{N+1}.
\]
Let \(Q\in\mathbb R^{N^2\times N^2}\) be the standard sparse
five-point discretization of the negative Laplacian. The discrete problem is
\begin{equation}
\min_{w\geq\phi}\frac12 w^\top Qw.
\label{eq:eop-discrete}
\end{equation}
The matrix \(Q\) is symmetric positive definite and has only about five
nonzero entries per row.

To express the constraint in the standard nonnegative form, we introduce the
gap between the membrane and the obstacle,
\[
v:=w-\phi.
\]
Then \(w\geq\phi\) is equivalent to \(v\geq0\), and
\[
\frac12(\phi+v)^\top Q(\phi+v)
=
\frac12v^\top Qv+v^\top Q\phi+\frac12\phi^\top Q\phi.
\]
After removing the constant term, we obtain the nonnegative quadratic problem
\begin{equation}
\min_{v\geq0}
\Phi(v)
:=
\frac12v^\top Qv-p^\top v,
\qquad
p:=-Q\phi.
\label{eq:eop-shifted}
\end{equation}
Its gradient and Hessian are
\[
\nabla\Phi(v)=Qv-p,
\qquad
\nabla^2\Phi(v)=Q.
\]
The objective values reported below correspond to the shifted objective in
\eqref{eq:eop-shifted}; they may therefore be negative even though the
original elastic energy is nonnegative.

\paragraph{Dedicated sparse \gravidy-Pos solver.}
We use the exponential representation
\[
v=e^z,
\]
where the exponential is applied componentwise. The corresponding parameter
dynamics are
\[
\dot z=-(Qe^z-p).
\]
At outer iteration \(k\), the implicit update is defined by
\begin{equation}
F_k(z)
:=
z-z_k+\eta_k(Qe^z-p)
=
0.
\label{eq:eop-implicit-residual}
\end{equation}
This equation is solved by the transformed Newton method introduced in
Section~\ref{sec:algorithms}. At an inner iterate \(z\), let
\[
v=e^z,
\qquad
D=\operatorname{Diag}(v).
\]
Instead of solving the generally nonsymmetric Newton system directly, we solve
the equivalent symmetric positive definite system
\begin{equation}
\left(
I+\eta_kD^{1/2}QD^{1/2}
\right)y
=
-D^{1/2}F_k(z),
\qquad
\delta z=D^{-1/2}y.
\label{eq:eop-transformed-newton}
\end{equation}
All matrix-vector products use the sparse matrix \(Q\). For the unmodified
exponential map, the linear systems are solved by conjugate gradients with the
Jacobi preconditioner
\[
M_{ii}=1+\eta_kv_iQ_{ii}.
\]
The PCG tolerance is relaxed when the nonlinear residual is large and
tightened as the root is approached.

For numerical safety, the implementation uses the safeguarded map
\[
\widehat v_i(z_i)
=
\exp\bigl(\max\{z_i,z_{\min}\}\bigr),
\qquad
z_{\min}=\log(10^{-16}),
\]
and uses the following selected generalized derivative,
\[
d_i(z_i)
=
\begin{cases}
\widehat v_i(z_i), & z_i\geq z_{\min},\\
0, & z_i<z_{\min}.
\end{cases}
\]
The implemented residual is
\[
\widehat F_k(z)
=
z-z_k+\eta_k\bigl(Q\widehat v(z)-p\bigr).
\]
Let
\[
\mathcal I=\{i:d_i>0\},
\qquad
\mathcal C_0=\{i:d_i=0\},
\]
and let \(D_{\mathcal I}\) contain the positive entries \(d_i\). The
transformed generalized Newton equation is solved on \(\mathcal I\):
\[
\left(
I+\eta_kD_{\mathcal I}^{1/2}
Q_{\mathcal I\mathcal I}
D_{\mathcal I}^{1/2}
\right)y_{\mathcal I}
=
-D_{\mathcal I}^{1/2}\widehat F_{k,\mathcal I}(z),
\qquad
\delta z_{\mathcal I}
=
D_{\mathcal I}^{-1/2}y_{\mathcal I}.
\]
The zero-derivative rows are recovered from the original Newton equation:
\[
\delta z_{\mathcal C_0}
=
-\widehat F_{k,\mathcal C_0}(z)
-\eta_k
Q_{\mathcal C_0\mathcal I}
D_{\mathcal I}\delta z_{\mathcal I}.
\]
For the safeguarded system, the Jacobi diagonal is
\[
M_{ii}=1+\eta_kd_iQ_{ii}.
\]

We use a target outer stepsize \(\eta_{\max}=300\). If an implicit solve fails,
the trial stepsize is divided by two. After a successful solve, it is increased
by a factor \(1.5\), up to \(\eta_{\max}\). Each accepted implicit step
satisfies
\[
\|\widehat F_k(z_{k+1})\|_2\leq10^{-8},
\]
and the inner Newton method uses at most \(20\) iterations. An Armijo line
search is applied to the merit function
\[
z\mapsto\frac12\|\widehat F_k(z)\|_2^2.
\]
Trial points exceeding the upper exponential safety threshold are rejected.

\paragraph{Competitors and protocol.}
We compare the sparse \gravidy-Pos Newton--PCG implementation with two
projected first-order methods.

The first method is projected gradient descent,
\begin{equation}
v_{k+1}
=
\Pi_{\mathbb R_+^{N^2}}
\left(
v_k-\frac1L(Qv_k-p)
\right),
\label{eq:eop-pgd}
\end{equation}
where the exact largest eigenvalue of the Kronecker Laplacian is
\[
L
=
\frac{8}{h^2}
\cos^2\left(\frac{\pi}{2(N+1)}\right).
\]
The second competitor is a projected Barzilai--Borwein method equipped with
monotone backtracking.

We use \(N=100\), corresponding to \(10\,000\) variables. The projected
methods start from \(v_0=0\), while \gravidy-Pos starts from
\[
v_0=10^{-8}\mathbf1,
\]
because its parameter representation requires a strictly positive point.
Projected gradient and projected BB are allowed, respectively, \(10\,000\)
and \(4\,000\) iterations. The stopping measure is the projected KKT
residual
\begin{equation}
R_+(v)
:=
\left\|
v-\Pi_{\mathbb R_+^{N^2}}
\bigl(v-\nabla\Phi(v)\bigr)
\right\|_2.
\label{eq:eop-kkt-residual}
\end{equation}
Solver times are medians over three runs after one warm-up run. Problem
construction, figure generation, and file writing are excluded from the
reported times. The construction of the sparse problem takes approximately
\(3.6\times10^{-3}\) seconds.

\begin{table}[ht!]
\centering
\caption{Results for the sparse elastic obstacle problem with
\(N=100\), corresponding to \(10\,000\) variables. The objective is the
shifted quadratic objective in \eqref{eq:eop-shifted}. For \gravidy-Pos, the
number of steps refers to implicit outer steps; for the projected methods, it
refers to first-order iterations.}
\label{tab:eop2d-results}
\begin{tabular}{lrrrr}
\toprule
Method
&
Final \(R_+(v)\)
&
Final \(\Phi(v)\)
&
Steps
&
Time [s]
\\
\midrule
\gravidy-Pos Newton--PCG
&
\(2.18\times10^{-9}\)
&
\(-2.90831\times10^{2}\)
&
\(78\)
&
\(2.115\)
\\
Projected gradient
&
\(3.03\times10^{-1}\)
&
\(-2.90676\times10^{2}\)
&
\(10\,000\)
&
\(0.664\)
\\
Projected BB
&
\(9.73\times10^{-6}\)
&
\(-2.90831\times10^{2}\)
&
\(4\,000\)
&
\(3.474\)
\\
\bottomrule
\end{tabular}
\end{table}

\begin{figure}[ht!]
\centering
\begin{minipage}[t]{0.48\textwidth}
\centering
\includegraphics[width=\linewidth]{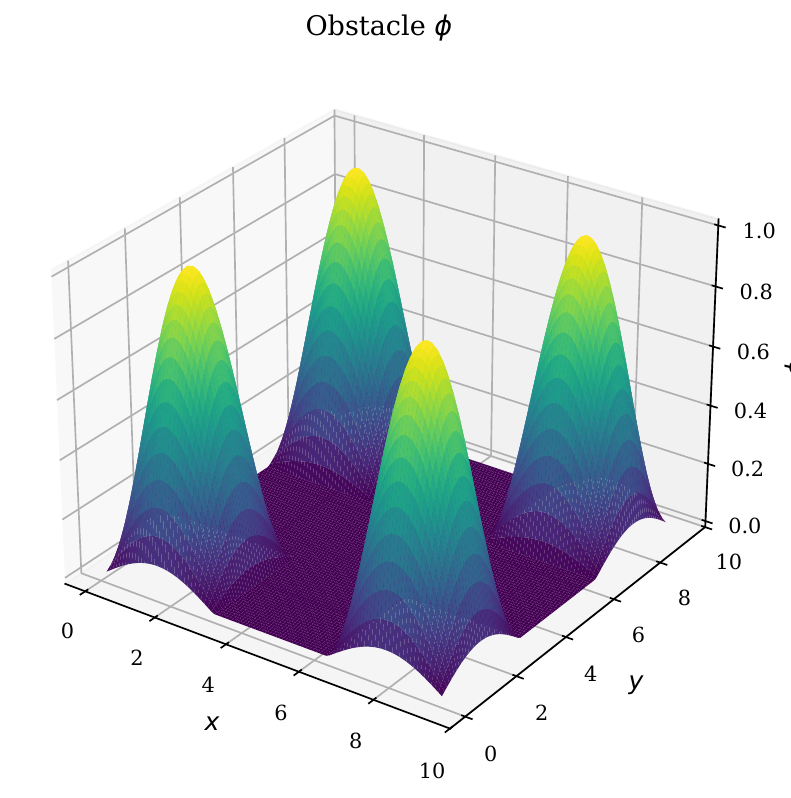}
\par\smallskip
{\small (a) Obstacle \(\phi\).}
\end{minipage}
\hfill
\begin{minipage}[t]{0.48\textwidth}
\centering
\includegraphics[width=\linewidth]{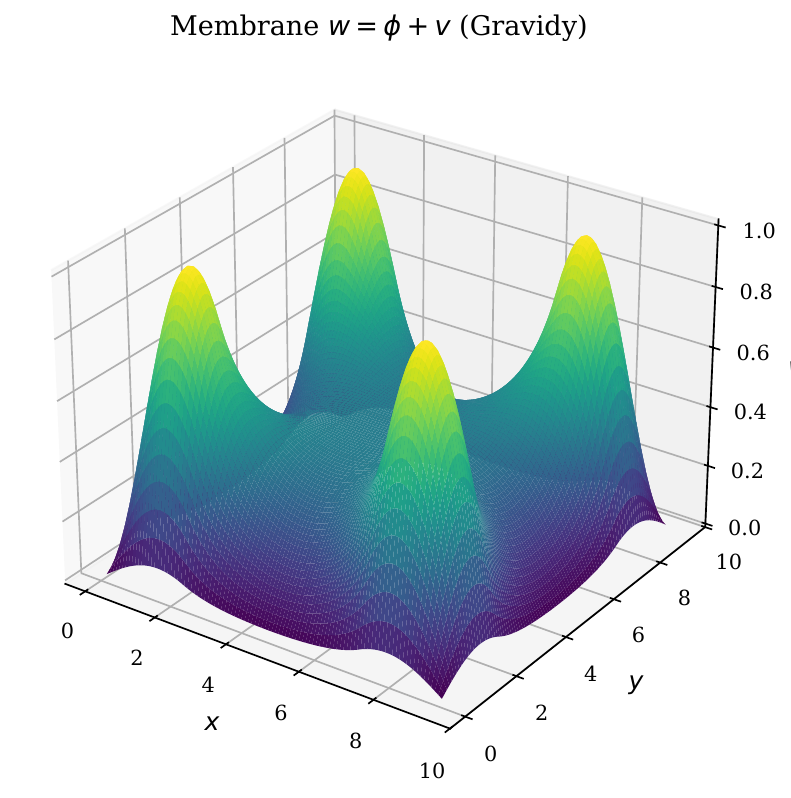}
\par\smallskip
{\small (b) Computed membrane \(w=\phi+v\).}
\end{minipage}

\medskip

\begin{minipage}[t]{0.48\textwidth}
\centering
\includegraphics[width=\linewidth]{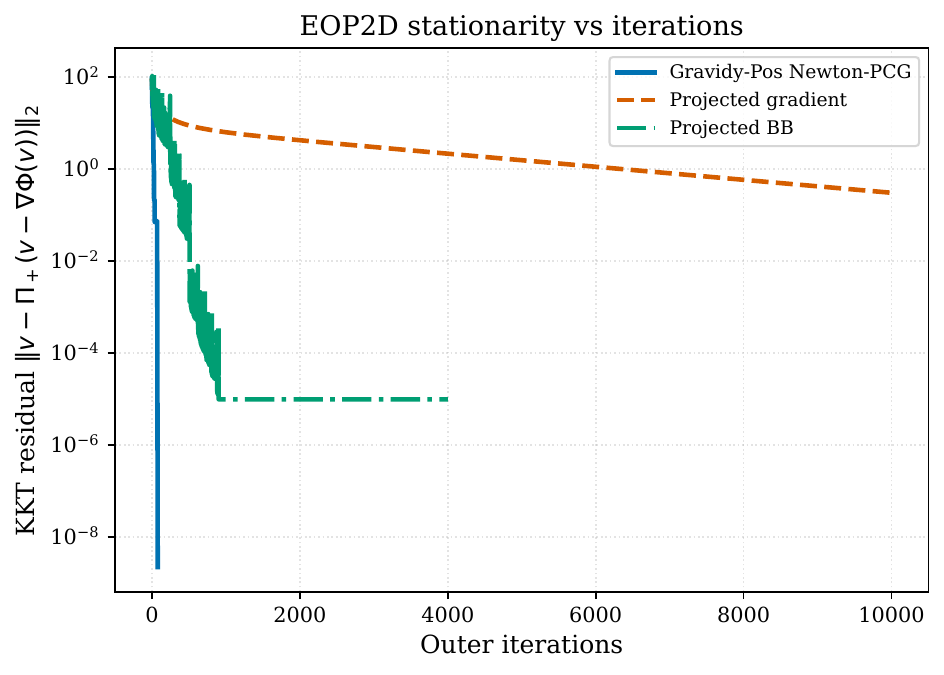}
\par\smallskip
{\small (c) Projected KKT residual versus iterations.}
\end{minipage}
\hfill
\begin{minipage}[t]{0.48\textwidth}
\centering
\includegraphics[width=\linewidth]{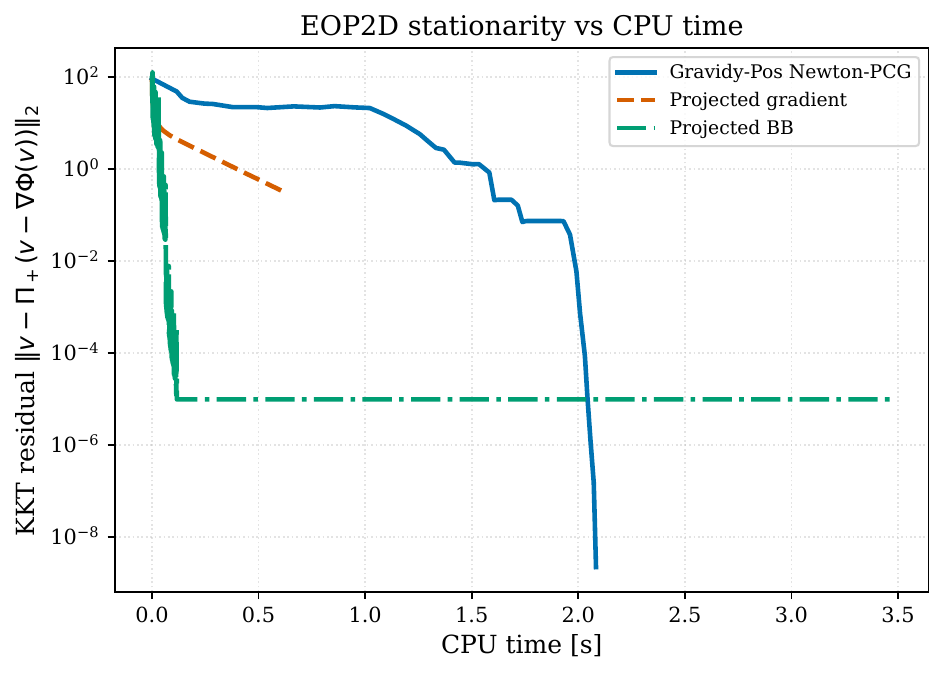}
\par\smallskip
{\small (d) Projected KKT residual versus CPU time.}
\end{minipage}

\caption{Sparse elastic obstacle problem with \(10\,000\) variables.
The optimized membrane remains above the obstacle and smoothly connects its
different peaks. The projected methods make inexpensive first-order updates,
whereas each \gravidy-Pos outer step requires the solution of a sparse
nonlinear system. Nevertheless, \gravidy-Pos reaches a projected KKT residual
below \(10^{-8}\) in \(78\) outer steps and \(2.115\) seconds.}
\label{fig:eop2d-results}
\end{figure}

\paragraph{Results.}
Figure~\ref{fig:eop2d-results} shows the obstacle, the membrane computed by
\gravidy-Pos, and the convergence histories. The membrane remains above the
obstacle by construction. It touches the obstacle on the active region and
smoothly bridges the lower regions between the four peaks.

Projected gradient has the lowest total runtime, but remains far from
stationarity after \(10\,000\) iterations, with
\[
R_+(v)\approx3.0\times10^{-1}.
\]
Projected BB is substantially more effective and reaches the same objective
value as \gravidy-Pos to the displayed precision. It reduces the residual
quickly at moderate accuracy, but reaches its iteration limit with
\[
R_+(v)\approx9.7\times10^{-6}.
\]

The sparse \gravidy-Pos method requires more work per outer step, but continues
to high accuracy. It reaches
\[
R_+(v)\approx2.2\times10^{-9}
\]
in \(78\) outer steps. Its median runtime is \(2.115\) seconds, compared with
\(3.474\) seconds for projected BB. All accepted implicit equations satisfy
the prescribed tolerance: their median residual is
\(3.71\times10^{-10}\), and their maximum residual is
\(9.42\times10^{-9}\). The complete run uses \(451\) nonlinear Newton
iterations and \(31\,884\) sparse PCG iterations.

This example therefore complements the dense synthetic tests. It shows that
the implicit orthant method can exploit a sparse structured Hessian directly,
without forming a dense matrix or a square-root factorization, and can reach
high stationarity accuracy on a problem with \(10\,000\) variables.

\subsection{Reproducibility}
\label{subsec:reproducibility}

The definitive experiments are generated by a single manuscript-aligned configuration and the script \texttt{run\_definitive\_paper.sh}. The result folder contains the complete trajectories in compressed numeric \texttt{.npz} files, the fixed environment and experiment configurations, execution logs, vector PDF figures, and CSV, JSON, and \LaTeX\ summaries. The implementation and reproduction scripts are available in the \gravidy\ repository at \url{https://github.com/vleplat/GRAVIDY}. The repository also includes separate visual benchmark scripts for the four geometries.

\section{Conclusion}\label{sec:conclusion}

We introduced \gravidy, a framework for constrained optimization in which the constraints are incorporated directly into the geometry of the dynamics. The resulting motion remains feasible by construction. The precise geometric mechanism depends on the constraint set. For the nonnegative orthant and box constraints, feasibility is encoded through smooth reparameterizations and the corresponding Hessian--Riemannian geometries. On the simplex, the dynamics follows the Fisher--Shahshahani geometry and its KL-proximal discretization. On the Stiefel manifold, the method is defined directly from the intrinsic Riemannian structure.

For the vector domains, the implicit steps admit an exact Bregman-proximal interpretation in the primal variable. This provides a direct link between the continuous flow, the implicit algorithm, and the convergence analysis. In the convex case, the objective decreases monotonically for every positive stepsize, and the last iterate satisfies the usual sublinear estimate. Relative strong convexity yields a linear contraction without an upper restriction on the stepsize. In the nonconvex case, we obtained convergence under the Kurdyka--Łojasiewicz property when the iterates remain in a compact subset of the interior and the accepted steps satisfy the required decrease and relative-error conditions.

We also studied how first-order optimality appears from the dynamics. For the orthant, simplex, and box constraints, convergent trajectories generated from the interior satisfy the corresponding KKT conditions. The same conclusion holds for convergent implicit iterates generated from the interior when their stepsizes are bounded away from zero. The boundary inequalities are recovered from the behavior of the parameter or dual variables, while complementarity follows from the loss of mobility at active constraints. On the Stiefel manifold, the vanishing of the canonical Riemannian gradient is equivalent to the usual equality-constrained first-order condition.

The numerical experiments support the practical interest of the framework.
On the orthant, simplex, box, and Stiefel test problems, the implicit methods
remain stable with large steps and reach high accuracy in a small number of
outer iterations. In particular, the simplex variants reach near machine
precision in only a few implicit steps, while the orthant and box methods
provide substantial reductions of the KKT residual. On the Stiefel manifold,
the implicit Cayley method preserves orthogonality up to machine precision and
reaches a small canonical-gradient norm in a limited number of outer
iterations. The elastic-obstacle experiment further shows that the orthant
construction extends naturally to a large sparse structured problem and can
exploit the discrete Laplacian directly. These results also show that the
choice of the inner solver matters: direct Newton methods are particularly
effective for the moderate dense problems considered here, while transformed
Newton--PCG provides a natural implementation for larger sparse systems.
Modified Gauss--Newton and Newton--Krylov variants remain useful when
additional regularization or matrix-free computations are preferred.

Several directions remain open. On the vector domains, it would be useful to design adaptive rules for the implicit stepsize, based for example on the nonlinear residual, the predicted decrease, or the progress of the inner solver. Inexact inner solves could also be controlled more closely through relative stopping criteria that preserve the convergence guarantees. Other natural directions include inertial variants, more general mirror geometries, and applications where the constraint-induced metric carries additional statistical or physical information.

The Stiefel construction suggests a broader geometric extension. The skew-symmetric matrix
\[
A(X)
=
\nabla\Phi(X)X^\top
-
X\nabla\Phi(X)^\top
\]
belongs to the Lie algebra \(\mathfrak{so}(n)\) of the orthogonal group. It can be viewed as an infinitesimal generator of a feasible motion: it describes, at first order, how the iterate should rotate while remaining on the constraint manifold. The Cayley transform then converts this infinitesimal generator into a finite orthogonal transformation. In this way, the continuous direction and the discrete feasible update are linked through the Lie algebra.

This observation opens a natural path beyond the Stiefel manifold. For other matrix Lie groups or homogeneous spaces, one may first identify an algebra-valued generator that represents the constrained descent direction, and then use an implicit geometric integrator to map this direction back to the feasible set. The main questions would be to determine which generators lead to descent, how their stationary equations relate to constrained optimality, and which implicit integrators preserve both the geometry and the desired stability properties. Developing this viewpoint systematically could extend \gravidy\ to a much broader family of structured optimization problems.

\section*{Acknowledgements}
I thank Andersen Man Shun Ang for an early assessment of this work and for pointers to the literature on constrained ODEs. I also thank Masoud Ahookhosh for insightful discussions on inner solvers (root-finding for the implicit steps) and for pointers to the literature on Levenberg--Marquardt methods.

\bibliographystyle{plain}
\bibliography{refs}

@article{AlvarezBolteBrahic2004,
  author  = {Alvarez, Felipe and Bolte, J{\'e}r{\^o}me and Brahic, Olivier},
  title   = {Hessian Riemannian Gradient Flows in Convex Programming},
  journal = {SIAM Journal on Control and Optimization},
  volume  = {43},
  number  = {2},
  pages   = {477--501},
  year    = {2004},
  doi     = {10.1137/S0363012902419977}
}

@book{Shahshahani1979,
  author    = {Shahshahani, Siavash},
  title     = {A New Mathematical Framework for the Study of Linkage and Selection},
  series    = {Memoirs of the American Mathematical Society},
  number    = {211},
  publisher = {American Mathematical Society},
  address   = {Providence, RI},
  year      = {1979},
  doi       = {10.1090/memo/0211}
}

@article{rockafellar1976monotone,
  title={Monotone operators and the proximal point algorithm},
  author={Rockafellar, R Tyrrell},
  journal={SIAM journal on control and optimization},
  volume={14},
  number={5},
  pages={877--898},
  year={1976},
  publisher={SIAM}
}

@article{VL2021,
  author  = {Leplat, Valentin and Gillis, Nicolas and F{\'e}votte, C{\'e}dric},
  title   = {Multi-Resolution $\beta$-Divergence {NMF} for Blind Spectral Unmixing},
  journal = {Signal Processing},
  volume  = {193},
  pages   = {108428},
  year    = {2022},
  doi     = {10.1016/j.sigpro.2021.108428}
}

@article{Fevotte_betadiv,
  title={Algorithms for nonnegative matrix factorization with the $\beta$-divergence},
  author={F{\'e}votte, C. and Idier, J.}, 
  journal={Neural computation},
  volume={23},
  number={9},
  pages={2421--2456},
  year={2011},
  publisher={MIT Press}
}

@article{WangWuchen2020,
  author  = {Wang, Yifei and Li, Wuchen},
  title   = {Accelerated Information Gradient Flow},
  journal = {Journal of Scientific Computing},
  volume  = {90},
  number  = {1},
  pages   = {11},
  year    = {2022},
  doi     = {10.1007/s10915-021-01709-3}
}

@inproceedings{Taghvaei2019,
  author    = {Taghvaei, Amirhossein and Mehta, Prashant G.},
  title     = {Accelerated Flow for Probability Distributions},
  booktitle = {Proceedings of the 36th International Conference on Machine Learning},
  series    = {Proceedings of Machine Learning Research},
  volume    = {97},
  pages     = {6076--6085},
  year      = {2019}
}

@InProceedings{pmlr_v97_liu19i,
  title = 	 {Understanding and Accelerating Particle-Based Variational Inference},
  author =       {Liu, Chang and Zhuo, Jingwei and Cheng, Pengyu and Zhang, Ruiyi and Zhu, Jun},
  booktitle = 	 {Proceedings of the 36th International Conference on Machine Learning},
  pages = 	 {4082--4092},
  year = 	 {2019},
  editor = 	 {Chaudhuri, Kamalika and Salakhutdinov, Ruslan},
  volume = 	 {97},
  series = 	 {Proceedings of Machine Learning Research},
  month = 	 {09--15 Jun},
  publisher =    {PMLR},
}

@article{EDELMAN1998,
  title={The Geometry of Algorithms with Orthogonality Constraints},
  author={Edelman, A. and Arias, T.A. and Smith, S.T.},
  journal={SIAM J. MATRIX ANAL. APPL.},
  volume = 	 {20},
  pages={303--353},
  year={1998},
  publisher={SIAM}
}

@article{dahlquist1963special,
  title={A special stability problem for linear multistep methods},
  author={Dahlquist, Germund G},
  journal={BIT Numerical Mathematics},
  volume={3},
  number={1},
  pages={27--43},
  year={1963},
  publisher={Springer}
}

@article{doi:10.1137/0916069,
author = {Byrd, Richard H. and Lu, Peihuang and Nocedal, Jorge and Zhu, Ciyou},
title = {A Limited Memory Algorithm for Bound Constrained Optimization},
journal = {SIAM Journal on Scientific Computing},
volume = {16},
number = {5},
pages = {1190-1208},
year = {1995},
doi = {10.1137/0916069},}

@article{bro1997fast,
  title={A fast non-negativity-constrained least squares algorithm},
  author={Bro, Rasmus and De Jong, Sijmen},
  journal={Journal of Chemometrics: A Journal of the Chemometrics Society},
  volume={11},
  number={5},
  pages={393--401},
  year={1997},
  publisher={Wiley Online Library}
}

@article{10.2307/2153286,
author = {Portugal, Lu\'{\i}s F. and J\'{u}dice, Joaquim J. and Vicente, Lu\'{\i}s N.},
title = {A Comparison of Block Pivoting and Interior-Point Algorithms for Linear Least Squares Problems with Nonnegative Variables},
year = {1994},
issue_date = {Oct. 1994},
publisher = {American Mathematical Society},
address = {USA},
volume = {63},
number = {208},
issn = {0025-5718},
url = {https://doi.org/10.2307/2153286},
doi = {10.2307/2153286},
journal = {Math. Comput.},
month = {oct},
pages = {625–643},
numpages = {19}
}

@book{lawson1995solving,
  title={Solving Least Squares Problems},
  author={Lawson, C.L. and Hanson, R.J.},
  isbn={9780898713565},
  lccn={95035178},
  series={Classics in Applied Mathematics},
  url={https://books.google.ru/books?id=AEwDbHp50FgC},
  year={1995},
  publisher={Society for Industrial and Applied Mathematics}
}

@article{bolte2007lojasiewicz,
  title={The {\L}ojasiewicz inequality for nonsmooth subanalytic functions with applications to subgradient dynamical systems},
  author={Bolte, J{\'e}r{\^o}me and Daniilidis, Aris and Lewis, Adrian},
  journal={SIAM Journal on Optimization},
  volume={17},
  number={4},
  pages={1205--1223},
  year={2007},
  publisher={SIAM}
}

@article{boggs1977algorithm,
  title={An algorithm, based on singular perturbation theory, for ill-conditioned minimization problems},
  author={Boggs, Paul T},
  journal={SIAM Journal on Numerical Analysis},
  volume={14},
  number={5},
  pages={830--843},
  year={1977},
  publisher={SIAM}
}

@phdthesis{zghier1981use,
  title={The use of differential equations in optimization},
  author={Zghier, Abbas K},
  year={1981},
  school={Loughborough University}
}

@article{brown1989some,
  title={Some effective methods for unconstrained optimization based on the solution of systems of ordinary differential equations},
  author={Brown, Andrew A and Bartholomew-Biggs, Michael C},
  journal={Journal of Optimization Theory and Applications},
  volume={62},
  pages={211--224},
  year={1989},
  publisher={Springer}
}

@article{attouch1996dynamical,
  title={A dynamical approach to convex minimization coupling approximation with the steepest descent method},
  author={Attouch, Hedy and Cominetti, Roberto},
  journal={Journal of Differential Equations},
  volume={128},
  number={2},
  pages={519--540},
  year={1996},
  publisher={Elsevier}
}

@article{higham1999trust,
  title={Trust region algorithms and timestep selection},
  author={Higham, Desmond J},
  journal={SIAM Journal on Numerical Analysis},
  volume={37},
  number={1},
  pages={194--210},
  year={1999},
  publisher={SIAM}
}

@article{attouch2004singular,
  title={Singular Riemannian Barrier Methods and Gradient-Projection Dynamical Systems for Constrained Optimization},
  author={Attouch, Hedy and Bolte, J{\'e}r{\^o}me and Redont, Patrick and Teboulle, Marc},
  journal={Optimization},
  volume={53},
  number={5--6},
  pages={435--454},
  year={2004},
  doi={10.1080/02331930412331327184},
  publisher={Taylor \& Francis}
}

@article{attouch2011continuous,
  title={A continuous dynamical Newton-like approach to solving monotone inclusions},
  author={Attouch, Hedy and Svaiter, Benar Fux},
  journal={SIAM Journal on Control and Optimization},
  volume={49},
  number={2},
  pages={574--598},
  year={2011},
  publisher={SIAM}
}

@article{su2014differential,
  title={A differential equation for modeling Nesterov’s accelerated gradient method: theory and insights},
  author={Su, Weijie and Boyd, Stephen and Candes, Emmanuel},
  journal={Advances in neural information processing systems},
  volume={27},
  year={2014}
}

@article{attouch2014dynamical,
  title={A dynamical approach to an inertial forward-backward algorithm for convex minimization},
  author={Attouch, H{\'e}dy and Peypouquet, Juan and Redont, Patrick},
  journal={SIAM Journal on Optimization},
  volume={24},
  number={1},
  pages={232--256},
  year={2014},
  publisher={SIAM}
}

@article{wibisono2016variational,
  title={A variational perspective on accelerated methods in optimization},
  author={Wibisono, Andre and Wilson, Ashia C and Jordan, Michael I},
  journal={proceedings of the National Academy of Sciences},
  volume={113},
  number={47},
  pages={E7351--E7358},
  year={2016},
  publisher={National Acad Sciences}
}

@article{scieur2017integration,
  title={Integration methods and optimization algorithms},
  author={Scieur, Damien and Roulet, Vincent and Bach, Francis and d'Aspremont, Alexandre},
  journal={Advances in Neural Information Processing Systems},
  volume={30},
  year={2017}
}

@inproceedings{franca2018admm,
  title={ADMM and accelerated ADMM as continuous dynamical systems},
  author={Franca, Guilherme and Robinson, Daniel and Vidal, Rene},
  booktitle={International Conference on Machine Learning},
  pages={1559--1567},
  year={2018},
  organization={PMLR}
}

@article{hassan2021proximal,
  title={Proximal gradient flow and Douglas--Rachford splitting dynamics: Global exponential stability via integral quadratic constraints},
  author={Hassan-Moghaddam, Sepideh and Jovanovi{\'c}, Mihailo R},
  journal={Automatica},
  volume={123},
  pages={109311},
  year={2021},
  publisher={Elsevier}
}

@article{shi2021understanding,
  author  = {Shi, Bin and Du, Simon S. and Jordan, Michael I. and Su, Weijie J.},
  title   = {Understanding the Acceleration Phenomenon via High-Resolution Differential Equations},
  journal = {Mathematical Programming},
  volume  = {195},
  pages   = {79--148},
  year    = {2022},
  doi     = {10.1007/s10107-021-01681-8}
}

@inproceedings{suh2022continuous,
  title={Continuous-time analysis of accelerated gradient methods via conservation laws in dilated coordinate systems},
  author={Suh, Jaewook J and Roh, Gyumin and Ryu, Ernest K},
  booktitle={International Conference on Machine Learning},
  pages={20640--20667},
  year={2022},
  organization={PMLR}
}

@article{luo2022differential,
  title={From differential equation solvers to accelerated first-order methods for convex optimization},
  author={Luo, Hao and Chen, Long},
  journal={Mathematical Programming},
  volume={195},
  number={1-2},
  pages={735--781},
  year={2022},
  publisher={Springer}
}

@article{doi:10.1080/10556788.2020.1712602,
author = {Ahookhosh, Masoud and Fleming, Ronan M. T. and Vuong, Phan T.},
title = {Finding zeros of Hölder metrically subregular mappings via globally convergent Levenberg–Marquardt methods},
journal = {Optimization Methods and Software},
volume = {37},
number = {1},
pages = {113--149},
year = {2022},
publisher = {Taylor \& Francis},
doi = {10.1080/10556788.2020.1712602}
}

@book{AbsilMahonySepulchre2008,
    author = {Absil, P.-A. and Mahony, R. and Sepulchre, R.},
    title = {Optimization Algorithms on Matrix Manifolds},
    publisher = {Princeton University Press},
    year = {2008}
}

@book{Boumal2023,
    author = {Boumal, N.},
    title = {An Introduction to Optimization on Smooth Manifolds},
    publisher = {Cambridge University Press},
    year = {2023}
}

@article{WenYin2013,
    author = {Wen, Z. and Yin, W.},
    title = {A feasible method for optimization with orthogonality constraints},
    journal = {Mathematical Programming},
    year = {2013},
    volume = {142},
    pages = {397--434}
}

@book{nesterov2018lectures,
  author    = {Nesterov, Yurii},
  title     = {Lectures on Convex Optimization},
  year      = {2018},
  edition   = {2},
  series    = {Springer Optimization and Its Applications},
  volume    = {137},
  publisher = {Springer},
  address   = {Cham},
  doi       = {10.1007/978-3-319-91578-4},
  isbn      = {978-3-319-91578-4},
  url       = {https://link.springer.com/book/10.1007/978-3-319-91578-4}
}

@article{BirginMartinezRaydan2000,
  author  = {Birgin, Ernesto G. and Mart{\'\i}nez, Jos{\'e} Mario and Raydan, Marcos},
  title   = {Nonmonotone spectral projected gradient methods on convex sets},
  journal = {SIAM Journal on Optimization},
  volume  = {10},
  number  = {4},
  pages   = {1196--1211},
  year    = {2000},
  doi     = {10.1137/S1052623497330963}
}

@article{BeckTeboulle2003,
  author  = {Beck, Amir and Teboulle, Marc},
  title   = {Mirror Descent and Nonlinear Projected Subgradient Methods for Convex Optimization},
  journal = {Operations Research Letters},
  volume  = {31},
  number  = {3},
  pages   = {167--175},
  year    = {2003},
  doi     = {10.1016/S0167-6377(02)00231-6}
}

@book{NemirovskyYudin1983,
  author    = {Nemirovsky, Arkadi S. and Yudin, David B.},
  title     = {Problem Complexity and Method Efficiency in Optimization},
  publisher = {Wiley},
  address   = {New York},
  year      = {1983}
}

@article{KivinenWarmuth1997,
  author  = {Kivinen, Jyrki and Warmuth, Manfred K.},
  title   = {Additive versus Exponentiated Gradient Updates for Linear Prediction},
  journal = {Information and Computation},
  volume  = {132},
  number  = {1},
  pages   = {1--63},
  year    = {1997},
  doi     = {10.1006/inco.1996.2612}
}

@article{lee1999learning,
  title   = {Learning the parts of objects by non-negative matrix factorization},
  author  = {Lee, Daniel D. and Seung, H. Sebastian},
  journal = {Nature},
  volume  = {401},
  number  = {6755},
  pages   = {788--791},
  year    = {1999}
}

@article{AttouchBolteSvaiter2013,
  author    = {Hedy Attouch and J{\'e}r{\^o}me Bolte and Benar F. Svaiter},
  title     = {Convergence of Descent Methods for Semi-Algebraic and Tame Problems: Proximal Algorithms, Forward--Backward Splitting, and Regularized Gauss--Seidel Methods},
  journal   = {Mathematical Programming},
  year      = {2013},
  volume    = {137},
  number    = {1--2},
  pages     = {91--129},
  doi       = {10.1007/s10107-011-0484-9}
}

@article{BolteSabachTeboulle2014,
  author  = {Bolte, J{\'e}r{\^o}me and Sabach, Shoham and Teboulle, Marc},
  title   = {Proximal Alternating Linearized Minimization for Nonconvex and Nonsmooth Problems},
  journal = {Mathematical Programming},
  volume  = {146},
  pages   = {459--494},
  year    = {2014},
  doi     = {10.1007/s10107-013-0701-9}
}

\appendix

\section{One-Dimensional Exponential Map: Exact Solution, Descent, and Limits}
\label{app:positive-flow}

We consider the one-dimensional nonnegative least-squares problem
\[
\min_{x\geq 0}\ \Phi(x):=\frac12(ax-b)^2,
\]
where $a,b\in\mathbb R$. We use the exponential representation $x=e^u$ and start from $x_0>0$. With
\[
q:=a^2\geq0,
\qquad
c:=ab,
\]
the parameter flow is
\begin{equation}
\dot u(t)=c-qe^{u(t)},
\qquad u(0)=\log x_0.
\label{eq:app:ode}
\end{equation}
Equivalently, the primal variable satisfies
\begin{equation}
\dot x(t)=x(t)\bigl(c-qx(t)\bigr)
=-x(t)\Phi'(x(t)).
\label{eq:app:logistic}
\end{equation}
In particular, $x(t)>0$ for every finite $t\geq0$. The boundary $x=0$ can only be reached as a limit.

\paragraph{Exact solution.}
Assume first that $q>0$. If $c\neq0$, solving the scalar equation gives
\begin{equation}
\frac{1}{x(t)}
=
\frac{q}{c}
+
\left(\frac{1}{x_0}-\frac{q}{c}\right)e^{-ct},
\qquad
x(t)
=
\left[
\frac{q}{c}
+
\left(\frac{1}{x_0}-\frac{q}{c}\right)e^{-ct}
\right]^{-1}.
\label{eq:app:solution-nonzero-c}
\end{equation}
If $c=0$, then
\begin{equation}
 x(t)=\frac{x_0}{1+qx_0t}.
\label{eq:app:solution-zero-c}
\end{equation}
These formulas cover the three possible signs of $c$. In particular, no separate local argument is needed to determine the limit.

\paragraph{Descent.}
Since
\[
\Phi'(x)=qx-c,
\]
we obtain directly from \eqref{eq:app:logistic}
\begin{equation}
\frac{d}{dt}\Phi(x(t))
=
\Phi'(x(t))\dot x(t)
=
-x(t)\bigl(qx(t)-c\bigr)^2
\leq0.
\label{eq:app:descent}
\end{equation}
Thus the objective is nonincreasing along the flow. For a finite time $t$, we have $x(t)>0$, so equality in \eqref{eq:app:descent} holds exactly when $qx(t)-c=0$.

\paragraph{Limits and KKT conditions.}
Assume $q>0$, or equivalently $a\neq0$. The limit depends on the sign of $c$.
\begin{itemize}[itemsep=3pt]
\item If $c>0$, then
\[
 x(t)\longrightarrow \frac{c}{q}=\frac{b}{a}>0.
\]
\item If $c=0$, then \eqref{eq:app:solution-zero-c} gives $x(t)\downarrow0$.
\item If $c<0$, then \eqref{eq:app:solution-nonzero-c} gives $x(t)\downarrow0$.
\end{itemize}
Hence, in all cases,
\begin{equation}
 x(t)\longrightarrow x^\star
 =\max\left\{\frac{b}{a},0\right\},
\label{eq:app:nnls-limit}
\end{equation}
which is the unique minimizer of the one-dimensional NNLS problem.

The KKT conditions can be written without introducing an additional sign convention for the multiplier:
\begin{equation}
 x^\star\geq0,
\qquad
 \Phi'(x^\star)=qx^\star-c\geq0,
\qquad
 x^\star\Phi'(x^\star)=0.
\label{eq:app:1d-kkt}
\end{equation}
If $c>0$, then $x^\star=c/q$ and $\Phi'(x^\star)=0$. If $c\leq0$, then $x^\star=0$ and
\[
\Phi'(0)=-c\geq0.
\]
Thus the limit satisfies the KKT conditions in every case.

\paragraph{Behavior of the parameter variable.}
This example also explains why a vanishing primal velocity at the boundary is not, by itself, sufficient to recover the KKT sign condition. Since
\[
\dot u(t)=c-qx(t),
\]
we have the following three behaviors:
\begin{itemize}[itemsep=3pt]
\item if $c>0$, then $u(t)\to\log(c/q)$ and $\dot u(t)\to0$;
\item if $c=0$, then $u(t)\to-\infty$ and $\dot u(t)=-qx(t)\to0$;
\item if $c<0$, then $u(t)\to-\infty$ and $\dot u(t)\to c<0$.
\end{itemize}
In the last case, the primal velocity still satisfies $\dot x(t)\to0$ because of the factor $x(t)$, although the parameter velocity does not vanish. The boundary KKT inequality follows from the sign of the limiting parameter dynamics, not from the equation $\dot x=0$ alone. This is the simple mechanism used in Section~\ref{sec:kkt-from-flow}.

\paragraph{Rates.}
The exact formulas also give the asymptotic rates. If $c>0$, then
\[
 x(t)-x^\star=O(e^{-ct}),
 \qquad
 \Phi(x(t))-\Phi(x^\star)=O(e^{-2ct}).
\]
If $c=0$, then
\[
 x(t)=O(t^{-1}),
 \qquad
 \Phi(x(t))-\Phi(0)=O(t^{-2}).
\]
If $c<0$, then
\[
 x(t)=O(e^{ct}),
 \qquad
 \Phi(x(t))-\Phi(0)=O(e^{ct}),
\]
where $e^{ct}\to0$ because $c<0$.

\paragraph{Degenerate case.}
If $a=0$, then $q=c=0$. The objective is constant,
\[
\Phi(x)=\frac12b^2,
\]
and \eqref{eq:app:ode} gives $\dot u=0$. Thus $x(t)\equiv x_0$, and every feasible point is optimal.

\section{Two Simplex Inner-Solver Variants}
\label{app:simplex-variants}

\subsection{A Hessian-Free Fixed-Point Solver}
\label{app:kl-fixed-point}

This appendix gives a simple fixed-point method for the simplex equation \eqref{eq:kl-implicit-kkt}. It only requires evaluations of $\nabla\Phi$.

For $x\in\operatorname{ri}(\Delta_n)$, define
\begin{equation}
T_k(x)
:=
\frac{x_k\odot\exp(-\eta_k\nabla\Phi(x))}
{\mathbf 1^\top\!\bigl(x_k\odot\exp(-\eta_k\nabla\Phi(x))\bigr)}.
\label{eq:simplex-picard-map}
\end{equation}
A point $x$ satisfies $T_k(x)=x$ if and only if there exists a scalar $\nu$ such that
\[
\log x-\log x_k+\eta_k\nabla\Phi(x)+\nu\mathbf 1=0.
\]
Thus the fixed points of $T_k$ are exactly the solutions of the KL-proximal first-order system.

The direction $d(x)=T_k(x)-x$ is a descent direction for
\[
R_k(x)=\KL(x\|x_k)+\eta_k\Phi(x)
\]
unless $x$ is already a fixed point. Indeed, the definition of $T_k(x)$ gives
\[
\log T_k(x)-\log x_k+\eta_k\nabla\Phi(x)+\nu(x)\mathbf 1=0
\]
for some $\nu(x)$. Since $\mathbf 1^\top d(x)=0$,
\begin{align*}
\langle\nabla R_k(x),d(x)\rangle
&=\langle\log x-\log T_k(x),T_k(x)-x\rangle\\
&=-\sum_{i=1}^n
\bigl(\log x_i-\log T_k(x)_i\bigr)
\bigl(x_i-T_k(x)_i\bigr)
\leq0.
\end{align*}
The inequality is strict when $T_k(x)\neq x$, because the logarithm is strictly increasing.

This observation justifies a relaxed iteration
\[
x^+=x+\lambda\bigl(T_k(x)-x\bigr),
\qquad \lambda\in(0,1],
\]
combined with an Armijo line search on $R_k$. The update remains in the simplex interior for every $\lambda\in(0,1]$.

\begin{algorithm}[H]
\DontPrintSemicolon
\caption{Hessian-free fixed-point solver for the simplex KL step}
\label{alg:KL-fixedpoint}
\KwIn{$x_k\in\operatorname{ri}(\Delta_n)$, $\eta_k>0$, and an inner tolerance $\varepsilon_k$}
$x\leftarrow x_k$\;
\For{$j=0,1,\ldots$}{
  Compute $T\leftarrow T_k(x)$ from \eqref{eq:simplex-picard-map}\;
  $d\leftarrow T-x$\;
  \lIf{$\|d\|\leq\varepsilon_k$}{\textbf{break}}
  Choose $\lambda\in(0,1]$ by Armijo backtracking on $R_k(x+\lambda d)$\;
  $x\leftarrow x+\lambda d$\;
}
\KwOut{$x_{k+1}\leftarrow x$}
\end{algorithm}

The descent calculation above guarantees that a sufficiently small Armijo step exists at every nonstationary inner iterate. A complete global convergence statement would require the usual additional assumptions, for example boundedness away from the boundary of the inner level set or an appropriate compactness argument. We therefore use this method as a practical Hessian-free alternative and do not assign it the local quadratic rate of Newton--KKT.

\subsection{Reduced-Logit Regularized Gauss--Newton Solver}
\label{app:kl-mgn-variant}

Let
\[
v_i=u_i-u_n,
\qquad i=1,\ldots,n-1,
\]
and write
\[
x(v)=\softmax([v;0]).
\]
Define
\[
P=\begin{bmatrix}I_{n-1}&-\mathbf 1_{n-1}\end{bmatrix}.
\]
The reduced form of the KL-proximal first-order system is
\begin{equation}
F_k^{\Delta}(v)
=v-v_k+\eta_kP\nabla\Phi(x(v))=0.
\label{eq:app-reduced-logit-residual}
\end{equation}
The matrix $P$ subtracts the $n$th component from the first $n-1$ components. This removes the softmax gauge and is necessary for equivalence with \eqref{eq:kl-implicit-kkt}.

Let
\[
J_v(v)=J_{\softmax}([v;0])_{:,1:n-1}.
\]
If $\Phi\in C^2$, then
\begin{equation}
DF_k^{\Delta}(v)
=I_{n-1}+\eta_kP\nabla^2\Phi(x(v))J_v(v).
\label{eq:app-reduced-logit-jacobian}
\end{equation}
A regularized Gauss--Newton direction is obtained from
\begin{equation}
\left(DF_k^{\Delta}(v)^\top DF_k^{\Delta}(v)+\lambda I\right)d
=-DF_k^{\Delta}(v)^\top F_k^{\Delta}(v),
\qquad \lambda>0.
\label{eq:app-reduced-logit-lm}
\end{equation}

\begin{algorithm}[H]
\DontPrintSemicolon
\caption{Reduced-logit regularized Gauss--Newton solver}
\label{alg:KL-MGN}
\KwIn{$v_k\in\mathbb R^{n-1}$, $\eta_k>0$, damping $\lambda>0$, and tolerance $\varepsilon_k$}
$v\leftarrow v_k$\;
\For{$j=0,1,\ldots$}{
  $x\leftarrow\softmax([v;0])$\;
  $F\leftarrow v-v_k+\eta_kP\nabla\Phi(x)$\;
  \lIf{$\|F\|\leq\varepsilon_k$}{\textbf{break}}
  $J_v\leftarrow J_{\softmax}([v;0])_{:,1:n-1}$\;
  $J\leftarrow I_{n-1}+\eta_kP\nabla^2\Phi(x)J_v$\;
  Solve $(J^\top J+\lambda I)d=-J^\top F$\;
  Choose $\alpha\in(0,1]$ by backtracking on $\|F_k^{\Delta}(v+\alpha d)\|$\;
  $v\leftarrow v+\alpha d$\;
  Optionally adapt $\lambda$ according to the accepted reduction\;
}
\KwOut{$v_{k+1}\leftarrow v$, $x_{k+1}\leftarrow\softmax([v_{k+1};0])$}
\end{algorithm}

At a root, \eqref{eq:app-reduced-logit-residual} is equivalent to the KL-proximal stationarity system, so this variant computes the same outer step as Newton--KKT. The formula without the matrix $P$ does not have this property.

\section{Fr\'echet Derivative of the Stiefel Cayley Residual}
\label{app:stiefel-frechet}

We derive the derivative used by the Newton--Krylov and dense Newton inner solvers. The derivation is first given for a general smooth objective and then specialized to the quadratic problem used in the experiments.

Let
\[
G(Y)=\nabla\Phi(Y),
\qquad
A(Y)=G(Y)Y^\top-YG(Y)^\top,
\]
and define
\begin{equation}
F_k(Y)
=\left(I+c_kA(Y)\right)Y
-\left(I-c_kA(Y)\right)X_k,
\qquad c_k=\frac{\eta_k}{2}.
\label{eq:app:F}
\end{equation}
Assume that $G$ is Fr\'echet differentiable and write
\[
\mathcal H_Y[H]=DG(Y)[H].
\]
Then
\begin{equation}
DA(Y)[H]
=\mathcal H_Y[H]Y^\top+G(Y)H^\top
-HG(Y)^\top-Y\mathcal H_Y[H]^\top.
\label{eq:app:dA}
\end{equation}
Differentiating \eqref{eq:app:F} gives
\begin{equation}
DF_k(Y)[H]
=\left(I+c_kA(Y)\right)H
+c_kDA(Y)[H](Y+X_k).
\label{eq:app:Frechet}
\end{equation}
This is the exact Jacobian--vector product for the residual in \eqref{eq:ics-fixedpoint}.

\paragraph{Quadratic model.}
For
\begin{equation}
\Phi(X)=\frac12\sum_{j=1}^p
\left\langle Q^{(j)}x^{(j)},x^{(j)}\right\rangle,
\qquad Q^{(j)}=Q^{(j)\top},
\label{eq:app:stiefel-problem}
\end{equation}
we have
\[
G(Y)=\left[Q^{(1)}y^{(1)},\ldots,Q^{(p)}y^{(p)}\right],
\qquad
\mathcal H_Y[H]=G(H).
\]
Hence
\begin{align}
DA(Y)[H]
={}&G(H)Y^\top+G(Y)H^\top
-HG(Y)^\top-YG(H)^\top,
\label{eq:app:dA-quadratic}\\
DF_k(Y)[H]
={}&\left(I+c_kA(Y)\right)H\\
&\quad+c_k\left(G(H)Y^\top+G(Y)H^\top-HG(Y)^\top-YG(H)^\top\right)(Y+X_k).
\label{eq:app:Frechet-quadratic}
\end{align}
The derivative must be evaluated at the current inner point $Y$. Replacing $Y$ and $G(Y)$ in \eqref{eq:app:dA-quadratic} by $X_k$ and $G(X_k)$ gives a frozen quasi-Newton approximation, not the exact derivative.

\paragraph{Vectorized form.}
Let $\operatorname{vec}$ stack the columns of a matrix and let $K_{np}$ be the commutation matrix satisfying
\[
\operatorname{vec}(H^\top)=K_{np}\operatorname{vec}(H).
\]
Define
\[
T=\operatorname{blkdiag}(Q^{(1)},\ldots,Q^{(p)}),
\qquad
X_+=Y+X_k.
\]
Then
\[
\operatorname{vec}(G(H))=T\operatorname{vec}(H).
\]
The matrix representation of \eqref{eq:app:Frechet-quadratic} is
\begin{equation}
\operatorname{vec}(DF_k(Y)[H])
=\mathcal J_k(Y)\operatorname{vec}(H),
\label{eq:app-vectorized-derivative}
\end{equation}
where
\begin{align}
\mathcal J_k(Y)
={}&I_p\otimes\left(I+c_kA(Y)\right)\\
&+c_k\Big[
\left(X_+^\top Y\otimes I_n\right)T
+\left(X_+^\top\otimes G(Y)\right)K_{np}\\
&\hspace{2.7cm}
-\left(X_+^\top G(Y)\otimes I_n\right)
-\left(X_+^\top\otimes Y\right)K_{np}T
\Big].
\label{eq:app-dense-J}
\end{align}
This formula is used only by the dense validation solver. The Newton--Krylov method applies \eqref{eq:app:Frechet} directly and does not form $\mathcal J_k(Y)$.

\paragraph{Local Newton statement.}
Suppose $Y^\star$ is a root of $F_k$ and $DF_k(Y^\star)$ is invertible. Then the standard ambient Newton iteration based on the exact derivative is locally quadratically convergent when the linear systems are solved exactly and $DF_k$ is locally Lipschitz. An inexact Newton--Krylov method has the usual superlinear behavior under an appropriate forcing sequence. These statements apply to the ambient Newton iteration. They do not automatically apply after projecting every correction back to the manifold.

\section{Feasible Predictor and Projected Safeguard for the Cayley Step}
\label{app:f-ics}

This appendix describes two practical devices used by the fast Stiefel implementation.

\subsection{Feasible Cayley Predictor}

Let $\widehat X\in\St(n,p)$ be a predictor, for example $\widehat X=X_k$ or a retracted extrapolation from the two previous outer iterates. Define
\[
Q_{\widehat X}
=\left(I-c_kA(\widehat X)\right)
 \left(I+c_kA(\widehat X)\right)^{-1},
\qquad c_k=\frac{\eta_k}{2},
\]
and set
\[
Y_0=Q_{\widehat X}X_k.
\]
Since $A(\widehat X)$ is skew-symmetric, $Q_{\widehat X}$ is orthogonal and $Y_0\in\St(n,p)$. This predictor is only an initial guess for the nonlinear equation unless $A(\widehat X)=A(Y_0)$.

The inverse of $I+c_kA(\widehat X)$ can be applied without forming an $n\times n$ factorization. Writing
\[
A(\widehat X)=UV^\top-VU^\top,
\qquad U=G(\widehat X),\quad V=\widehat X,
\]
Woodbury reduces the solve to a $2p\times2p$ system. The cost is $O(np^2+p^3)$.

\subsection{Ambient Newton Correction}

The mathematically direct corrector solves
\[
DF_k(Y_j)[H_j]=-F_k(Y_j)
\]
and sets
\[
Y_{j+1}=Y_j+\alpha_jH_j,
\]
where $\alpha_j$ is chosen by backtracking on $\frac12\|F_k(Y)\|_F^2$. The inner iterates need not remain on the Stiefel manifold, but every exact root is feasible. Under the assumptions stated in Appendix~\ref{app:stiefel-frechet}, this is the variant to which the standard Newton convergence theory applies.

\subsection{Projected Safeguard}

A practical alternative is
\[
Y_{j+1}=\operatorname{Proj}_{\St}(Y_j+\alpha_jH_j),
\]
where $\operatorname{Proj}_{\St}$ is computed by a polar factor or a thin QR factorization. This keeps all inner iterates feasible. It is useful when a large ambient correction would otherwise produce a poorly scaled gradient or Hessian evaluation.

The projection changes the Newton map, because the correction $H_j$ obtained from the ambient equation is not necessarily tangent. We therefore use this update only as a safeguard and do not infer quadratic convergence from the ambient Newton theorem. The residual $\|F_k(Y_{j+1})\|_F$ must always be recomputed after the projection.

\begin{algorithm}[H]
\DontPrintSemicolon
\caption{Newton--Krylov corrector with an optional feasibility safeguard}
\label{alg:ICS_res}
\KwIn{$X_k\in\St(n,p)$, $\eta_k>0$, residual tolerance $\varepsilon_k$}
Construct the feasible predictor $Y\leftarrow Y_0$\;
\For{$j=0,1,\ldots$}{
  $R\leftarrow F_k(Y)$\;
  \lIf{$\|R\|_F\leq\varepsilon_k$}{\textbf{break}}
  Solve $DF_k(Y)[H]=-R$ approximately by GMRES\;
  Choose $\alpha\in(0,1]$ by backtracking on $\frac12\|F_k(Y+\alpha H)\|_F^2$\;
  Set $\widetilde Y\leftarrow Y+\alpha H$\;
  Optionally set $\widetilde Y\leftarrow\operatorname{Proj}_{\St}(\widetilde Y)$\;
  $Y\leftarrow\widetilde Y$\;
}
\KwOut{An approximate root $Y$ and its residual $\|F_k(Y)\|_F$}
\end{algorithm}

Before the returned point is tested for outer acceptance, one final polar
projection onto the Stiefel manifold is applied, and both the residual and the
objective are recomputed. The projected point is accepted as the next outer
iterate only after the tests in Section~\ref{sec:stiefel-ics} have been
checked.

\section{Dense Newton Solver for the Cayley Equation}
\label{app:dense-nr}

For moderate values of $np$, the exact Jacobian $\mathcal J_k(Y)$ in \eqref{eq:app-dense-J} can be assembled explicitly. Equivalently, one may build it column by column by applying the exact action \eqref{eq:app:Frechet} to the canonical basis matrices. At inner iteration $j$, solve
\[
\mathcal J_k(Y_j)\operatorname{vec}(H_j)
=-\operatorname{vec}(F_k(Y_j))
\]
by LU, and globalize the correction with a line search on
\[
m_k(Y)=\frac12\|F_k(Y)\|_F^2.
\]

\begin{algorithm}[H]
\DontPrintSemicolon
\caption{Dense Newton solver for the Stiefel Cayley residual}
\label{alg:ICS-NR}
\KwIn{$X_k\in\St(n,p)$, $\eta_k>0$, and residual tolerance $\varepsilon_k$}
Choose an initial guess $Y$, for example $Y=X_k$ or the predictor of Appendix~\ref{app:f-ics}\;
\For{$j=0,1,\ldots$}{
  $R\leftarrow F_k(Y)$\;
  \lIf{$\|R\|_F\leq\varepsilon_k$}{\textbf{break}}
  Assemble the exact Jacobian $\mathcal J_k(Y)$ from \eqref{eq:app-dense-J}\;
  Solve $\mathcal J_k(Y)h=-\operatorname{vec}(R)$ by LU and set $H\leftarrow\operatorname{unvec}(h)$\;
  Choose $\alpha\in(0,1]$ by backtracking on $m_k(Y+\alpha H)$\;
  $Y\leftarrow Y+\alpha H$\;
}
\KwOut{An approximate root $Y$ and its residual $\|F_k(Y)\|_F$}
\end{algorithm}

If the residual is small, the feasibility defect is also small because every
exact root is a Cayley transform of $X_k$. Before outer acceptance, one final
polar projection is applied to remove the remaining feasibility error, and
both the residual and the objective are recomputed. If either outer test
fails, the step is rejected and recomputed with a smaller $\eta_k$.



\section{Complete Proofs for Section~\ref{sec:global-convergence}}
\label{app:proofs-global}

We use the notation of Section~\ref{sec:global-convergence}. In particular, $\mathcal C=\operatorname{ri}(\mathcal X)$, all gradients are taken in the affine hull of $\mathcal X$, and
\[
D_h(x,y)=h(x)-h(y)-\langle\nabla h(y),x-y\rangle.
\]
The second argument of $D_h$ is always an interior point. The first argument may be a boundary point when $h$ has a continuous extension there.

\subsection{Basic Identities and the Exact Implicit Step}

\begin{lemma}[Three-point identity]
\label{lem:three-point}
For $x,y\in\mathcal C$ and every admissible $u\in\mathcal X$,
\begin{equation}
\langle\nabla h(y)-\nabla h(x),u-y\rangle
=
D_h(u,x)-D_h(u,y)-D_h(y,x).
\label{eq:appendix-three-point}
\end{equation}
\end{lemma}

\begin{proof}
Expanding the three Bregman distances gives
\begin{align*}
&D_h(u,x)-D_h(u,y)-D_h(y,x)\\
&=h(u)-h(x)-\langle\nabla h(x),u-x\rangle\\
&\quad-h(u)+h(y)+\langle\nabla h(y),u-y\rangle\\
&\quad-h(y)+h(x)+\langle\nabla h(x),y-x\rangle\\
&=\langle\nabla h(y)-\nabla h(x),u-y\rangle.
\end{align*}
\end{proof}

\begin{lemma}[Optimality equation for the Bregman step]
\label{lem:prox-eq}
Assume that $f$ is convex and differentiable on $\mathcal C$. If the problem
\[
\min_{x\in\mathcal X}
\left\{f(x)+\frac{1}{\eta_k}D_h(x,x_k)\right\}
\]
has a minimizer $x_{k+1}\in\mathcal C$, then this minimizer is unique and satisfies
\begin{equation}
\nabla h(x_{k+1})-\nabla h(x_k)+\eta_k\nabla f(x_{k+1})=0.
\label{eq:appendix-optimality-equation}
\end{equation}
Conversely, any interior point satisfying \eqref{eq:appendix-optimality-equation} is the unique minimizer.
\end{lemma}

\begin{proof}
The objective
\[
\psi_k(x)=f(x)+\frac{1}{\eta_k}D_h(x,x_k)
\]
is strictly convex on $\mathcal C$, because $f$ is convex and $h$ is strictly convex. Its gradient is
\[
\nabla\psi_k(x)
=\nabla f(x)+\frac{1}{\eta_k}\bigl(\nabla h(x)-\nabla h(x_k)\bigr).
\]
An interior minimizer therefore satisfies \eqref{eq:appendix-optimality-equation}. Conversely, for a differentiable convex function, a zero gradient is sufficient for global optimality. Strict convexity gives uniqueness.
\end{proof}

For a nonconvex $f$, equation \eqref{eq:appendix-optimality-equation} is only a first-order condition. It does not identify a global minimizer. This is why the nonconvex analysis in Section~\ref{subsec:nonconvex} explicitly assumes either a global solution of the subproblem or the accepted inexact conditions \eqref{eq:accepted-inexact-conditions}.

\subsection{Convex Descent and the \texorpdfstring{$O(1/S_K)$}{O(1/SK)} Rate}

\begin{proof}[Proof of Proposition~\ref{prop:one-step}]
By convexity of $f$,
\[
f(x_{k+1})-f(x)
\le
\langle\nabla f(x_{k+1}),x_{k+1}-x\rangle.
\]
Using \eqref{eq:global-dual-step},
\begin{align*}
\eta_k\bigl(f(x_{k+1})-f(x)\bigr)
&\le
-\langle\nabla h(x_{k+1})-\nabla h(x_k),x_{k+1}-x\rangle\\
&=
\langle\nabla h(x_{k+1})-\nabla h(x_k),x-x_{k+1}\rangle.
\end{align*}
Apply Lemma~\ref{lem:three-point} with $y=x_{k+1}$, $x=x_k$, and $u=x$. This gives \eqref{eq:one-step}. Setting $x=x_k$ gives
\[
\eta_k\bigl(f(x_{k+1})-f(x_k)\bigr)
\le -D_h(x_{k+1},x_k),
\]
which is \eqref{eq:convex-descent}.

If the comparison point $x$ lies on the boundary, the same inequality follows by taking an interior sequence $x^{(j)}\to x$ and passing to the limit, provided $D_h(x,x_k)$ is finite and $f$ is continuous at $x$.
\end{proof}

\begin{proof}[Proof of Theorem~\ref{thm:global-convex}]
Apply \eqref{eq:one-step} with $x=x^\star$:
\begin{equation}
\eta_k\bigl(f(x_{k+1})-f^\star\bigr)
\le
D_h(x^\star,x_k)-D_h(x^\star,x_{k+1})-D_h(x_{k+1},x_k).
\label{eq:appendix-convex-telescope}
\end{equation}
Summing from $k=0$ to $K$ and dropping the last nonnegative terms gives
\[
\sum_{k=0}^{K}\eta_k\bigl(f(x_{k+1})-f^\star\bigr)
\le
D_h(x^\star,x_0)-D_h(x^\star,x_{K+1})
\le D_h(x^\star,x_0).
\]
This is \eqref{eq:convex-weighted-sum}.

By Proposition~\ref{prop:one-step}, the sequence $f(x_k)$ is nonincreasing. Hence
\[
S_K\bigl(f(x_{K+1})-f^\star\bigr)
\le
\sum_{k=0}^{K}\eta_k\bigl(f(x_{k+1})-f^\star\bigr),
\]
which proves \eqref{eq:convex-last-iterate}.

Since the left-hand side of \eqref{eq:appendix-convex-telescope} is nonnegative, we also obtain
\[
D_h(x^\star,x_{k+1})\le D_h(x^\star,x_k).
\]
Finally, if $S_K\to\infty$, then \eqref{eq:convex-last-iterate} implies $f(x_k)\to f^\star$. By continuity, every cluster point is a minimizer.
\end{proof}

\subsection{Relative Strong Convexity}

\begin{proof}[Proof of Theorem~\ref{thm:linear}]
Apply relative strong convexity \eqref{eq:relative-strong-convexity} with $x=x_{k+1}$ and $y=x^\star$:
\[
f^\star
\ge
f(x_{k+1})
+\langle\nabla f(x_{k+1}),x^\star-x_{k+1}\rangle
+\mu D_h(x^\star,x_{k+1}).
\]
Rearranging gives
\begin{align*}
f(x_{k+1})-f^\star
&\le
\langle\nabla f(x_{k+1}),x_{k+1}-x^\star\rangle
-\mu D_h(x^\star,x_{k+1})\\
&=
\frac{1}{\eta_k}
\Bigl(
D_h(x^\star,x_k)-D_h(x^\star,x_{k+1})-D_h(x_{k+1},x_k)
\Bigr)\\
&\quad-\mu D_h(x^\star,x_{k+1}),
\end{align*}
where the second equality follows from \eqref{eq:global-dual-step} and the three-point identity. Since $f(x_{k+1})-f^\star\ge0$, we can drop this term from the left and obtain
\[
(1+\eta_k\mu)D_h(x^\star,x_{k+1})
+D_h(x_{k+1},x_k)
\le D_h(x^\star,x_k).
\]
This is \eqref{eq:relative-strong-contraction}. Iterating the inequality after dropping $D_h(x_{k+1},x_k)$ gives \eqref{eq:relative-strong-product}.

If $h$ is $\sigma_h$-strongly convex on the region visited by the iterates, then
\[
D_h(x^\star,x_k)\ge\frac{\sigma_h}{2}\|x_k-x^\star\|^2.
\]
The product bound therefore gives linear convergence in norm for a constant stepsize.
\end{proof}

Suppose now that $x^\star\in\mathcal C$, that $f$ is relatively smooth with constant $L_{\rm rel}$, and that
\[
D_h(x,x^\star)\le c_hD_h(x^\star,x)
\]
on the relevant set. Since the gradient of $f$ vanishes in the affine hull at the interior minimizer,
\[
f(x_k)-f^\star\le L_{\rm rel}D_h(x_k,x^\star)
\le c_hL_{\rm rel}D_h(x^\star,x_k).
\]
Combining this with \eqref{eq:relative-strong-product} gives the objective $R$-linear estimate mentioned after Theorem~\ref{thm:linear}.

\subsection{Inexact Convex Steps}

\begin{proof}[Proof of Proposition~\ref{prop:inexact}]
From the definition \eqref{eq:dual-residual-inexact},
\[
\eta_k\nabla f(x_{k+1})
=r_{k+1}-\bigl(\nabla h(x_{k+1})-\nabla h(x_k)\bigr).
\]
By convexity,
\begin{align*}
\eta_k\bigl(f(x_{k+1})-f(x)\bigr)
&\le
\eta_k\langle\nabla f(x_{k+1}),x_{k+1}-x\rangle\\
&=
\langle r_{k+1},x_{k+1}-x\rangle\\
&\quad+
\langle\nabla h(x_{k+1})-\nabla h(x_k),x-x_{k+1}\rangle.
\end{align*}
The three-point identity gives \eqref{eq:inexact-fundamental}. Notice that the residual term is
\[
\langle r_{k+1},x_{k+1}-x\rangle,
\]
not its negative.

Set $x=x^\star$, sum from $k=0$ to $K$, and drop the nonnegative Bregman terms. If $\|x_{k+1}-x^\star\|\le R$, then
\[
\langle r_{k+1},x_{k+1}-x^\star\rangle
\le R\|r_{k+1}\|_*.
\]
This proves \eqref{eq:inexact-convex-sum}. Dividing by $S_K$ gives the best weighted-iterate bound. If the objective values are also nonincreasing, the last objective value is no larger than the weighted average, and the same estimate holds for $x_{K+1}$.
\end{proof}

Under relative strong convexity, repeat the proof of Theorem~\ref{thm:linear} while retaining the residual term. This gives
\begin{align*}
(1+\eta_k\mu)D_h(x^\star,x_{k+1})
+D_h(x_{k+1},x_k)
&\le
D_h(x^\star,x_k)\\
&\quad+\langle r_{k+1},x_{k+1}-x^\star\rangle.
\end{align*}
Dropping the step distance gives \eqref{eq:inexact-linear-recursion}. This recursion also explains why summability alone does not give the exact contraction factor. For example, if the residual term is bounded by $C\gamma^k$ with $\gamma\in(0,1)$ and $\eta_k\equiv\eta$, then a standard linear-recursion argument gives an $R$-linear rate governed by $\max\{(1+\eta\mu)^{-1},\gamma\}$, up to the usual additional factor $k$ when the two numbers coincide.

\subsection{Nonconvex Descent, KL Convergence, and PL Rate}

\begin{proof}[Proof of Proposition~\ref{prop:noncvx-decrease}]
Since $x_{k+1}$ is a global minimizer of \eqref{eq:global-bregman-step}, comparison with $x=x_k$ gives
\[
f(x_{k+1})+\frac{1}{\eta_k}D_h(x_{k+1},x_k)
\le f(x_k).
\]
By the lower Hessian bound in \eqref{eq:nonconvex-h-bounds},
\[
D_h(x_{k+1},x_k)
\ge\frac{\sigma_h}{2}\|x_{k+1}-x_k\|^2.
\]
Using $\eta_k\le\overline\eta$ gives \eqref{eq:nonconvex-sufficient-decrease}. Summing it yields
\[
\frac{\sigma_h}{2\overline\eta}
\sum_{k=0}^{K}\|x_{k+1}-x_k\|^2
\le f(x_0)-f(x_{K+1}),
\]
so the squared step lengths are summable.

The first-order condition gives
\[
\nabla f(x_{k+1})
=-\frac{1}{\eta_k}\bigl(\nabla h(x_{k+1})-\nabla h(x_k)\bigr).
\]
Since $\nabla h$ is $L_h$-Lipschitz on $\mathcal K$,
\[
\|\nabla f(x_{k+1})\|
\le\frac{L_h}{\eta_k}\|x_{k+1}-x_k\|
\le\frac{L_h}{\underline\eta}\|x_{k+1}-x_k\|.
\]
This proves \eqref{eq:nonconvex-relative-error} and \eqref{eq:nonconvex-gradient-zero}.

For the square-summability estimate, combine the Lipschitz bound on $\nabla h$ with the lower Hessian bound:
\begin{align*}
\eta_k\|\nabla f(x_{k+1})\|^2
&\le
\frac{L_h^2}{\eta_k}
\|x_{k+1}-x_k\|^2\\
&\le
\frac{2L_h^2}{\sigma_h\eta_k}
D_h(x_{k+1},x_k)\\
&\le
\frac{2L_h^2}{\sigma_h}
\bigl(f(x_k)-f(x_{k+1})\bigr).
\end{align*}
Summation gives \eqref{eq:gradient-square-sum}.
\end{proof}

\begin{proof}[Proof of Theorem~\ref{thm:KL-conv}]
Set
\[
d_k:=\|x_{k+1}-x_k\|,
\qquad
f_k:=f(x_k).
\]
Proposition~\ref{prop:noncvx-decrease} gives two constants
\[
a:=\frac{\sigma_h}{2\overline\eta}>0,
\qquad
b:=\frac{L_h}{\underline\eta}>0,
\]
such that
\begin{equation}
 f_k-f_{k+1}\ge a d_k^2,
 \qquad
 \|\nabla f(x_{k+1})\|\le b d_k.
\label{eq:appendix-H1-H2}
\end{equation}
The objective values decrease and are bounded below, so $f_k\to f_\infty$.

The sequence is contained in the compact set $\mathcal K$, hence it has a nonempty compact set of cluster points. Since $d_k\to0$ and $\nabla f(x_k)\to0$, every cluster point is critical and has objective value $f_\infty$. The KL property can therefore be uniformized on this compact cluster set: there exist an index $k_0$, a neighborhood of the cluster set, and a concave $C^1$ function $\varphi$ with $\varphi'>0$ such that
\begin{equation}
 \varphi'(f_k-f_\infty)\|\nabla f(x_k)\|\ge1,
 \qquad k\ge k_0.
\label{eq:uniform-KL}
\end{equation}
This is the standard uniformized KL inequality used in descent analyses; see, for example, \cite{AttouchBolteSvaiter2013,BolteSabachTeboulle2014}.

Let $s_k=f_k-f_\infty$. By concavity of $\varphi$,
\begin{align*}
\varphi(s_k)-\varphi(s_{k+1})
&\ge \varphi'(s_k)(s_k-s_{k+1})\\
&\ge a\varphi'(s_k)d_k^2.
\end{align*}
If $d_{k-1}=0$, the second estimate in \eqref{eq:appendix-H1-H2} gives $\nabla f(x_k)=0$. The KL inequality then implies $s_k=0$, and the convergence conclusion is immediate. We may therefore consider the case $d_{k-1}>0$. Using \eqref{eq:uniform-KL} and the second estimate in \eqref{eq:appendix-H1-H2} at index $k-1$ gives
\[
\varphi'(s_k)\ge\frac{1}{\|\nabla f(x_k)\|}
\ge\frac{1}{bd_{k-1}}.
\]
Hence
\begin{equation}
 d_k^2
 \le \frac{b}{a}d_{k-1}
 \bigl(\varphi(s_k)-\varphi(s_{k+1})\bigr).
\label{eq:finite-length-key}
\end{equation}
Using $2\sqrt{uv}\le u+v$ in \eqref{eq:finite-length-key},
\[
2d_k
\le d_{k-1}
+\frac{b}{a}\bigl(\varphi(s_k)-\varphi(s_{k+1})\bigr).
\]
Summing this inequality from $k=k_0$ to $m$ gives
\[
\sum_{k=k_0}^{m}d_k+d_m
\le d_{k_0-1}+\frac{b}{a}\varphi(s_{k_0}).
\]
Letting $m\to\infty$ proves $\sum_kd_k<\infty$. Thus $(x_k)$ is Cauchy and converges to one critical point $x^\star$.

For the rates, assume that the KL function can be chosen as
\[
\varphi(s)=c s^{1-\theta},
\qquad \theta\in[0,1).
\]
The KL inequality and \eqref{eq:appendix-H1-H2} imply, for all large $k$,
\[
 d_k\ge c_1 s_{k+1}^{\theta}
\]
for some $c_1>0$. Therefore
\begin{equation}
 s_k-s_{k+1}
 \ge c_2s_{k+1}^{2\theta}
\label{eq:KL-rate-recursion}
\end{equation}
for some $c_2>0$.

If $\theta=0$, \eqref{eq:KL-rate-recursion} gives a fixed positive decrease whenever $s_{k+1}>0$, so the sequence terminates in finitely many steps. If $0<\theta\le1/2$, then $s_{k+1}^{2\theta}\ge s_{k+1}$ once $s_{k+1}\le1$, and \eqref{eq:KL-rate-recursion} gives a linear contraction. If $1/2<\theta<1$, the standard discrete comparison for \eqref{eq:KL-rate-recursion} yields
\[
s_k=O\bigl(k^{-1/(2\theta-1)}\bigr).
\]
This proves the stated rate trichotomy.
\end{proof}

\paragraph{Accepted inexact steps.}
Suppose an approximate sequence satisfies \eqref{eq:accepted-inexact-conditions}. Then the proof above applies verbatim with the constants $a$ and $b$ appearing in that condition. To connect the second condition with the inner residual, write
\[
\nabla f(x_{k+1})
=-\frac{1}{\eta_k}\bigl(\nabla h(x_{k+1})-\nabla h(x_k)\bigr)+e_{k+1}.
\]
If $\|e_{k+1}\|\le\gamma\|x_{k+1}-x_k\|$, then
\[
\|\nabla f(x_{k+1})\|
\le
\left(\frac{L_h}{\underline\eta}+\gamma\right)
\|x_{k+1}-x_k\|.
\]
Thus a relative inner stopping rule, together with sufficient decrease, gives exactly the two estimates used in the KL proof.

\begin{proof}[Proof of Theorem~\ref{thm:PL-linear}]
By the optimality equation and the two Hessian bounds on $h$,
\begin{align*}
\|\nabla f(x_{k+1})\|^2
&\le
\frac{L_h^2}{\eta^2}
\|x_{k+1}-x_k\|^2\\
&\le
\frac{2L_h^2}{\sigma_h\eta^2}
D_h(x_{k+1},x_k).
\end{align*}
The sufficient-decrease inequality gives
\[
D_h(x_{k+1},x_k)
\le\eta\bigl(f(x_k)-f(x_{k+1})\bigr).
\]
Consequently,
\[
\|\nabla f(x_{k+1})\|^2
\le
\frac{2L_h^2}{\sigma_h\eta}
\bigl(f(x_k)-f(x_{k+1})\bigr).
\]
Apply \eqref{eq:PL-inequality} at $x_{k+1}$:
\[
f(x_{k+1})-f^\star
\le
\frac{L_h^2}{\mu_{\rm PL}\sigma_h\eta}
\bigl(f(x_k)-f(x_{k+1})\bigr).
\]
Rearranging gives \eqref{eq:PL-rate}.
\end{proof}

\subsection{A--Stability for a Fixed-Metric Linear--Quadratic Model}

We finish with the linear--quadratic model used to motivate the stability discussion in Section~\ref{sec:stability}. This calculation is exact and concerns a fixed metric. It does not replace the nonlinear convergence results proved above.

\begin{proposition}[Backward Euler for a fixed-metric quadratic]
\label{prop:linear-quadratic-stability}
Let
\[
f(x)=\frac12x^\top Qx-c^\top x,
\qquad
Q\succeq0,
\qquad
G\succ0,
\]
and assume that the solution set
\[
\mathcal X^\star:=\{x\in\mathbb R^n:Qx=c\}
\]
is nonempty. For a fixed stepsize $\eta>0$, consider
\begin{equation}
G\frac{x_{k+1}-x_k}{\eta}+Qx_{k+1}-c=0.
\label{eq:F6-backward-euler}
\end{equation}
Then the following statements hold.

\begin{enumerate}[label=\textnormal{(\roman*)},itemsep=4pt]
\item For every $x^\star\in\mathcal X^\star$,
\begin{equation}
x_{k+1}-x^\star
=(G+\eta Q)^{-1}G(x_k-x^\star).
\label{eq:F6-error-iteration}
\end{equation}
The eigenvalues of the iteration matrix belong to $(0,1]$. They all belong to $(0,1)$ when $Q\succ0$.

\item The sequence $(x_k)$ converges to a point $x_\infty\in\mathcal X^\star$. If $Q\succ0$, then $\mathcal X^\star=\{x^\star\}$ and
\begin{equation}
\|x_{k+1}-x^\star\|_G
\leq
\frac{1}{1+\eta\lambda_{\min}(S)}
\|x_k-x^\star\|_G,
\qquad
S:=G^{-1/2}QG^{-1/2}.
\label{eq:F6-G-contraction}
\end{equation}

\item The objective satisfies the exact descent identity
\begin{equation}
f(x_k)-f(x_{k+1})
=
\frac{1}{\eta}\|x_{k+1}-x_k\|_G^2
+
\frac12\|x_{k+1}-x_k\|_Q^2
\geq0.
\label{eq:F6-exact-descent}
\end{equation}
\end{enumerate}

If $Q\neq0$ and $\lambda_+(S)$ denotes the smallest positive eigenvalue of $S$, then
\begin{equation}
f(x_k)-f^\star
\leq
\left(\frac{1}{1+\eta\lambda_+(S)}\right)^{2k}
\bigl(f(x_0)-f^\star\bigr).
\label{eq:F6-objective-rate}
\end{equation}
\end{proposition}

\begin{proof}
Let $x^\star\in\mathcal X^\star$. Since $Qx^\star=c$, subtracting this identity from \eqref{eq:F6-backward-euler} gives \eqref{eq:F6-error-iteration}.

Set
\[
y_k:=G^{1/2}(x_k-x^\star),
\qquad
S:=G^{-1/2}QG^{-1/2}\succeq0.
\]
Then \eqref{eq:F6-error-iteration} becomes
\begin{equation}
y_{k+1}=(I+\eta S)^{-1}y_k.
\label{eq:F6-transformed-iteration}
\end{equation}
The matrix $S$ is symmetric positive semidefinite. If $\lambda_i\geq0$ is one of its eigenvalues, the corresponding amplification factor is
\[
\frac{1}{1+\eta\lambda_i}\in(0,1].
\]
This proves the spectral statement in part (i). If $Q\succ0$, then $S\succ0$, and all amplification factors are strictly smaller than one. Taking the largest one gives \eqref{eq:F6-G-contraction}.

If $Q$ is singular, the components of $y_k$ in $\ker S$ remain unchanged, while all components associated with positive eigenvalues converge to zero. Therefore $y_k$ converges to a point $y_\infty\in\ker S$. Equivalently, $x_k$ converges to $x_\infty=x^\star+G^{-1/2}y_\infty$, which satisfies
\[
Qx_\infty=c.
\]
This proves part (ii).

Let $d_k:=x_{k+1}-x_k$. Since $f$ is quadratic,
\[
f(x_k)
=
f(x_{k+1})
+\langle\nabla f(x_{k+1}),x_k-x_{k+1}\rangle
+\frac12d_k^\top Qd_k.
\]
Equation \eqref{eq:F6-backward-euler} gives
\[
\nabla f(x_{k+1})
=Qx_{k+1}-c
=-\frac{1}{\eta}Gd_k.
\]
Substitution yields \eqref{eq:F6-exact-descent}.

Finally, for any $x^\star\in\mathcal X^\star$,
\[
f(x_k)-f^\star
=
\frac12(x_k-x^\star)^\top Q(x_k-x^\star)
=
\frac12y_k^\top Sy_k.
\]
In an orthonormal eigenbasis of $S$, every component associated with a positive eigenvalue is multiplied at step $k$ by $(1+\eta\lambda_i)^{-k}$. Bounding these factors by the one associated with $\lambda_+(S)$ gives \eqref{eq:F6-objective-rate}.
\end{proof}

\paragraph{Relation with A--stability.}
The continuous fixed-metric flow is
\[
G\dot x+Qx-c=0.
\]
In the coordinates $y=G^{1/2}(x-x^\star)$, it becomes
\[
\dot y=-Sy.
\]
Each eigenmode therefore satisfies the scalar equation $\dot z=-\lambda_i z$. Backward Euler has amplification factor
\[
R(-\eta\lambda_i)=\frac{1}{1+\eta\lambda_i},
\]
which has modulus at most one for every $\eta>0$. Proposition~\ref{prop:linear-quadratic-stability} is the matrix form of this classical A--stability calculation. The conclusion is limited to the fixed linear model. The nonlinear descent and convergence results of this paper follow instead from the Bregman arguments in Section~\ref{sec:global-convergence}.

\end{document}